\documentclass[11pt]{amsart}

\usepackage{amsmath}
\usepackage{amssymb}
\usepackage{bbm}
\usepackage{pdfsync}
\usepackage{esint} %integrals
\usepackage[pagebackref]{hyperref}
\usepackage[utf8]{inputenc}
\usepackage[T1]{fontenc}
\usepackage{mathtools}

\usepackage{etoolbox}
\patchcmd{\section}{\scshape}{\scshape\bfseries}{}{}
\renewcommand*{\backref}[1]{}
\renewcommand*{\backrefalt}[4]{%
	\ifcase #1 (Not cited.)%
	\or        (Cited on page~#2)%
	\else      (Cited on pages~#2)%
	\fi}

\usepackage[capitalise]{cleveref} %after hyperref. Use \cref{label}
\crefname{equation}{}{} %this removes eq. in eq. (1.1) for example
\crefname{claim}{Claim}{Claims} %this produces Claims X and Y when \cref{claim1,claim2}
\crefname{page}{p.}{pp.}
\crefname{rem}{Remark}{Remarks} %this produces Remarks X and Y when \cref{remark1,remark2}
\crefname{coro}{Corollary}{Corollaries} %this produces Corollaries X and Y when \cref{coro1,coro2}
\crefname{enumi}{item}{items}
\crefname{figure}{fig.}{figs.}  % For \cref (lowercase, short)
\Crefname{figure}{Figure}{Figures} % For \Cref (capitalized, full word)

\usepackage{enumitem}
\usepackage{aligned-overset}

\hypersetup{ 
	linktoc=page, %linktoc=all (to have link in title and page number in the table of content
	colorlinks=true,
	linkcolor={red!90!black},
	citecolor={blue!90!black},
	urlcolor=magenta,%{blue!90!black}, %link (defaul:magenta)
	filecolor=teal,
}

\let\oldtocsection=\tocsection
\let\oldtocsubsection=\tocsubsection
\let\oldtocsubsubsection=\tocsubsubsection
\renewcommand{\tocsection}[2]{\hspace{0em}\oldtocsection{#1}{#2}}
\renewcommand{\tocsubsection}[2]{\hspace{1em}\oldtocsubsection{#1}{#2}}
\renewcommand{\tocsubsubsection}[2]{\hspace{2em}\oldtocsubsubsection{#1}{#2}}

\makeatletter
\newcommand\@dotsep{4.5}
\def\@tocline#1#2#3#4#5#6#7{\relax
	\ifnum #1>\c@tocdepth % then omit
	\else
	\par \addpenalty\@secpenalty\addvspace{#2}%
	\begingroup \hyphenpenalty\@M
	\@ifempty{#4}{%
		\@tempdima\csname r@tocindent\number#1\endcsname\relax
	}{%
		\@tempdima#4\relax
	}%
	\parindent\z@ \leftskip#3\relax
	\advance\leftskip\@tempdima\relax
	\rightskip\@pnumwidth plus1em \parfillskip-\@pnumwidth
	#5\leavevmode\hskip-\@tempdima #6\relax
	\leaders\hbox{$\m@th
		\mkern \@dotsep mu\hbox{.}\mkern \@dotsep mu$}\hfill
	\hbox to\@pnumwidth{\@tocpagenum{#7}}\par
	\nobreak
	\endgroup
	\fi}
\makeatother

\newcommand{\R}{{\mathbb R}}       % Field of real numbers

\newcommand{\Z}{{\mathbb Z}}       % Ring of integer numbers
\newcommand{\DD}{{\mathcal D}}

\newcommand{\HH}{{\mathcal H}}

\newcommand{\WW}{{\mathcal W}}

\newcommand{\EE}{{\mathcal E}}

\newcommand{\OO}{{\mathcal O}}

\newcommand{\diam}{{\rm diam}}
\newcommand{\dist}{{\rm dist}}

\newcommand{\rf}[1]{{(\ref{#1})}}

\newcommand{\supp}{\operatorname{supp}}

\newcommand{\vv}{{\vspace{2mm}}}

\usepackage{dsfont}
\newcommand{\characteristic}{\mathbf{1}} %characteristic function \chi or \mathbf{1} or \mathds{1}

\newcommand{\good}{{\mathsf{Good}}}

\newcommand{\HD}{{\mathsf{HD}}}

\newcommand{\TS}{{\mathsf{TS}}} %touching set (only in this article)

\newcommand{\divv}{\operatorname{div}}
\newcommand{\Capacity}{\operatorname{Cap}}

\newcommand{\loc}{\operatorname{loc}} %locally

\newcommand{\rom}[1]{%
	\textup{\uppercase\expandafter{\romannumeral#1}}%
}

\def\Xint#1{\mathchoice
	{\XXint\displaystyle\textstyle{#1}}%
	{\XXint\textstyle\scriptstyle{#1}}%
	{\XXint\scriptstyle\scriptscriptstyle{#1}}%
	{\XXint\scriptscriptstyle\scriptscriptstyle{#1}}%
	\!\int}
\def\XXint#1#2#3{{\setbox0=\hbox{$#1{#2#3}{\int}$ }
		\vcenter{\hbox{$#2#3$ }}\kern-.58\wd0}}

\def\avint{\;\Xint-}

\makeatletter
\newcommand{\doublewidetilde}[1]{{%
		\mathpalette\double@widetilde{#1}%
}}
\newcommand{\double@widetilde}[2]{%
	\sbox\z@{$\m@th#1\widetilde{#2}$}%
	\ht\z@=.8\ht\z@
	\widetilde{\box\z@}%
}
\makeatother

\usepackage{todonotes}

\usepackage{pgf,tikz} %tikz
\usetikzlibrary{lindenmayersystems}

\usepackage{graphicx}
\usepackage{caption}
\usepackage{subcaption}

\newtheorem{theorem}{Theorem}[section]
\newtheorem{lemma}[theorem]{Lemma}
\newtheorem{mlemma}[theorem]{{M}ain Lemma}

\newtheorem{coro}[theorem]{Corollary}

\newtheorem{claim}[theorem]{Claim}
\newtheorem*{claim*}{Claim}

\newtheorem*{theorem*}{Theorem}

\theoremstyle{definition}
\newtheorem{definition}[theorem]{Definition}
\newtheorem*{notation}{Notation}
\newtheorem{example}[theorem]{Example}

\newtheorem{rem}[theorem]{Remark} %use this or the following two lines
\numberwithin{equation}{section}

\pgfdeclarelindenmayersystem{Koch}
{\rule {F -> F+F--F+F}}

\definecolor{taronja}{rgb}{0.9,0.5,0.05}
\usepackage{xcolor}

\usepackage{stackengine}
\newcommand\ithposition{\mathop{\stackMath\stackinset{c}{0pt}{c}{-0.6ex}{{\scriptscriptstyle\smile}}{{\scriptscriptstyle i}}}}
\begin{document}
%[short title]{long title}
\title[The dimension of planar Lipschitz elliptic measures in Reifenberg flat domains]{The dimension of planar elliptic measures arising from Lipschitz matrices in Reifenberg flat domains}

\author[I. Guillén-Mola]{Ignasi Guillén-Mola}
%info Ignasi
\address{Ignasi Guillén-Mola,
	Departament de Matem\`atiques, Universitat Aut\`onoma de Barcelona.
}
\email{ignasi.guillen@uab.cat}

\author[M. Prats]{Martí Prats}
%info Martí
\address{Martí Prats,
	Departament de Matem\`atiques, Universitat Aut\`onoma de Barcelona, and Centre de Recerca Matemàtica, Barcelona, Catalonia.}
\email{marti.prats@uab.cat}

\author[X. Tolsa]{Xavier Tolsa}
%info Xavi
\address{Xavier Tolsa, ICREA, Barcelona, Departament de Matemàtiques, Universitat Autònoma de Barcelona, and Centre de Recerca Matemàtica, Barcelona, Catalonia.}
\email{xavier.tolsa@uab.cat}

\date{\today}



\begin{abstract}	
	In this paper we show that, given a planar Reifenberg flat domain with small constant and a divergence form operator associated to a real (not necessarily symmetric) uniformly elliptic matrix with Lipschitz coefficients, the Hausdorff dimension of its elliptic measure is at most $1$. More precisely, we prove that there exists a subset of the boundary with full elliptic measure and with $\sigma$-finite one-dimensional Hausdorff measure. For Reifenberg flat domains, this result extends a previous work of Thomas H. Wolff for the harmonic measure.
\end{abstract}

\maketitle

{
	\tableofcontents
}

% ********************************************************************************
% ********************************************************************************

\section{Introduction and Main results}

We study the dimension of planar elliptic measures in Reifenberg flat domains with small constant, assuming also Lipschitz continuity of the coefficients of the matrix. In fact, that regularity is only needed near the boundary.

Let $L_A u\coloneqq -\divv (A\nabla u)$ be a second order operator in divergence form associated to the $(n+1) \times (n+1)$ real matrix $A(\cdot)=(a_{ij} (\cdot))_{1\leq i,j \leq n+1}$ such that there exists $\lambda \geq 1$ with
\begin{subequations}
\begin{align}
	\lambda^{-1}|\xi|^{2} \leq\langle A(x) \xi, \xi\rangle& \text{ for all } \xi \in \R^{n+1} \text{ and a.e.\ } x\in \R^{n+1},\label{elliptic1}\\
	\langle A(x) \xi, \eta\rangle \leq \lambda |\xi||\eta|& \text{ for all } \xi,\eta \in \R^{n+1} \text{ and a.e.\ } x\in \R^{n+1}.\label{elliptic2}
\end{align}
\end{subequations}
We say that a measurable matrix is uniformly elliptic with ellipticity constant $\lambda \geq 1$ if \rf{elliptic1} and \rf{elliptic2} are satisfied. The uniform ellipticity condition implies that the matrix has bounded coefficients.

Let $\Omega \subset \R^2$ be an open set and fix a point $p\in \Omega$. Consider the operator
$$
\begin{array}{cccc}
T:	& C(\partial \Omega) & \to & \mathbb{R}\\
	& f &\mapsto & u_f^{L_A} (p),
\end{array}
$$
where $u_f^{L_A}$ is the $L_A$-harmonic extension of $f$ in $\Omega$ via Perron's method. By the maximum principle for $L_A$-harmonic functions with uniformly elliptic matrices (see \cite[p.~46]{Gilbarg2001} for example), the operator $T$ is linear, bounded and positive. For bounded open sets $\Omega\subset \R^{n+1}$ with $n\geq 2$ we do the same construction.

The elliptic ($L_A$-harmonic) measure in $\Omega$ with pole $p\in\Omega$ is the unique Radon probability measure $\omega_{\Omega, A}^p$ (by the Riesz representation theorem) such that
$$
u_f^{L_A} (p) = \int_{\partial \Omega} f (\xi)\, d\omega_{\Omega, A}^p (\xi) \text{ for every } f\in C(\partial \Omega) .
$$
Via the Perron method construction, although the characteristic function $\characteristic_E$ of a set $E\subset \partial \Omega$ is not continuous in $\partial \Omega$, the function $x\mapsto \omega_{\Omega, A}^x (E)$ is $L_A$-harmonic in $\Omega$. For a detailed construction of elliptic/harmonic measures, see \cite[Section 11]{Heinonen2006}. If the set $\Omega$, the matrix $A$, or the pole $p\in \Omega$ is clear from the context, we will omit it in $\omega_{\Omega, A}^p$.

In this work, we study the Hausdorff dimension of elliptic measures arising from uniformly elliptic matrices. This is defined as
$$
\dim_\HH \omega_{\Omega, A}^p \coloneqq \inf \{\dim_\HH F : \omega_{\Omega, A}^p (F^c) = 0\}.
$$
Naturally $\dim_\HH \omega \leq \dim_\HH \partial \Omega$. The Hausdorff dimension of the harmonic measure (i.e., with $\omega = \omega_{Id}$) has been studied by several authors, both in the plane and in higher dimensions, showing a different behavior in each case.

In the plane, Carleson proved $\dim_\HH \omega \leq 1$ for ``snowflake type'' sets and $\dim_\HH \omega <1$ for some self similar ``Cantor type'' sets, both results in \cite{Carleson1985}. More precise results were obtained for simply connected domains by Makarov in \cite{Makarov1985} by showing $\dim_\HH \omega =1$. The upper bound $\dim_\HH \omega \leq 1$ with no assumptions on the domain was shown by Jones and Wolff in \cite{Jones1988}, and later it was improved by Wolff in \cite{Wolff1993} by proving that there is a subset $F\subset \partial \Omega$ with full harmonic measure $\omega (F)=1$ and $\sigma$-finite length.

In the same direction for higher dimensions, in $\R^{n+1}$ with $n\geq 2$, Bourgain proved in \cite{Bourgain1987} that there exists a dimensional constant $b_n >0$ such that $\dim_\HH \omega \leq n+1-b_n$. By the results on the previous paragraph, we can take $b_1 = 1$ and this choice is optimal. For $n\geq 2$, Wolff constructed in \cite{Wolff1995} a domain $\Omega_n \subset \R^{n+1}$ such that $\dim_\HH \omega_{\Omega_n} > n$, and in particular the Bourgain constant can't equal $1$, i.e., $b_n < 1$. In a recent work \cite{Badger2024}, Badger and Genshaw refined the proof in \cite{Bourgain1987} to find estimates on the Bourgain constant $b_n \in (0,1)$.

From the results in the previous paragraphs we have that $\dim_\HH \omega < n+1$ (with some precise gap), but possibly $\dim_\HH \omega = \dim_\HH \partial \Omega$. In fact, the situation $\dim_\HH \omega < \dim_\HH \partial \Omega$ holds for many non-trivial domains. This phenomenon is frequently called the ``dimension drop'' for harmonic measure. The dimension drop is closely related to the results mentioned above. Indeed, similar techniques to the ones in \cite{Bourgain1987,Jones1988} are used in some of the following articles. The first work in this direction is due to Kaufman and Wu \cite{Kaufman-Wu-1985} for the planar $\log 4/\log 3$-dimensional Koch snowflake. Subsequently, Carleson \cite{Carleson1985} extended this result to self-similar Cantor sets in the plane. For this type of domains, see Makarov and Volberg \cite{Makarov1986}, Batakis \cite{Batakis1996}, and also \cite{Volberg1992,Volberg1993}. Jones and Wolff showed that the dimension drop happens for some uniform and disconnected planar domains, see Theorem 2.1 in \cite[Section \rom{9}.2]{Garnett2005}. For IFS domains (iterated functions systems), see Urba\'{n}ski and Zdunik \cite{Urbanski2002}, and Batakis and Zdunik \cite{Batakis2015}.

For general AD-regular domains of fractional dimension, one may ask if the dimension drop happens. It holds on uniformly ``non-flat'' AD-regular domains with codimension smaller than $1$, as shown by Azzam in \cite{Azzam2020}. Recently, for $n\geq 1$ and $s\in [n-\frac{1}{2},n)$, the third author showed in \cite{Tolsa2024} that the dimension drop occurs for higher codimensional $s$-AD-regular subsets of $C^1$ $n$-dimensional manifolds in $\R^{n+1}$. In contrast, David, Jeznach and Julia proved in \cite{David2023} that this last result may fail for $s$ close enough to $n-1$.

The situation is different for uniformly elliptic matrices. In \cite{Sweezy1992} in the planar case and in \cite{Sweezy1994} in higher dimensions, for every $\varepsilon > 0$, Sweezy constructed a domain $\Omega\subset \R^{n+1}$ and an elliptic operator in divergence form $L_A$ whose associated elliptic measure $\omega_{\Omega,A}$ satisfies $\dim_\HH \omega_{\Omega,A}\geq n+1-\varepsilon$. However, as $\varepsilon$ becomes smaller, the ellipticity constant of the resulting matrix increases. Such planar domains and elliptic measures are constructed by the push forward under quasiconformal mappings, and the higher dimensional analog is deduced from the planar case.

Using a new approach in the planar case, David and Mayboroda in \cite{David2021} constructed an elliptic operator $L=-\divv a\nabla$, where $a$ is a uniformly elliptic and continous scalar function on the complementary of the four corner Cantor set of dimension $1$, whose elliptic measure $\omega_{aId}$ is proportional to the one-dimensional Hausdorff measure on the Cantor set. Operators of this form are the so-called ``good elliptic operators''. Following the same strategy, Perstneva in \cite{Perstneva2025} constructed ``good elliptic operators'' on the complementary of the planar $d$-dimensional Koch-type snowflake with $1<d<\log4 /\log3$, whose elliptic measure is proportional to the $d$-dimension Hausdorff measure on the Koch snowflake.

After Sweezy's results, it is natural to ask which conditions on the matrix imply the analogous result of \cite{Jones1988,Wolff1993} for the elliptic case. In this paper we study the metric properties of the elliptic measure when assuming regularity conditions on the matrix $A$. For other results in this line, see for example \cite{Prat2021}, about the rectifiability of the elliptic measure for Hölder matrices.

Suppose from now on that the domain is $(\delta,r_0)$-Reifenberg flat and the coefficients of the matrix are also Lipschitz. Roughly speaking, a domain $\Omega\subset\R^2$ is $(\delta,r_0)$-Reifenberg flat if for every ball $B$ centered at $\partial\Omega$ and with radius smaller than $r_0$, the $\delta$-neighborhood of a line through the center of $B$ contains $B\cap\partial\Omega$. See \cref{def:Reifenberg flat} for the details. In this situation we show that the dimension of the elliptic measure in Reifenberg flat domains with small constant is at most $1$. More precisely, for this type of domains we obtain the analogous result of \cite[Theorem 1]{Wolff1993} in the following theorem.

\begin{theorem}\label{dimension elliptic measure}
	Let $\Omega \subset \R^2$ be a $(\delta, r_0)$-Reifenberg flat domain, $p\in \Omega$, and $A$ be a real uniformly elliptic (not necessarily symmetric) matrix with ellipticity constant $\lambda$. Suppose also that $A$ is Lipschitz. Then there exists $\delta_0 = \delta_0 (\lambda) >0$ such that if $0<\delta \leq \delta_0$ then there is a set $F\subset \partial \Omega$ satisfying $\omega_{\Omega,A} (F)=1$ and with $\sigma$-finite one-dimensional Hausdorff measure. In particular $\dim_\HH \omega_{\Omega,A} \leq 1$.
\end{theorem}

Despite we are requiring to work with $\delta$-Reifenberg flat domains with small enough constant $\delta$, such sets can be constructed with Hausdorff dimension strictly larger than 1. For example, a suitable variant of the Koch snowflake can be constructed to be $\delta$-Reifenberg flat.

It is well-known that the Hausdorff dimension of elliptic measures only depends on how the matrix is near the boundary and hence it is only necessary to assume the Lipschitz regularity around the boundary. In fact, in our proof we use the regularity only around the boundary, and so the theorem is still true assuming Lipschitz continuity in a small neighbourhood of $\partial \Omega$.

Similarly as in \cite[Proof of Theorem 1]{Wolff1993}, \cref{dimension elliptic measure} follows from a more quantitative result involving a good covering of a subset of the boundary with big elliptic measure, see \cref{only lambda dependence covering elliptic measure} for the precise statement.

One might think that this result could be obtained by the application of quasiconformal mappings. For symmetric matrices with determinant 1, the principal solution of the associated Beltrami equation, which depends only on the matrix, can act as a change of variables which inherits the extra regularity of the coefficients. Then, we obtain that the elliptic measure is the pushforward of the harmonic measure in the image domain. In the general case, we can obtain the elliptic measure as a pushforward of a harmonic measure using a quasiconformal change of variables which depends also on the Green function of the domain, as long as it satisfies the capacity density condition. See the article \cite{Guillen-quasiconformal} by the first author. In this case, the extra dependence of the Green function does not allow us to obtain extra regularity estimates for the change of variables, and the $\sigma$-finiteness of length can not be attained.

A key point, and the main difficulty in this paper is to obtain the lower bound
$$
\int_{\partial \widetilde \Omega} \log \frac{d\widetilde \omega^p}{d\sigma} (\xi) \, d\widetilde \omega^p (\xi) \geq -\textrm{const} > -\infty,
$$
with a bound independent of the smoothness, where $\widetilde \Omega$ is the modified domain appearing in the proof of \cref{covering elliptic measure} and $\widetilde \omega^p$ the elliptic measure in $\widetilde \Omega$ with respect to the matrix $A$. This lower bound is known in the harmonic case for general smooth domains in the plane. This fact is the key point in the study of the dimension in the planar case, see for example \cite{Carleson1985}, \cite{Jones1988}, \cite{Wolff1993}. Actually, the behavior of this integral is also crucial in the study of the dimension of harmonic measures in higher dimensions in \cite{Wolff1995}.

The first occurrence (as far as we know) of the use of this lower bound in the study of the dimension of the harmonic measure is in \cite{Carleson1985}. For the proof see \cite{Jones1988}, and for further details see \cite{Cufi2019-arxiv}.

The previous lower bound is used in the proof of \cref{covering elliptic measure} to obtain a subset with big elliptic measure as it is done in \cite{Jones1988} and \cite{Wolff1993}, and it can be deduced from
 \begin{equation}\label{key integral log bound}
 \int_{\partial \widetilde \Omega} \log |S\nabla g_p^T (\xi)|^2 \, d\widetilde \omega^p (\xi) \geq -\textrm{const} > -\infty,
\end{equation}
where $g_p^T$ is the Green function in $\widetilde \Omega$ with respect to the matrix $A^T$, and $S$ is the square root matrix of the symmetric part $A_0=(A+A^T)/2$, i.e., $S^TS=A_0$. In \cref{sec:log integral} we obtain also an upper bound for the integral above, see \cref{key integral log bound absolute lemma}.

A fundamental tool to obtain the estimate \rf{key integral log bound} is the relation
	\begin{equation*}
		 |\nabla g_p^T(y)| \gtrsim\frac{g_p^T(y)}{\dist (y, \partial \widetilde \Omega)} 
	\end{equation*}
near the boundary. In Reifenberg flat domains this estimate was obtained by Lewis, Lundström and Nyström in \cite{Lewis2008}, see \cref{comparability function gradient} below. The converse inequality is well-known for general domains and Hölder matrices and follows by Schauder estimates. Obtaining an inequality of this type for other domains may allow to apply the techniques exposed in this paper.

Very similar results to the one in \cref{dimension elliptic measure} about the dimension are true for $p$-harmonic measures, i.e., when the associated operator is the $p$-Laplacian $\divv \left( |\nabla u|^{p-2} \nabla u\right)$ for Reifenberg flat sets with small constant (see \cite{Lewis2013}), and simply connected sets in the plane (see \cite{Lewis2011}).

\section{Preliminaries and definitions}

\subsection{Notation}

\begin{itemize}
	\item We use $c,C\geq 1$ to denote constants that may depend only on the dimension and the constants appearing in the hypotheses of the results, and whose values may change at each occurrence.
	
	\item We write $a\lesssim b$ if there exists a constant $C\geq 1$ such that $a\leq Cb$, and $a\approx b$ if $C^{-1} b \leq a \leq C b$.
	
	\item If we want to stress the dependence of the constant on a parameter $\eta$, we write $a\lesssim_\eta b$ or $a\approx_\eta b$ meaning that $C=C(\eta)=C_\eta$.
	
	\item The ambient space is $\R^2$. However, some auxiliary results will be stated in general, i.e., in $\R^{n+1}$ for $n\geq 1$.
	
	\item The diameter of a set $E\subset \R^{n+1}$ is denoted by $\diam E \coloneqq \sup_{x,y\in E} |x-y|$.
	
	\item We denote by $B_r (x)$ or $B(x,r)$ the open ball with center $x$ and radius $r$, i.e., $B_r (x)=B(x,r)=\{y\in \R^{n+1} : |y-x|<r\}$. We denote $B_r \coloneqq B_r (0)$.
	
	\item Given a ball $B$, we denote by $r_B$ or $r(B)$ its radius, and by $c_B$ or $c(B)$ its center.
	
	\item We say that a matrix $A$ is Hölder continuous with exponent $\alpha \in (0,1]$ in a set $U$, or briefly $C^{0,\alpha} (U)$, if its coefficients are Hölder continuous with exponent $\alpha$. That is, there exists a constant $C_\alpha >0$ (called the Hölder seminorm) such that
	\begin{equation*}
		\left|a_{i j}(x)-a_{i j}(y)\right| \leq C_\alpha |x-y|^{\alpha} \text { for all } x, y \in U \text { and } 1 \leq i, j \leq n+1 .
	\end{equation*}
	For shortness we write $C^\alpha$ instead of $C^{0,\alpha}$ if $\alpha \in (0,1)$, and when $\alpha = 1$ we say ``Lipschitz continuous''. In this case we write $C_L$ instead of $C_1$, i.e.,
	\begin{equation*}
		\left|a_{i j}(x)-a_{i j}(y)\right| \leq C_L |x-y| \text { for all } x, y \in U \text { and } 1 \leq i, j \leq n+1 .
	\end{equation*}
	
	\item We say that a function $f$ is $\kappa$-Lipschitz in $U$ if $|f(x)-f(y)| \leq \kappa |x-y|$ for all $x,y\in U$.
	
	\item We denote the characteristic function of a set $E$ by $\characteristic_E$.
	
	\item Denote $\DD(\R^{n+1})$ the standard dyadic grid. That is, $\DD(\R^{n+1}) = \bigcup_{k\in \Z} \DD_k (\R^{n+1})$ where $\DD_k (\R^{n+1})$ is the collection of all cubes of the form
	$$
	\{ x \in \R^{n+1} : m_i 2^{-k} \leq x_i < (m_i + 1)2^{-k} \text{ for } i=1, \ldots, n+1\},
	$$
	where $m_i \in \Z$.
	
	\item Given $t>0$ and a set $E\subset \R^{n+1}$, we write $U_t (E) \coloneqq \{x\in \R^{n+1} : \dist (x,E) < t \}$ for the $t$-neighborhood $E$.
\end{itemize}

\subsection{CDC, NTA and Reifenberg flat domains}\label{sec:nta-reifenberg domains}

In this subsection we introduce the capacity density condition (CDC), non-tangentially accessible domains (NTA) and Reifenberg flat domains, the main object of our study.

\begin{definition}
	Let $K$ be a compact subset of $\Omega$. Its capacity is
	$$
	\Capacity (K,\Omega) = \inf \left\{\int_{\Omega} |\nabla u|^2 \, dx : u\in C^\infty_0 (\Omega) \text{ and } u\geq 1 \text{ on } K  \right\}.
	$$
\end{definition}

Just from the definition of capacity we obtain the following facts. For more properties see {\cite[2.2. Theorem]{Heinonen2006}}.

\begin{lemma}\label{basic properties condenser capacity}
	The set function $E\mapsto \Capacity (E,\Omega)$, $E\subset \Omega$, enjoys the following properties:
	\begin{enumerate}
		\item If $E_1 \subset E_2$, then $\Capacity (E_1, \Omega) \leq \Capacity (E_2, \Omega)$.
		\item If $\Omega_1 \subset \Omega_2$ are open and $E\subset \Omega_1$, then $\Capacity (E, \Omega_2) \leq \Capacity (E, \Omega_1)$.
	\end{enumerate}
\end{lemma}

\begin{definition}[CDC domain. {\cite[(11.20), (2.13)]{Heinonen2006}}]\label{def:general CDC}
	A domain $\Omega \subset \R^{n+1}$, $n\geq 1$, satisfies the capacity density condition, CDC for short, if there exist constants $c_0, r_0 >0$ such that 
	$$
	\Capacity \left(\overline{B(x_0,r)} \cap \Omega^c, B(x_0,2r) \right) \geq c_0 r^{n-1},
	$$
	for all $x_0\in \partial \Omega$, $x_0 \not = \infty$ and $r\leq r_0$.
\end{definition}

In order to define NTA domains we need to introduce some concepts. Given a domain $\Omega \subset \R^{n+1}$, $n\geq 1$, and a fixed constant $C$, we define:
\begin{itemize}
	\item \textit{$C$-Whitney ball}: A ball $B(x,r) \subset \Omega$ is a $C$-Whitney ball in $\Omega$ if $C^{-1}r < \dist (B(x,r), \partial \Omega) < C r$.
	\item \textit{$C$-Harnack chain}: For $p_1, p_2 \in \Omega$, a $C$-Harnack chain from $p_1$ to $p_2$ in $\Omega$ is a sequence of $C$-Whitney balls such that the first ball contains $p_1$, the last contains $p_2$, and such that consecutive balls intersect. The number of balls is called the length of the $C$-Harnack chain.
\end{itemize}

Consecutive balls in a $C$-Harnack chain must have comparable radius. Given a positive $L_A$-harmonic function in $\Omega$, a $C$-Harnack chain between two points $p_1, p_2 \in \Omega$ allows us (via Harnack's inequality) to obtain $u(p_1) \approx u(p_2)$, where the constant involved depends on $C$ and the length of the $C$-Harnack chain.

\begin{definition}[NTA domain. {\cite[Section 3]{Jerison1982}}]
	A domain $\Omega \subset \R^{n+1}$, $n\geq 1$, is called non-tangentially accessible, NTA for short, if there exists constants $C> 1$ and $r_0 >0$ such that:
	\begin{enumerate}
		\item \textit{Interior Corkscrew condition}: For any $\xi \in \partial \Omega$, $0<r< r_0$, there exists a point $A_r (\xi) \in \Omega$ such that $|A_r (\xi)-\xi|<r$ and $\dist(A_r (\xi), \partial \Omega)>C^{-1}r$. The point $A_r (\xi)=A(\xi,r)$ is called the Corkscrew point of the point $x$ at radius $r$.
		\item \textit{Exterior Corkscrew condition}: $\overline\Omega^c$ satisfies the interior Corkscrew condition.
		\item \textit{Harnack chain condition}: If $\varepsilon >0$ and $p_1,p_2\in \Omega$ satisfy that $p_1, p_2 \in \Omega \cap B(\xi, r_0 / 4)$ for some $\xi \in \partial \Omega$, $\dist(p_j, \partial \Omega)>\varepsilon$, and $|p_1 - p_2| < 2^k \varepsilon$, then there exists a Harnack chain from $p_1$ to $p_2$ of length $Ck$ and such that the diameter of each ball is bounded below by $C^{-1} \min_{j=1,2} \dist (p_j, \partial \Omega)$.
	\end{enumerate}
\end{definition}

\begin{rem}\label{rem:NTA are CDC}
	A domain with the exterior Corkscrew condition satisfies the capacity density condition, see \cite[Theorem 6.31]{Heinonen2006}. In particular, NTA domains satisfy the capacity density condition.
\end{rem}

\begin{definition}[Hausdorff distance]
	Given nonempty sets $A$ and $B$, we denote their Hausdorff distance $\dist_\HH (A,B)$ by
	$$
	\dist_\HH (A,B) = \max \left\{ \sup_{x\in A} \dist (x,B),\sup_{y\in B} \dist (y,A) \right\} .
	$$
\end{definition}

\begin{definition}[Reifenberg flat domain]\label{def:Reifenberg flat}
	Let $\Omega \subset \R^{n+1}$, $n\geq 1$, be an open set, and let $0 < \delta < 1/2$ and $r_0 >0$. We say that $\Omega$ is a $(\delta, r_0)$-Reifenberg flat domain if it satisfies the following conditions:
	\begin{enumerate}
		\item\label{def:Reifenberg flat 1}For every $x\in \partial \Omega$ and every $0<r\leq r_0$, there exists a hyperplane $\mathcal P (x,r)$ containing $x$ such that
		$$
		\dist_\HH \left(\partial \Omega \cap B(x,r), \mathcal P (x,r) \cap B(x,r)\right) \leq \delta r .
		$$
		
		\item\label{def:Reifenberg flat 2} For every $x\in \partial \Omega$ and every $0<r\leq r_0$, one of the connected components of 
		$$
		B(x,r) \cap \left\{ y\in \R^{n+1} : \dist (y, \mathcal P (x,r)) \geq 2\delta r \right\}
		$$
		is contained in $\Omega$ and the other is contained in $\R^{n+1} \setminus \Omega$.
	\end{enumerate}
\end{definition}

For small enough $\delta >0$ we have that a $(\delta, r_0)$-Reifenberg flat domain is also an NTA domain (see \cite[Section 3]{Kenig1997}) and it satisfies the capacity density condition, see \cref{rem:NTA are CDC}.

\subsection{Partial differential equations}\label{sec:PDE}

We want to study the elliptic equation
\begin{equation}\label{pde}
	L_{A} u(x) \coloneqq -\divv(A(\cdot) \nabla u(\cdot))(x) =0,
\end{equation}
which should be understood in the distributional sense. We simply write $L$ instead of $L_A$ when the matrix is clear from the context.

\begin{definition}\label{def:elliptic solution}
	We say that a function $u\in W_{\loc}^{1,2} (\Omega)$ is a solution of \rf{pde}, or $L_A$-harmonic, in an open set $\Omega \subset \R^{n+1}$, $n\geq 1$, if
	$$
	\int \langle A(y) \nabla u (y) , \nabla\varphi (y) \rangle \, dy = 0, \text{ for all } \varphi\in C_c^\infty (\Omega).
	$$
\end{definition}

By the De Giorgi-Nash-Moser theorem a solution $u\in W_{\loc}^{1,2} (\Omega)$ of \rf{pde} is locally Hölder continuous. If the matrix has Hölder regularity then the solution is locally $C^{1,\beta}$ for some $\beta \in (0,1)$, see \cite[Theorem 3.13]{Han2011}. Assuming Lipschitz regularity of the coefficients, $L_A$-harmonic functions enjoy more regularity. More precisely, weak solutions of $-\divv A\nabla u=0$ with $A$ Lipschitz are twice weakly differentiable, and there is a ``Caccioppoli type'' inequality for second derivatives.

\begin{theorem}[See {\cite[Theorem 8.8]{Gilbarg2001}}]\label{twice weak dif solutions}
	Let $u\in W^{1,2} (\Omega)$ be a weak solution of the equation $L_A u = 0$, see \rf{pde}, in $\Omega \subset \R^{n+1}$, $n\geq 1$, where $A$ is uniformly elliptic and Lipschitz continuous in $\Omega$. Then for any subdomain $\Omega^\prime \subset \subset \Omega$, we have $u\in W^{2,2} (\Omega^\prime)$ and
	\begin{equation}\label{cacciopoli type ineq for weak twice diff solutions}
	\|u\|_{W^{2,2} (\Omega^\prime)} \leq C \|u\|_{W^{1,2} (\Omega)},
	\end{equation}
	where $C$ depends on the dimension, the ellipticity constant and Lipschitz constant of $A$, and the value $\dist (\Omega^\prime , \partial \Omega)$.
\end{theorem}

If, in addition, the domain has boundary of class $C^2$, then $L_A$-harmonic functions are globally in $W^{2,2}$.

\begin{theorem}[See {\cite[Theorem 8.12]{Gilbarg2001}}]\label{global twice weak dif solutions}
	Assume, in addition to the hypotheses of \cref{twice weak dif solutions}, that $\Omega$ is bounded with $\partial \Omega$ of class $C^2$ and that there exists a function $\varphi \in W^{2,2} (\Omega)$ for which $u-\varphi \in W^{1,2}_0 (\Omega)$. Then we have also $u\in W^{2,2} (\Omega)$ and
	$$
	\|u\|_{W^{2,2} (\Omega)} \leq C \left(\|u\|_{L^2(\Omega)} + \|\varphi\|_{W^{2,2}(\Omega)}\right),
	$$
	where $C$ depends on the dimension, the ellipticity constant and Lipschitz constant of $A$, and $\partial \Omega$.
\end{theorem}

\begin{rem}\label{pointwise A-solution}
	Let $U\subset \Omega$ an open subset. Assuming $A\in C^{0,1} (U)$, then any weak $L_A$-harmonic function $u\in W^{1,2}_{\loc} (\Omega)$ (in fact $u\in W^{2,2}_{\loc} (U)$ by \cref{twice weak dif solutions}) satisfies $\divv A\nabla u=0$ a.e.\ in $U$. Indeed, since the matrix $A$ is differentiable a.e.\ (by Rademacher's theorem) and $u\in W^{2,2}_{\loc} (U)$ (by \cref{twice weak dif solutions}), we have $\divv A\nabla u \in L^1_{\loc} (U)$. Moreover, since $u$ is $L_A$-harmonic then $\int \divv A\nabla u \cdot \psi = \int A\nabla u \nabla \psi =0$ for any $\psi \in C^\infty_c (U) \subset C^\infty_c (\Omega)$. By the fundamental lemma of calculus of variations\footnote{The fundamental lemma of calculus of variations: If $U$ open set, $f\in L^1_{\loc}(U)$ and $\int f \phi=0$ for any $\phi\in C^\infty_c (U)$, then $f=0$ a.e.\ in $U$.} we conclude $\divv A\nabla u=0$ a.e.\ in $U$.
\end{rem}

The following theorem about the Hölder continuity of $L_A$-harmonic functions up to the boundary of regular enough domains will allow us to bound the elliptic measure on a specific domain (an annulus) by means of studying the Green function near the boundary.

\begin{theorem}[See {\cite[Corollary 8.36]{Gilbarg2001}}]\label{C^a to the boundary}
	Let $T$ be a (possibly empty) $C^{1,\alpha}$ boundary portion of a bounded domain $\Omega \subset \R^{n+1}$, $n\geq 1$, and suppose $u\in W^{1,2} (\Omega)$ is a weak solution of \rf{pde} in $\Omega$, where $A\in C^\alpha$ $(0<\alpha <1)$, such that $u= 0$ on $T$ (in the sense of $W^{1,2} (\Omega)$). Then $u\in C^{1,\alpha} (\Omega \cup T)$, and for any $\Omega' \subset \subset \Omega \cup T$ we have
	$$
	\|u\|_{C^{1,\alpha} (\Omega')} \coloneqq \|u\|_{1;\Omega'} + [u]_{1,\alpha; \Omega'} \leq C \sup_{\Omega} |u|,
	$$
	where
	$$
	\|u\|_{1;\Omega'} \coloneqq \sup_{\Omega'} |u| + \sup_{\Omega'} |\nabla u|, \quad 
	[u]_{1,\alpha; \Omega'} \coloneqq \sup_{\substack{x,y\in \Omega' \\ x\not=y}} \frac{|\nabla u (x)-\nabla u(y) |}{|x-y|^\alpha} ,
	$$
	for $C$ depending on $n$, the value of $\dist(\Omega', \partial \Omega \setminus T)$, the $C^{1,\alpha}$ character of $T$, and the ellipticity constants and the Hölder norm of the matrix $A$.
\end{theorem}

Since bounded Lipschitz functions are also Hölder continuous for any exponent $\alpha \in (0,1)$, the previous theorem remains true for matrices with bounded Lipschitz coefficients.

\subsubsection{Non-degeneracy of \texorpdfstring{$|\nabla u|$}{the gradient of solutions of PDE's} in Reifenberg flat domains with small constant}

Despite the following result applies to more general matrices (see \cite[Definition 1.1 and Lemma 3.35]{Lewis2008}), it is only stated in our setting for our purposes, i.e., for uniformly elliptic \rf{elliptic1}-\rf{elliptic2} and Hölder continuous matrices. In fact, the Hölder continuity is only used near the boundary.

\begin{lemma}[{\cite[Lemma 3.35]{Lewis2008}}]\label{comparability function gradient}
	Let $\Omega \subset \R^{n+1}$, $n\geq 1$, be a $(\delta , r_0)$-Reifenberg flat domain, $\xi\in \partial \Omega$, and $0<r<\min\{r_0 , 1\}$. Let $0<\alpha <1$ and $A\in C^\alpha (\{x : \dist(x,\partial \Omega)<10r_0\} )$ be a real uniformly elliptic (not necessarily symmetric) matrix with ellipticity constant $\lambda$ and Hölder seminorm $C_\alpha$. Suppose that $u$ is a positive $L_A$-harmonic function in $\Omega \cap B(\xi, 4r)$, that is continuous in $\overline{\Omega} \cap B(\xi, 4r)$, and that $u=0$ on $\partial \Omega \cap B(\xi, 4r)$. There exist $\hat \delta = \hat \delta (n,\lambda, C_\alpha, \alpha)$,  $ \gamma =  \gamma (n,\lambda, C_\alpha, \alpha)$ and  $\hat c = \hat c (n,\lambda, C_\alpha, \alpha)$ such that if $0<\delta \leq \hat \delta$, then
	$$
	\gamma^{-1} \frac{u(y)}{\dist(y, \partial \Omega)} \leq \left|\nabla u (y)\right| \leq  \gamma \frac{u(y)}{\dist(y, \partial \Omega)} \text{ whenever } y \in \Omega \cap B(\xi, r/\hat c) .
	$$
\end{lemma}

\begin{rem}\label{comparability function gradient Lipschitz}
	Note that if the matrix has bounded Lipschitz coefficients, i.e., $\alpha=1$, then the same result holds with constants depending on $n$, the ellipticity constant $\lambda$ and the value $C_L \|A\|_{L^\infty (\R^{n+1})}$, where $C_L$ is the Lipschitz seminorm of $A$. The value $C_L \|A\|_{L^\infty (\R^{n+1})}$ comes from the fact that bounded Lipschitz functions are Hölder continuous for any exponent. Indeed, a quick computation shows that the matrix $A$ is Hölder continuous with exponent $1/2$ with Hölder seminorm $C_{1/2} \coloneqq \left(2C_L \|A\|_{L^\infty (\R^{n+1})} \right)^{1/2}$.
\end{rem}

The comparability in \cref{comparability function gradient} will allow to bound the error terms in the study of the key term in \rf{key integral log bound absolute}.

\subsection{The fundamental solution and the Green function}

We denote by $\EE_x^A (y)$ the fundamental solution with pole at $x$ for $L_A$ in $\R^{n+1}$, $n\geq 1$, so that $L_A \EE_x^A (\cdot) = \delta_x$ in the distributional sense, where $\delta_x$ is the Dirac mass at the point $x\in \R^{n+1}$. We write $\EE_x (y)$ when the matrix $A$ is clear from the context. For a construction of the fundamental solution for real and uniformly elliptic matrices we refer to \cite{Hofmann2007} for higher dimensions, $\R^{n+1}$ with $n\geq 2$, and \cite[Appendix]{Kenig1985} for the planar case.

In higher dimensions the fundamental solution behaves ``in many senses'' as in the harmonic case, see \cite{Hofmann2007} for more details, but in the plane the situation is more delicate due to the change of sign of $\EE_0^{Id} (x) = \log |x|$ in $|x|=1$, the fundamental solution of the Laplacian in the plane.
 
Here we collect the formal definition and some properties of the fundamental solution in the plane.

\begin{definition}[{\cite[Definition 2.5]{Kenig2009}}]
	A function $\EE_x : \R^2 \to \R$ is called a fundamental solution for $L_A = \divv A\nabla \cdot$ with pole at $x$ if 
	\begin{enumerate}
		\item $\EE_x \in W_{\loc}^{1,2} (\R ^2\setminus \{x\}) \cap W_{\loc}^{1,p} (\R^2)$ for all $p<2$, and
		$$
		\int_{\R^2} \langle A(z) \nabla \EE_{x}(z) , \nabla \varphi(z) \rangle \, d z=-\varphi(x), \text{ for all } \varphi \in C_c^\infty (\R^2),
		$$
		\item $\left|\EE_{x}(y)\right|=\OO (\log |x-y|)$ as $|y| \rightarrow \infty$.
	\end{enumerate}
\end{definition}

The following result controls the fundamental solution far from the pole similarly as the fundamental solution for the harmonic case.

\begin{theorem}[{\cite[Theorem 2.6]{Kenig2009}}]\label{existence_fund_sol_pla}
	For each $x\in \R^2$ there exists a unique (modulo an additive constant) fundamental solution $\EE_x$ for $L_A$ with pole at $x$, and positive constants $C_1$, $C_2$, $R_1 < 1 < R_2$, which depend only on $\lambda$, such that
	\begin{equation*}
		%in columns
		\begin{aligned}C_{1} \log (1 /|x-y|)  &\leq{}\\ C_{1} \log (|x-y|) &\leq{} \end{aligned}\!
		\begin{gathered}-\EE_{x}(y)\\ \EE_{x}(y) \end{gathered}\!
		\begin{aligned}{}&\leq C_{2} \log (1 /|x-y|)\\ {}&\leq C_{2} \log (|x-y|)\end{aligned}
		\begin{aligned}{}& \text { for }|x-y|<R_{1}, \text{ and}\\ {}& \text{ for }|x-y|>R_{2}.\end{aligned}
	\end{equation*}
\end{theorem}

From the previous result and the maximum principle we obtain the following pointwise bound.

\begin{coro}\label{pointwise_bound_log}
	$|\EE_x (y)| \lesssim 1+ \left|\log |x-y|\right|$ for all $x,y\in \R^2$, where the constant depends on $\lambda$.
	
	\begin{proof}
		For $|x-y|<R_1$ and $|x-y|>R_2$, \cref{existence_fund_sol_pla} gives $|\EE_x (y)| \lesssim \left|\log |x-y|\right| \leq 1+ \left|\log |x-y|\right| $. In $\mathcal A=\{y \in \R^2 : R_1 < |x-y| < R_2\}$ we have $L\EE_x (y) = 0$, and hence by the maximum principle we obtain $C_2 \log R_1 \leq \EE_x (\cdot) \leq C_2 \log R_2$ in the annulus $\mathcal A$. So $|\EE_x (y)| \lesssim 1 \leq 1+ \left|\log |x-y|\right|$.
	\end{proof}
\end{coro}

We have the following relation between the fundamental solutions of the operators with matrices $A$ and $A^T$. The same holds in higher dimensions, see \cite[(3.43)]{Hofmann2007}.

\begin{lemma}[{\cite[Lemma 2.7]{Kenig2009}}]\label{transpose_sol}
	Fix $x,y\in \R^2$. Let $\EE_x$ be the fundamental solution for an elliptic operator $L_A$ with pole at $x$, and $\EE_y^T$ be the fundamental solution to the adjoint operator $L_{A^T}$ with pole at $y$. Then $\EE_x (y) = \EE_y^T (x)$.
\end{lemma}

Now we focus on Green's function. Given a bounded Wiener regular domain $\Omega \subset \R^{n+1}$, $n\geq 1$, and a uniformly elliptic matrix $A$, let $g_x = g_x^{\Omega,A}$ denote the Green function, in $\Omega$ with pole at $x\in \Omega$ with respect to the matrix $A$, constructed in \cite[Theorem 2.12]{Dong2009} in the planar case, and \cite[Theorem 1.1]{Gruter1982} in higher dimensions. We denote $g_x^T=g_x^{\Omega,A^T}$ the Green function with respect to the matrix $A^T$. In particular, $g_x$ satisfies
\begin{equation}\label{dirac delta definition green function}
	\int_\Omega   A (z) \nabla g_x (z) \nabla \varphi (z) \, dz = \varphi (x), \text{ for all } \varphi \in C_c^\infty (\Omega),
\end{equation}
and the following:
\begin{enumerate}
	\item $g_x (y) = g_y^T (x)$ for all $x,y\in \Omega$ and $x\not = y$. See \cite[Theorem 2.12 (2.18)]{Dong2009} in the plane, and \cite[Theorem 1.3]{Gruter1982} in higher dimensions.
	\item\label{sobolev 1-2 regularity green function} For each $x\in \Omega$ and any $0<r<\dist(x,\partial \Omega)$, $g_x \in W^{1,2} (\Omega \setminus B_r (x))$. See \cite[Theorem 2.12 (2.15)]{Dong2009} in the plane, and \cite[Theorem 1.1 (1.3)]{Gruter1982} in higher dimensions.
\end{enumerate}
If the domain $\Omega$ has boundary of class $C^2$ and the matrix is Lipschitz continuous in a neighborhood $U_{2s} (\partial \Omega) = \{x\in \R^{n+1} : \dist(x,\partial \Omega) < 2s \}$, then the Green function also satisfies:
\begin{enumerate}
	\setcounter{enumi}{2}
	\item For each $x\in \Omega$ and any $0<r<\min\{s, \dist(x,\partial \Omega)\}$, $g_x \in W^{2,2} (\Omega \cap U_r (\partial\Omega))$. This is a consequence of \rf{sobolev 1-2 regularity green function} by \cref{twice weak dif solutions,global twice weak dif solutions}.
\end{enumerate}

Next we show that $L_{A^T} g^T_x = \omega_{\Omega,A}^x - \delta_x$ in the distributional sense. With this identity we can move from integrating on the boundary to the interior of the set.

\begin{lemma}
	Let $\Omega \subset \R^{n+1}$, $n\geq 1$, be a bounded Wiener regular domain and $\varphi\in C(\overline \Omega) \cap W^{1,2}(\Omega)$. Then
	\begin{equation}\label{boundary integral to interior integral}
		\int_{\partial \Omega} \varphi (\xi) \, d \omega^x_{\Omega,A} (\xi) -\varphi(x) = -\int_{\Omega} A^T (z)  \nabla g_x^T (z) \nabla \varphi (z) \, dz, \text{ for a.e.\ } x\in\Omega.
	\end{equation}
	\begin{proof}[Sketch of proof]
		In higher dimensions this is proved in \cite[(2.6)]{Azzam2022}. Here we detail the differences in the planar case.
		
		As $\Omega$ is bounded Wiener regular and $\varphi \in C(\overline\Omega)\cap W^{1,2} (\Omega$), the $L_A$-harmonic function $u$ solving the Dirichlet problem with boundary data $\varphi |_{\partial \Omega}$ can be taken to be in $C(\overline \Omega)\cap W^{1,2} (\Omega)$ and $u-\varphi \in W^{1,2}_0(\Omega)$.  Indeed, by the Lax-Milgram theorem in $W_0^{1,2} (\Omega)$ there is a unique function $v\in W_0^{1,2} (\Omega)$ with $\int_\Omega A \nabla v \nabla \vartheta = - \int_\Omega A \nabla \varphi \nabla \vartheta$ for any $\vartheta \in W_0^{1,2} (\Omega)$. Taking $u=v+\varphi$ it is clear that $L_A u = 0$, $u\in W^{1,2} (\Omega)$ and $u-\varphi \in W^{1,2}_0(\Omega)$. On the other hand, since $\Omega$ is bounded Wiener regular, the $L_A$-harmonic extension $u$ of $\varphi |_{\partial\Omega}$ is continuous up to the boundary, see \cite[Theorem 6.27]{Heinonen2006}. Moreover, by the definition of elliptic measure and the uniqueness of solutions we have
        $$
        u(x) =\int_{\partial\Omega} \varphi (\xi) \, d\omega^x_{\Omega,A} (\xi) .
        $$
		
        Write
        $$
        \int_\Omega A^T \nabla g_x^T \nabla \varphi =  
        \int_\Omega A^T \nabla g_x^T \nabla u 
        + \int_\Omega A^T \nabla g_x^T \nabla (\varphi - u)
        \eqqcolon\rom{1} + \rom{2} .
        $$
        The same proof of \cite[(2.10) and (2.12)]{Azzam2022} applies also in the plane to have that the left-hand side integral is absolutely convergent and $\rom{2}=\varphi(x)-u(x)$ for a.e.\ $x\in\Omega$, replacing the use of (2.8) and (2.9) in \cite{Azzam2022} by the fact that the Green function satisfies $\nabla g_z \in L^p (\Omega)$ for all $p\in [1,2)$ and $g_z\in W_{\loc}^{1,2}(\Omega\setminus \{z\})$, see \cite[Remark 2.19 and (3.66)]{Dong2009} and \cref{sobolev 1-2 regularity green function} in \cpageref{sobolev 1-2 regularity green function} respectively. Using that $|g_x^T (z)|\lesssim \left|\log|z-x|\right|$ when $|x-z|\leq 2\varepsilon$ is small enough, see \cite[(2.17)]{Dong2009}, in the planar case the term $\rom{1}_\varepsilon^2$ defined in \cite[p. 10855]{Azzam2022} is controlled by
        $$
        |\rom{1}_\varepsilon^2| \lesssim \frac{\left|\log\varepsilon\right|}{\varepsilon} \int_{B_{2\varepsilon}(x)} |\nabla u| \lesssim \varepsilon\left|\log\varepsilon\right| \mathcal{M}(\nabla u\characteristic_{\Omega})(x),
        $$
        and as $\varepsilon \left|\log \varepsilon\right| \to 0$ as $\varepsilon \to 0$, the same proof there implies that $\rom{1}=0$ for a.e.\ $x\in \Omega$. Therefore, \rf{boundary integral to interior integral} holds also in the planar case.
	\end{proof}
\end{lemma}

We will show below that from the equality \rf{boundary integral to interior integral} it follows that
\begin{equation}\label{def green function}
g_y (x) =
	-\EE_y (x) + \int_{\partial \Omega} \EE_y (\xi) \, d\omega_{\Omega,A}^x (\xi), \text{ for all } x,y\in \Omega.
\end{equation}
(Probably this is already known but we will show the full details in the plane for completeness). Recall that $x \mapsto \int_{\partial \Omega} \EE_y (\xi) \, d\omega_{\Omega,A}^x (\xi)$ is the $L_A$-harmonic extension of $\EE_y$ inside $\Omega$. Assuming that $g_y (x)=0$ if $y\not\in \Omega$, we have that \rf{def green function} also holds in this case, since $\Omega$ is Wiener regular and therefore the Green function is continuous through the boundary. Moreover, by \rf{def green function} and since $\Omega$ is bounded, the Green function also satisfies:
\begin{enumerate}
    \setcounter{enumi}{3}
    \item For each $x\in \Omega$, $g_x (y)\geq 0$ for any $y\in \Omega \setminus \{x\}$.
\end{enumerate}
This was already proved in higher dimension in \cite[Theorem 1.1]{Gruter1982}. However, since the situation is more delicate for unbounded planar domains due to the logarithmic behaviour of the fundamental solution, we only consider bounded planar Wiener regular domains.

\begin{proof}[Proof of \rf{def green function} in the planar case]
    Let $0<\varepsilon\ll\min\{|x-y|,\dist(x,\partial\Omega),\dist(y,\partial\Omega)\}$,  $\psi^y=\psi^y_\varepsilon \in C^\infty_c (B_{2\varepsilon}(y))$ such that $\psi^y_\varepsilon = 1$ in $B_\varepsilon(y)$ and $|\nabla \phi^y_\varepsilon|\lesssim 1/\varepsilon$, and $\psi^x = \psi^x_\varepsilon$ defined analogously.

    Applying \rf{boundary integral to interior integral} to $(1-\psi^y) \EE_y \in C(\overline{\Omega})\cap W^{1,2}(\Omega)$ and using that $1-\psi^y =1$ in $\partial \Omega \cup \{x\}$, we have
    $$
    \int_{\partial\Omega} \EE_y (\xi)\, d\omega_{\Omega,A}^x (\xi) - \EE_y (x) 
    = -\int_\Omega A^T(z) \nabla g_x^T (z)\nabla ((1-\psi^y)\EE_y)(z)\, dz.
    $$

    Write the right-hand side term as
    $$
    \begin{aligned}
    -\int_\Omega A^T \nabla g_x^T\nabla ((1-\psi^y)\EE_y)\, dz
    =&
    \int_{B_{2\varepsilon}(y)\setminus B_\varepsilon(y)} A^T \nabla g_x^T \nabla \psi^y \cdot \EE_y\, dz
    + \int_{B_{2\varepsilon}(y)} A \nabla \EE_y\nabla g_x^T\cdot\psi^y\, dz\\
    &-\int_{B_{2\varepsilon}(x)\setminus B_\varepsilon(x)} A\nabla \EE_y \nabla \psi^x \cdot g_x^T \, dz
    -\int_\Omega A \nabla \EE_y\nabla ((1-\psi^x)g_x^T)\, dz\\
    \eqqcolon&\,\rom{1}_\varepsilon+\rom{2}_\varepsilon+\rom{3}_\varepsilon+\rom{4}_\varepsilon.
    \end{aligned}
    $$
    Using that $|g_x^T (z)| \lesssim \left|\log |x-z|\right|$ and $|\EE_y(z)| \lesssim \left|\log |y-z|\right|$ when $|x-z|\leq 2\varepsilon$ and $|y-z|\leq 2\varepsilon$ for small enough $\varepsilon>0$, see \cite[(2.17)]{Dong2009} and \cref{existence_fund_sol_pla}, the terms $\rom{1}_\varepsilon$ and $\rom{3}_\varepsilon$ are bounded by
    \begin{equation}\label{terms I and III}
    |\rom{1}_\varepsilon|+|\rom{3}_\varepsilon| 
    \lesssim \varepsilon\left|\log \varepsilon\right|\left\{
    \left(\avint_{B_{2\varepsilon}(y)} |\nabla g_x^T|^2\, dz\right)^{1/2}
    + \left(\avint_{B_{2\varepsilon}(x)} |\nabla \EE_y|^2\, dz\right)^{1/2}
    \right\}.
    \end{equation}
    For the bound of $\rom{2}_\varepsilon$, since $g_x^T$ is $L_{A^T}$-harmonic in $\Omega\setminus \{x\}\supset B_{10\varepsilon}(y)$, there exists $p=p(\lambda)>2$ such that
    $$
    \left(\avint_{B_{2\varepsilon} (y)} |\nabla g_x^T|^p \, dz\right)^{1/p} \lesssim_\lambda \left(\avint_{B_{4\varepsilon} (y)} |\nabla g_x^T|^2 \, dz\right)^{1/2},
    $$
    see \cite[Lemma 1.1.12]{Kenig1994}, and let $1\leq q<2$ be its Hölder exponent conjugate, i.e., $1/p+1/q=1$. By Hölder's inequality and the choice of $p>2$, the term $\rom{2}_\varepsilon$ is controlled by
    \begin{equation*}
    |\rom{2}_\varepsilon|=\left|\int_{B_{2\varepsilon}(y)} A \nabla \EE_y\nabla g_x^T\cdot\psi^y\, dz\right| 
    \lesssim \varepsilon^2  
    \left(\avint_{B_{2\varepsilon} (y)} |\nabla \EE_y|^q \, dz\right)^{1/q}
    \left(\avint_{B_{4\varepsilon} (y)} |\nabla g_x^T|^2 \, dz \right)^{1/2}.
    \end{equation*}
    By \cite[Theorem 0.1]{Chanillo1992} we have $\left(\avint_{B_{2\varepsilon} (y)} |\nabla \EE_y|^q \, dz\right)^{1/q}\lesssim_\lambda 1/\varepsilon$, and so
    \begin{equation}\label{term II}
    |\rom{2}_\varepsilon|
    \lesssim \varepsilon\left(\avint_{B_{4\varepsilon} (x)} |\nabla g_x^T|^2 \, dz \right)^{1/2}.
    \end{equation}
    Since $g_x^T\in W_{\loc}^{1,2} (\Omega\setminus \{x\})$ and $\EE_y \in W_{\loc}^{1,2} (\Omega\setminus \{y\})$, in particular $|\nabla g_x^T|^2 \in L_{\loc}^1 (\Omega\setminus \{x\})$ and $|\nabla \EE_y|^2 \in L_{\loc}^1 (\Omega\setminus \{y\})$, by the Lebesgue differentiation theorem we have that $\avint_{B_{4\varepsilon} (y)} |\nabla g_x^T|^2\to |\nabla g_x^T (y)|^2$ for a.e.\ $y\in \Omega$ and $\avint_{B_{4\varepsilon} (x)} |\nabla \EE_y|^2\to |\nabla \EE_y (x)|^2$ for a.e.\ $x\in \Omega$ respectively. That is, by \rf{terms I and III} and \rf{term II} we have $|\rom{1}_\varepsilon| + |\rom{2}_\varepsilon| + |\rom{3}_\varepsilon|=0$ a.e.\ $x,y\in \Omega$.
    
    On the other hand, from the Dirac delta property of the fundamental solution, $(1-\psi^x)g_x^T\in W^{1,2}_0 (\Omega)$ and the density of $C^\infty_c (\Omega)\subset C^\infty_c (\R^{n+1})$ in $W^{1,2}_0 (\Omega)$, we obtain $\rom{4}_\varepsilon=(1-\psi^x (y))g_x^T(y)=g_y (x)$, and \rf{def green function} is proved for a.e.\ $x,y\in \Omega$. By continuity, it also holds for all $x,y\in \Omega$.
\end{proof}

To end this section, we see how the Green function is related to the density of the elliptic measure in smooth domains. Assume now $A$ is Lipschitz continuous in an open neighborhood of $\partial \Omega$, say $U_s (\partial \Omega) = \{x\in \R^{n+1} : \dist (x, \partial \Omega)<s\}$. Under this assumption $A$ is differentiable a.e.\ by Rademacher's theorem, and $L_A$-harmonic functions are in $W^{2,2}$ by \cref{global twice weak dif solutions}.

\begin{lemma}
    Let $\Omega \subset \R^{n+1}$, $n\geq 1$, be a bounded domain with smooth boundary (and hence Wiener regular) and $A\in C^{0,1} (U_s (\partial \Omega))$. The elliptic measure $\omega^p_{\Omega,A}$ can be written as 
    \begin{equation}\label{elliptic measure density}
    d\omega_{\Omega,A}^p = - \langle A^T \nabla g_p^T , \nu \rangle \, d\sigma , \text{ for a.e.\ } p \in \Omega \setminus U_{2s} (\partial \Omega),
    \end{equation}
    where $\nu$ is the unit outer normal to $\partial \Omega$ and $\sigma$ is the surface measure on $\partial \Omega$.
    \begin{proof}
        Let $\varphi \in C^\infty_c (U_s (\partial \Omega))$ and set $\phi_p (z) = g_p^T (z) + \EE_p^T (z)$ for $z\in \R^{n+1}$, see \rf{def green function} when $z\in \Omega$. Notice that $\phi_p (z) = -\EE^T_p (z)$ in $\Omega^c$, and hence $\phi_p$ is $L_{A^T}$-harmonic in $\Omega$ and $\overline{\Omega}^c$.
        
        The claim follows since the right-hand side of \rf{boundary integral to interior integral} is
        $$
        - \int_{\partial \Omega} \varphi(\xi) \langle A^T (\xi) \nabla g_p^T (\xi), \nu (\xi) \rangle \, d\sigma (\xi) - \varphi(p).
        $$
        Indeed, since $A$ is differentiable a.e.\ in $\supp \varphi$ (Rademacher's theorem), $g_p^T \in W^{2,2} (\supp \varphi)$ (\cref{global twice weak dif solutions}) and $\varphi \in C^\infty_c (U_s (\partial \Omega))$, in particular $A^T \nabla g_p^T \cdot \varphi \in W^{1,2} (\supp \varphi)$. As the integration by parts formula holds for $W^{1,2}$ functions, then
        \begin{equation*}
            \begin{aligned}
                - \int_\Omega A^T \nabla  g_p^T \nabla \varphi   = & \int_\Omega \varphi \divv \left( A^T \nabla g_p^T \right)  - \int_\Omega \divv \left( A^T \nabla g_p^T \cdot \varphi\right) \\
                =& \int_\Omega \varphi \divv \left( A^T  \nabla g_p^T \right) - \int_{\partial \Omega} \varphi \langle A^T \nabla g_p^T, \nu \rangle \, d\sigma  \\
                =&  \int_{\R^{n+1}} \varphi \divv \left( A^T \nabla \phi_p \right) 
                -\int_{\R^{n+1}} \varphi \divv \left( A^T \nabla \EE^T_p \right) - \int_{\partial \Omega} \varphi \langle A^T  \nabla g_p^T , \nu  \rangle \, d\sigma  \\
                =& \int_{\R^{n+1}} \varphi \divv \left( A^T \nabla \phi_p \right) - \varphi(p) - \int_{\partial \Omega} \varphi \langle A^T  \nabla g_p^T , \nu \rangle \, d\sigma,
            \end{aligned}
        \end{equation*}
        where in the last equality we used that $\varphi$ has compact support and $\int \varphi \divv \left( A^T \nabla \EE_p^T \right) = \varphi(p)$ by the definition of the fundamental solution. Since $\phi_p$ is $L_{A^T}$-harmonic in $\R^{n+1} \setminus \partial \Omega$, $\HH^{n+1} (\partial \Omega)=0$ and $A\in C^{0,1} (\supp \varphi)$, we have that
        $\divv \left( A^T \nabla \phi_p \right)=0$ a.e.\ in $\supp \varphi$ by \cref{pointwise A-solution}, and the claim follows.
    \end{proof}
\end{lemma}

\section{Main Lemma and preliminary reductions}\label{weakening of the main lemma}

As in \cite{Wolff1993}, \cref{dimension elliptic measure} will follow from the following more quantitative result.

\begin{mlemma}\label{only lambda dependence covering elliptic measure}
	Let $\Omega \subset \R^2$ be a bounded $(\delta, r_0)$-Reifenberg flat domain, a point $p\in \Omega$ with $\dist(p,\partial \Omega) > r_0$, and $A$ be a real uniformly elliptic (not necessarily symmetric) matrix with ellipticity constant $\lambda$, and suppose also that $A$ is $\kappa$-Lipschitz in $U_{r_0} (\partial\Omega)\coloneqq \{x\in \R^2 : \dist(x,\partial \Omega) <r_0\}$. For a given $0<r\leq 1$ satisfying $r\kappa \|A\|_{L^\infty (\R^2)} \leq 1$, there exists $\delta_0 = \delta_0 (\lambda) > 0$ such that for every $0<\delta \leq \delta_0$ we have the following:
	
	For any $0< \tau < 1$, sufficiently large $M$, and $\rho \in \left(0, r/M\right)$ there is a set $F \subset \partial \Omega$ such that $\omega_{\Omega,A}^p (F) \geq C^{-1} \tau$ and a countable covering $F\subset \bigcup_{i} B(z_i, r_i)$ where
	\begin{enumerate}
		\item $\sum_{i} r_i \leq C M^\tau$,
		\item $\sum_{\{i \,:\, r_i > \rho\}} r_i \leq C M^{-1}$,
	\end{enumerate}
	with universal constant $C$. 
\end{mlemma}

\begin{rem}
    Given $M$ sufficiently large to satisfy the conclusions of the lemma, the particular choice $\rho = r/(2M)$ yields that the number of balls $B(z_i,r_i)$ with $r_i > r/(2M)$ is universally bounded, that is, $\sum_{\{i \,:\, r_i > r/(2M)\}} 1 \leq 2C/r$.
\end{rem}

By means of a linear deformation of the plane (see \cref{sec:the change of variables} below) and a rescaling, we see that it suffices to prove the following weaker lemma to obtain \cref{only lambda dependence covering elliptic measure}.

\begin{lemma}[Weak form of \cref{only lambda dependence covering elliptic measure}]\label{covering elliptic measure}
	Let $\Omega\subset\R^2$, $p\in\Omega$ and $A$ as in \cref{only lambda dependence covering elliptic measure}. Suppose also that $A_0 = \frac{A+A^T}{2}$ is of the form $A_0=R^TBR$ with $R\in C^{0,1} (U_{r_0}(\partial \Omega))$ a rotation, and $B \in C^{0,1} (U_{r_0}(\partial \Omega))$ diagonal. Then there exists $\delta_0 = \delta_0 (\lambda, \kappa\|A\|_{L^\infty (\R^2)}) > 0$ such that for every $0<\delta \leq \delta_0$ we have the following:
	
	For any $0< \tau < 1$, sufficiently large $M$ (how large depends on $\tau$ and on the constants in the hypothesis), and $\rho \in (0,1/M)$ there is a set $F \subset \partial \Omega$ such that $\omega_{\Omega,A}^p (F) \geq C^{-1} \tau$ and a countable covering $F\subset \bigcup_i B(z_i, r_i)$ with
	\begin{enumerate}[label=({\arabic*}*),ref={\arabic*}*]
		\item\label{covering elliptic measure cond1} $\sum_i r_i \leq C M^\tau$,
		\item\label{covering elliptic measure cond2} $\sum_{\{i \, : \, r_i > \rho\}} r_i \leq C M^{-1}$,
	\end{enumerate}
	with universal constant $C$.
\end{lemma}

We remark that this weaker form replaces the assumption on the parameter $r$ by the additional assumption $A_0 = R^T B R$, and allows $\delta_0$ to depend also on $\kappa \|A\|_{L^\infty}$, and $\rho$ to be in $(0,1/M)$. For the proof of \cref{covering elliptic measure}, see \cref{proof of weak version of main lemma}.

\subsection{The change of variables in Reduction 1}\label{sec:the change of variables}

In this subsection we first collect some auxiliary results about changes of variables that will be useful to prove some technical lemmas, secondly we construct the precise linear deformation in the plane that allows the reduction from \cref{only lambda dependence covering elliptic measure} to \cref{covering elliptic measure}, and finally we see how it distorts planar Reifenberg flat domains.

\subsubsection{Linear changes of variables}

We will see how $L_A$-harmonic functions behave under linear changes of variables. See \cite[Lemmas 3.8 and 3.9]{Azzam2019} for a detailed proof of the following two results.

\begin{lemma}\label{linear change of variables}
	Let $D \in \R^{(n+1)\times (n+1)}$ be a constant matrix with $\det D \not = 0$, $n\geq 0$. A function $f$ is $L_A$-harmonic in $\Omega$ if and only if $\widetilde f = f \circ D$ is $L_{\widetilde A}$-harmonic in $D^{-1}(\Omega)$, where $\widetilde A (\cdot) = D^{-1} A(D \cdot) (D^{-1})^T$ and $D^{-1}(\Omega) = \{ D^{-1} x : x\in \Omega\}$.
\end{lemma}

By the definition of elliptic measure, the previous lemma implies the following relation of elliptic measures under a linear change of variables. 

\begin{coro}\label{elliptic measure deformation}
	Let $D\in\R^{(n+1)\times(n+1)}$ be a constant matrix such that $\det D \not = 0$, $n\geq 0$, and let $\Omega$ be a Wiener regular domain. Let $\omega = \omega_{\Omega, A}$ be the elliptic measure in $\Omega$ with matrix $A$, and $\widetilde \omega = \omega_{D^{-1}(\Omega), \widetilde A}$ where $\widetilde A (\cdot) = D^{-1} A(D\cdot) \left(D^{-1}\right)^T$. Then $\omega^x (E) = \widetilde \omega^{D^{-1}x} (D^{-1} (E))$ for every $x\in\Omega$ and $E\subset\partial\Omega$.
\end{coro}

\subsubsection{Lipschitz diagonalization of symmetric matrices in the plane}

In the study of the integral \rf{key integral log bound absolute} in \cref{sec:log integral} we will use that after a suitable linear change of variables $D$, the symmetric part of the matrix $\widetilde A$ in \cref{elliptic measure deformation} diagonalizes in the form $R^T BR$, where $R$ is a Lipschitz rotation and $B$ is Lipschitz diagonal. In this subsection we see that we can always reduce to this case.

We need to follow this strategy because in general it is not true that Lipschitz elliptic symmetric matrices diagonalize in the aforementioned form, as we can see in the following example.

\begin{example}
	Let $A_1, A_2$ be two constant symmetric matrices diagonalizing with different eigenvectors, and $f:\R^{n+1}\to \R$ be a Lipschitz function with $f\characteristic_{\{|x|<1\}} <0$, $f\characteristic_{\{|x|>1\}} >0$ and $\|f\|_\infty \leq \varepsilon$ for small enough fixed constant $\varepsilon >0$.
	
	Set $A(x) \coloneqq Id + f(x)\characteristic_{\{|x|<1\}} (x) A_1 + f(x) \characteristic_{\{|x|>1\}} (x) A_2$ and take $\varepsilon >0$ small enough to ensure the ellipticity condition on the matrix $A$. Moreover, with this choice of the function $f$ we have that the matrix $A$ has Lipschitz coefficients.
	
	Let $v_1$ be an eigenvector of $A_1$ with eigenvalue $\mu_1$, i.e., $A_1 v_1 = \mu_1 u_1$, and let $v_2$ be an eigenvector of $A_2$ with eigenvalue $\mu_2$, i.e., $A_2 v_2 = \mu_2 v_2$. Then, for $|x|<1$ the vector $u$ is an eigenvector of $A$,
	$$
	A(x) v_1 = (Id + f(x) A_1)v_1 = (1+f(x)\mu_1)v_1,
	$$
	and for $|x|>1$ the vector $v$ is an eigenvector of $A$,
	$$
	A(x) v_2 = (Id + f(x) A_2)v_2 = (1+f(x)\mu_2)v_2.
	$$
	From this we get that the matrix $A$ diagonalizes with the same basis as $A_1$ if $|x|<1$, and with the same basis as $A_2$ if $|x|>1$, whence we obtain that the basis is not continuous.
\end{example}

In the following lemma we see that we can avoid the situation seen in the previous example by using a linear change of variables.

\begin{lemma}\label{symmetric to RBR}
	Let $U\subset \R^2$ be a set. Let $A \in C^{0,1} (U)$ be a uniformly elliptic and symmetric $2\times 2$ matrix with ellipticity constant $\lambda$, and let $D =\begin{pmatrix}
		1/K & 0 \\
		0 & 1
	\end{pmatrix}$ with $K^2 \geq \lambda^2+\lambda$. Then the matrix $\widetilde A (\cdot) = D^{-1} A(D \cdot) D^{-1}$ is of the form $\widetilde A = R^T B R \in C^{0,1} (D^{-1}(U))$, with $B\in C^{0,1} (D^{-1}(U))$ diagonal and $R\in C^{0,1} (D^{-1}(U))$ a rotation.

\begin{proof}
	Denote the matrix $A(x) = \begin{pmatrix}
	a(x) & b(x) \\
	b(x) & d(x)
\end{pmatrix}$, and let
$$
\widetilde A  = D^{-1} (A \circ D) D^{-1} 
= \begin{pmatrix}
	K^2 \widetilde a  & K \widetilde b  \\
	K\widetilde b  & \widetilde d 
\end{pmatrix},
$$
where we write $\widetilde a (x) = a(Dx)$ and the analogous expressions for the other elements of the matrix. We want to see that when $K^2\geq \lambda^2 + \lambda$ we can write $\widetilde A=R^T B R$ where $B$ is Lipschitz and diagonal, and $R$ is a Lipschitz rotation matrix.

The eigenvalues $\lambda_\pm = \lambda_\pm (\cdot)$ of $\widetilde A$ are
\begin{equation}\label{eigenvalues definition}
\lambda_{\pm} = \frac{K^2 \widetilde a + \widetilde d \pm \sqrt {\left(K^2 \widetilde a - \widetilde d\right)^2 + 4K^2 \widetilde b ^2}}{2} ,
\end{equation}
and we want to see that they are Lipschitz if $K^2\geq \lambda^2 + \lambda$. Note that 
\begin{equation}\label{sum of eigenvalues property}
\lambda_+ + \lambda_- = K^2\widetilde a + \widetilde d.
\end{equation}

For shortness, let $f = \left(K^2 \widetilde a  - \widetilde d  \right)^2 + 4K^2 \widetilde b^2$ be an auxiliar function. Note that since $a,b,d \in C^{0,1} (U) \cap L^\infty (U)$, and so $\widetilde a, \widetilde b, \widetilde d \in C^{0,1} (D^{-1} U) \cap L^\infty (D^{-1} U)$, we have $f\in C^{0,1} (D^{-1} U) \cap L^\infty (D^{-1} U)$.

For $x,y\in D^{-1} U$, i.e., $Dx,Dy\in U$,
\begin{equation}\label{Lip vaps part 1}
\begin{aligned}
	2|\lambda_\pm (x) - \lambda_\pm (y)|  &\leq K^2 |\widetilde a (x) - \widetilde a (y)| + |\widetilde d (x) - \widetilde d (y)| + |\sqrt{f(x)}- \sqrt{f(y)}| \\
	 &\leq C|x-y| + |\sqrt{f(x)}- \sqrt{f(y)}| = C|x-y| + \frac{|f(x)- f(y)|}{|\sqrt{f(x)}+ \sqrt{f(y)}|} .
\end{aligned}
\end{equation}

Since $a \geq \lambda ^{-1}$ and $d \leq \lambda$ by ellipticity (indeed $a,d\approx_\lambda 1$), and so $\widetilde a \geq \lambda ^{-1}$ and $\widetilde d \leq \lambda$, in particular $K^2 \widetilde a - \widetilde d \geq K^2 \lambda^{-1}-\lambda$. Since $K^2 \geq \lambda^2+\lambda$, we obtain that
\begin{equation}\label{lower bound difference of elements new matrix}
	K^2 \widetilde a - \widetilde d \geq K^2 \lambda^{-1}-\lambda \geq 1 ,
\end{equation}
and with this we have
\begin{equation}\label{Lip vaps part 2}
\sqrt{f } = \sqrt{\left(K^2 \widetilde a  - \widetilde d  \right)^2 + 4K^2 \widetilde b^2} \geq K^2 \widetilde a - \widetilde d \geq 1 .
\end{equation}

Combining the estimates \rf{Lip vaps part 1} and \rf{Lip vaps part 2} we get
$$
|\lambda_\pm (x) - \lambda_\pm (y)| \lesssim |x-y| + |f(x)-f(y)| \lesssim |x-y|,
$$
i.e., $\lambda_\pm$ are Lipschitz.

It remains to see that the matrix diagonalizes in the form $R^TBR$ and that the eigenvectors are also Lipschitz. Let $u^\pm = (u_1^\pm, u_2^\pm)$ be the eigenvectors of the eigenvalues $\lambda_\pm$. Hence, $\left( \widetilde A - \lambda_\pm Id \right) u^\pm = 0$, i.e.,
\begin{equation}\label{eigenvalue}
\begin{cases}
	(K^2 \widetilde a - \lambda_\pm) u_1^\pm + K\widetilde b u_2^\pm =0, \\
	K\widetilde b u_1^\pm + (\widetilde d -\lambda_\pm) u_2^\pm = 0.
\end{cases}
\end{equation}

Consider the vectors $v^+ = \left( \lambda_+ - \widetilde d, K\widetilde b \right)$ and $v^- = \left( K\widetilde b, \lambda_- - K^2 \widetilde a \right)$, which are clearly Lipschitz by the preceding discussion. We claim that $v^+$ and $v^-$ satisfy \rf{eigenvalue}. Indeed, $v^+$ satisfy the second equality in \rf{eigenvalue} by the definition of $v^+$, and the first equality follows from the definition of the eigenvalues $\lambda_\pm$ and the equality $\lambda_+ - \widetilde d = K^2\widetilde a - \lambda_-$, see \rf{sum of eigenvalues property}. The vector $v^-$ satisfy \rf{eigenvalue} by the same reason.

By \rf{sum of eigenvalues property}, we can write $v^+ = \left( K^2\widetilde a - \lambda_- , K \widetilde b \right)$ (and so $v^-$ and $v^+$ are orthogonal), and hence $\|v^+\|=\|v^-\|$. Moreover,
\begin{equation*}
	\|v^\pm\|^2
 = \left( K^2 \widetilde a - \lambda_- \right)^2 + K^2 \widetilde b^2
 \geq \left( K^2 \widetilde a - \lambda_- \right)^2 
  = \left( \frac{K^2 \widetilde a - \widetilde d + \sqrt{f}}{2} \right)^2 \geq 1,
\end{equation*}
by \rf{lower bound difference of elements new matrix} and \rf{Lip vaps part 2}. To conclude, set the unitary vectors $u^\pm \coloneqq \frac{v^\pm}{\|v^\pm\|}$. They are orthonormal and hence we conclude that $\widetilde A = R^T B R$ with
$$
B=
\begin{pmatrix}
\lambda_- & 0 \\
0 & \lambda_+
\end{pmatrix}, \quad
R=
\begin{pmatrix}
	u_1^- & u_1^+ \\
	u_2^- & u_2^+
\end{pmatrix}
=
\begin{pmatrix}
	u_2^+ & u_1^+ \\
	-u_1^+ & u_2^+
\end{pmatrix},
$$
and $u^\pm$ are Lipschitz since $v^\pm$ are Lipschitz and $\|v^\pm\|\geq 1$.
\end{proof}
\end{lemma}

We also need to control how Reifenberg flat sets change under the linear planar deformation in the previous lemma.

\begin{lemma}\label{reifenberg flat dilation}
    Let $K\geq 1$, $D =\begin{pmatrix}
    1/K & 0 \\
    0 & 1
    \end{pmatrix}$ and $\Omega$ be a $(\delta, r_0)$-Reifenberg flat domain. If $\delta <  \sqrt{15}/(16 K)$ then $D^{-1} (\Omega)$ is a $\left( \frac{8\sqrt{15}}{15} K\delta, r_0 \right)$-Reifenberg flat domain.
    \begin{proof}
        Let $0<r\leq r_0$, $\xi \in \partial \Omega$ and $\mathcal P \coloneqq \mathcal P (\xi, r) \ni \xi$. Denote $\Omega^\prime = D^{-1} (\Omega)$, $\xi^\prime =D^{-1} \xi$ and $\mathcal P^\prime = D^{-1} (\mathcal P)$. We want to check the conditions in \cref{def:Reifenberg flat} with the point $\xi'$, radius $r$ and the hyperplane $\mathcal P '$. See \Cref{Rf deformation both}.

\begin{figure}%[h]
\centering
\begin{subfigure}[b]{.4\textwidth}
  \centering
\begin{tikzpicture}[scale=2]
\clip(-1.5,-1.5) rectangle (1.5,1.5); %picture region
\begin{scope}
    \clip (0,0) circle (1);
    \fill[blue!15] (0.08097604617428178,-2.0178807893297317) -- (-1.6114797824488898,-2.324464005046858) -- (-2.324464005046858,1.6114797824488898) -- (-0.6320081764236866,1.918062998166016) -- cycle;  
\end{scope}
\begin{scope}
    \clip (0,0) circle (1);
    \fill[orange!15] (0.6320081764236866,-1.918062998166016) -- (2.324464005046858,-1.6114797824488898) -- (1.6114797824488898,2.324464005046858) -- (-0.08097604617428178,2.0178807893297317) -- cycle; 
\end{scope}
\draw  (0,0) circle (1);
\draw plot(\x,{(-0-0.9839859468739369*\x)/0.17824605564949209});
\draw [dashed] plot(\x,{(--0.14-0.9839859468739369*\x)/0.17824605564949209});
\draw [dotted] plot(\x,{(-0.28-0.9839859468739369*\x)/0.17824605564949209});
\draw [dashed] plot(\x,{(-0.14-0.9839859468739369*\x)/0.17824605564949209});
\draw [dotted] plot(\x,{(--0.28-0.9839859468739369*\x)/0.17824605564949209});
\draw (0,0) node[circle,fill,inner sep=1]{} node[above]{$\xi\in\mathcal P\cap\partial\Omega$};
\draw (0.18,-0.98) node[below right]{$\mathcal P$};
\draw (-0.1,1.15) node[above]{$\delta r$};
\draw (-0.6,1.15) node{$2\delta r$};
\draw (0.7,0) node[right]{$\Omega$};
\draw (-1,-0.3) node[right]{$\R^2\setminus\overline \Omega$};
%--
\draw (-0.6030449590854325,2.9268808399279753)--(-0.6010703928913901,2.9147925659972023)--(-0.5990958266973476,2.9027042920664297)--(-0.5976550304415114,2.8804988091439157)--(-0.596214234185675,2.8582933262214016)--(-0.5953072078680449,2.8259706343071467)--(-0.5936163353117311,2.7973877456343645)--(-0.5919254627554171,2.768804856961582)--(-0.5894507439604192,2.74396177153027)--(-0.586487446064777,2.7218487566409664)--(-0.5835241481691348,2.6997357417516628)--(-0.5800722711728483,2.6803527974043657)--(-0.5764695782976975,2.6620772174524756)--(-0.5728668854225467,2.6438016375005855)--(-0.5691133766685317,2.626633421944102)--(-0.5654264599445209,2.6112313615025777)--(-0.56173954322051,2.5958293010610536)--(-0.5581192185265035,2.582193395734489)--(-0.5544127447793954,2.5664778273147455)--(-0.5507062710322874,2.5507622588950025)--(-0.5469136482320781,2.532967027382082)--(-0.5431285992462677,2.5139479657536565)--(-0.5393435502604572,2.4949289041252305)--(-0.5355660750890453,2.4746860123812984)--(-0.5323650753913203,2.4564456767141407)--(-0.5291640756935951,2.438205341046984)--(-0.5265395514695568,2.4219675614566025)--(-0.524555172624017,2.4066401730965157)--(-0.5225707937784773,2.391312784736429)--(-0.5212265603114362,2.3768957876066383)--(-0.5205529561846096,2.362750704102552)--(-0.5198793520577829,2.348605620598465)--(-0.5198763772711705,2.3347324507200815)--(-0.5205455745584112,2.323320583807128)--(-0.5212147718456519,2.311908716894174)--(-0.5225561412067457,2.30295815294665)--(-0.5224556366143817,2.2887261252829316)--(-0.5223551320220177,2.274494097619213)--(-0.5208127534761956,2.2549806062392985)--(-0.5192128323198612,2.237585168356536)--(-0.5176129111635268,2.220189730473773)--(-0.5159554473966801,2.204912346088162)--(-0.5141593690300647,2.190735678028266)--(-0.5123632906634491,2.1765590099683703)--(-0.5104285976970646,2.163483058234189)--(-0.5077993572055018,2.1481879636921066)--(-0.5051701167139392,2.1328928691500244)--(-0.501846328697198,2.1153786318000383)--(-0.49856040000851826,2.0992510937354747)--(-0.4952744713198385,2.083123555670911)--(-0.49202640195922004,2.06838271689177)--(-0.4887022229205864,2.0545931197493768)--(-0.4853780438819527,2.0408035226069843)--(-0.48197775516530394,2.0279651671013394)--(-0.47779278416686755,2.010328833545725)--(-0.47360781316843115,1.9926924999901094)--(-0.4686381598882074,1.9702581883845247)--(-0.4642315428016779,1.9481503416563655)--(-0.4598249257151484,1.9260424949282067)--(-0.4559813448223132,1.9042611130774736)--(-0.4528954639786779,1.886788297965508)--(-0.4498095831350427,1.8693154828535425)--(-0.4474814023406071,1.8561512344803432)--(-0.44521876680961947,1.8374668018980265)--(-0.44295613127863176,1.81878236931571)--(-0.44075904101109176,1.7945777525242745)--(-0.4383444830377406,1.7705200266503445)--(-0.4359299250643893,1.746462300776414)--(-0.4332978993852267,1.722551465819988)--(-0.43039544961126874,1.7007191149973808)--(-0.42749299983731076,1.6788867641747733)--(-0.4243201259685571,1.6591328974859825)--(-0.4208895035316772,1.6398385680785834) --(-0.41745888109479723,1.6205442386711844) --(-0.4137705100897908,1.6017094465451767) --(-0.41013442292057845,1.5704443358542208) --(-0.40649833575136607,1.5391792251632648) --(-0.4029145324179477,1.4954837959073606) --(-0.3986235988117951,1.4621620566775397) --(-0.39433266520564253,1.428840317447719) --(-0.38933460132675546,1.405892268243981) --(-0.38761235738340066,1.3960663061699194) --(-0.3858901134400458,1.386240344095858) --(-0.38744368943222307,1.3895364691514729) --(-0.3755525896807531,1.353738065500147) -- (-0.32475590026080275,1.2025715870119151) --(-0.31118608608743986,1.1620964445332773) --(-0.3093822332358315,1.15603909193427) --(-0.30657363235115964,1.1403155327898398) --(-0.3037650314664877,1.1245919736454097) --(-0.2999516825487521,1.099202207955556) --(-0.29754800755023475,1.0782300749908276) --(-0.29514433255171746,1.0572579420260992) --(-0.29415033147241837,1.0407034417864958) --(-0.2920701112736056,1.0249533114740514) --(-0.28998989107479295,1.009203181161607) --(-0.28682345175646673,0.9942574207763215) --(-0.2833484033926781,0.9784634300717531) --(-0.2798733550288895,0.9626694393671845) --(-0.27608969761963836,0.9460272183433321) --(-0.2709274430247548,0.9256821754491269) --(-0.2657651884298713,0.9053371325549218) --(-0.25922433664935507,0.881289267790363) --(-0.252538681259664,0.8584532928037435) --(-0.24585302586997287,0.835617317817124) --(-0.23902256687110668,0.8139932326084434) --(-0.23243481224064663,0.7893907570429841) --(-0.22584705761018659,0.7647882814775246) --(-0.2195020073481327,0.7372074155552866) --(-0.21478465079361134,0.7087555319236009) --(-0.21006729423908999,0.6803036482919151) --(-0.2069776313921011,0.650980746950781) --(-0.20189569827184983,0.6227354060056985) --(-0.19681376515159854,0.5944900650606161) --(-0.18973956175808493,0.567322284511585) --(-0.18469148771403487,0.5498243877585065) --(-0.1796434136699848,0.5323264910054281) --(-0.17662146897539827,0.5244984780483022) --(-0.17193813546679604,0.5165954932102697) --(-0.1672548019581938,0.5086925083722369) --(-0.1609100796355759,0.5007145516532976) --(-0.15208858058234292,0.48616630341282213) --(-0.1432670815291099,0.4716180551723466) --(-0.1319688057452618,0.4504995154103349) --(-0.1258886470396078,0.4390518447815668) --(-0.11980848833395383,0.4276041741527987) --(-0.1189464467064939,0.425827372657274) --(-0.11487173231595775,0.41551870543272895) --(-0.1107970179254216,0.40521003820818385) --(-0.10350963077180916,0.3863695052546184) --(-0.09727348932966873,0.36408180964968684) --(-0.09103734788752828,0.3417941140447553) --(-0.08585245215685973,0.31605925578845734) --(-0.08148690274099513,0.2904668322828608) --(-0.07712135332513054,0.2648744087772643) --(-0.07357515022406982,0.23942442002236908) --(-0.06750944450075766,0.21508102678804358) --(-0.06144373877744548,0.19073763355371812) --(-0.052858530431881844,0.16750083583996234)--(-0.04661050425031679,0.15341871694930578)--(-0.04036247806875174,0.13933659805864917)--(-0.03645163405118525,0.13440915799109168)--(-0.0322573699409791,0.1266484053337397)--(-0.028063105830772962,0.11888765267638768)--(-0.023585421627927167,0.1082935874292411)--(-0.01981649882811401,0.09685990694267137)--(-0.016047576028300852,0.08542622645610164)--(-0.012987414631520337,0.07315293073010873)--(-0.008977851754556442,0.052885646895535765)--(-0.004968288877592545,0.032618363060962795)--(-0.000009324520445271197,0.004357091117809775)--(0.003866787772041278,-0.019061004266560773)--(0.007742900064527828,-0.04247909965093133)--(0.010536160292353656,-0.06105401847651943)--(0.013364586426664434,-0.08015034710060176)--(0.016193012560975214,-0.0992466757246841)--(0.019056604601770926,-0.1188644141472606)--(0.025072876245770637,-0.14095817033367974)--(0.03108914788977035,-0.16305192652009887)--(0.04025809913697405,-0.18762170047036064)--(0.047332483687013455,-0.20667284133962308)--(0.05440686823705285,-0.22572398220888545)--(0.05938668608992796,-0.23925648999714866)--(0.06446889262717126,-0.2521059116395057)--(0.06955109916441454,-0.26495533328186277)--(0.074735694386026,-0.2771216687783138)--(0.07895601612058306,-0.28625043116880405)--(0.08317633785514014,-0.29537919355929426)--(0.08643238610264287,-0.3014703828438239)--(0.08935368133392918,-0.30659285298016103)--(0.0922749765652155,-0.3117153231164982)--(0.09486151878028541,-0.31586907410464277)--(0.09851146994405383,-0.3233894956093376)--(0.10216142110782223,-0.33090991711403245)--(0.10687478122028918,-0.3417970091352776)--(0.1125420909668621,-0.3589175407561555)--(0.11820940071343503,-0.3760380723770334)--(0.12483066009411393,-0.39939204359754377)--(0.13065933115897665,-0.4225268571242885)--(0.13648800222383936,-0.4456616706510332)--(0.14152408497288574,-0.4685773264840116)--(0.14589148857820017,-0.4903478745209332)--(0.15025889218351457,-0.5121184225578548)--(0.15395761664509686,-0.5327438627987191)--(0.15780912354840657,-0.5553444618980942)--(0.1616606304517163,-0.5779450609974692)--(0.1656649197967534,-0.6025208189553553)--(0.1661531966490663,-0.6251347897985335)--(0.1666414735013792,-0.6477487606417118)--(0.16361373786096786,-0.6684009443701822)--(0.16215918142138946,-0.6827290243915993)--(0.16070462498181107,-0.6970571044130163)--(0.1608232477430656,-0.7050610807273799)--(0.16147154168395933,-0.7163195133370062)--(0.16211983562485305,-0.7275779459466325)--(0.16329780074538597,-0.7420908348515215)--(0.16232914709960714,-0.7863074389127157) -- (0.1556677215120775,-0.9609718976704773) -- (0.14970173207813395,-1.0986094467786842) -- (0.1475667688595369,-1.21373359985481)--(0.14743836876487215,-1.2420667190060737)--(0.14783936844242113,-1.2539171852722393)--(0.1514754381987838,-1.2758313313331533)--(0.1551115079551465,-1.2977454773940673)--(0.1619826477903229,-1.3297233032497289)--(0.16649411510228945,-1.3478069325151383)--(0.171005582414256,-1.3658905617805475)--(0.17315737720301277,-1.370079994455704)--(0.17809900480440022,-1.3772092570872163) -- (0.1952737677714423,-1.4002163603048299)--(0.19977544272307868,-1.406025108303063)--(0.20104733221233306,-1.4075735117116612)--(0.2057319329565159,-1.418850644020651)--(0.21041653370069874,-1.4301277763296407)--(0.21851384569980997,-1.4511336375390211)--(0.22417092950789466,-1.4721169754770056)--(0.22982801331597935,-1.4931003134149898)--(0.23304486893303722,-1.5140611280815768)--(0.2351795293451019,-1.5314041641671792)--(0.2373141897571665,-1.5487472002527813)--(0.23836665496423753,-1.5624724577573974)--(0.239409905671156,-1.579215773931963)--(0.2404531563780745,-1.5959590901065284)--(0.24148719258484053,-1.615720464951043)--(0.24308673663775532,-1.6385099331303445)--(0.24468628069067003,-1.6612994013096463)--(0.24685133258973355,-1.6871169628237341)--(0.24915925384472343,-1.7089556523809215)--(0.2514671750997133,-1.7307943419381087)--(0.2539179657106294,-1.7486541595383946)--(0.2563238811951949,-1.7645691838631865)--(0.25872979667976037,-1.7804842081879781)--(0.26109083703797525,-1.794454439237276)--(0.26380828494310593,-1.8082343751085965)--(0.26652573284823666,-1.8220143109799174)--(0.2695995883002832,-1.835603951673261)--(0.2726650598883756,-1.849618416140483)--(0.275730531476468,-1.863632880607705)--(0.2787876192006062,-1.878072168848805)--(0.28266452536377157,-1.8983158948724086)--(0.28654143152693695,-1.9185596208960123)--(0.2912381561291293,-1.9446077847021186)--(0.2946047447482473,-1.9669873054812683)--(0.29797133336736537,-1.989366826260418)--(0.30000778600340927,-2.008077704012612)--(0.3034021948471276,-2.028504331626212)--(0.30679660369084594,-2.0489309592398124)--(0.3115489687422384,-2.071073336714817)--(0.3152101615895042,-2.088866247049523)--(0.31887135443677,-2.1066591573842293)--(0.32144137507990883,-2.120102600578636)--(0.323979798404189,-2.1340262891918695)--(0.3265182217284691,-2.1479499778051037)--(0.32902504773389063,-2.162353911837166)--(0.33215729293322505,-2.1773626742190917)--(0.3352895381325594,-2.1923714366010163)--(0.33904720252580645,-2.207985027332802)--(0.34339123825771595,-2.22563327575301)--(0.3477352739896255,-2.243281524173217)--(0.35266568106019747,-2.2629644302818455)--(0.357620117264222,-2.284946885491233)--(0.36257455346824663,-2.30692934070062)--(0.36755301880572366,-2.331211345010765)--(0.37216278832277205,-2.354569654457447)--(0.37677255783982033,-2.3779279639041286)--(0.3810136315364398,-2.4003625784873464)--(0.3853132293041256,-2.423834856853426)--(0.3896128270718114,-2.4473071352195053)--(0.39397094891056295,-2.471817077368444)--(0.3984516423788563,-2.4943217676800304)--(0.4029323358471496,-2.5168264579916175)--(0.4075356009449846,-2.5373258964658527)--(0.41058444651583836,-2.549949976286579)--(0.41363329208669214,-2.5625740561073043)--(0.4151277181305647,-2.5673227772745206)--(0.4228093905418595,-2.594077590186728) -- (0.4519375166797628,-2.6965262346289776)--(0.4595150516276543,-2.723459161501821)--(0.46080120274472014,-2.7285641105903093)--(0.4635075718574687,-2.7397562049716093)--(0.4662139409702173,-2.7509482993529097)--(0.4703405280786481,-2.7682275390270203)--(0.4746490514228066,-2.788246811331348)--(0.478957574766965,-2.8082660836356754)--(0.4834480343468509,-2.8310253885702195)--(0.48640952374016466,-2.846076643533358) -- (0.492236051546962,-2.8758143084799648);
%--
\end{tikzpicture}
  \caption{$\Omega$, $0<r\leq r_0$ and $\xi\in\mathcal P\cap \partial\Omega$.}
  \label{Rf deformation 1}
\end{subfigure}%
\hfill
\begin{subfigure}[b]{.6\textwidth}
  \centering
\begin{tikzpicture}[scale=2]
%Geogebra with K=2, alpha=1.75, delta=0.14
\clip(-2.1,-1.5) rectangle (2.1,1.5); %picture region
\begin{scope}
    \clip (0,0) ellipse (2 and 1);
    \fill[blue!15] (0.16195209234856356,-2.0178807893297317) -- (-3.2229595648977796,-2.324464005046858) -- (-4.648928010093716,1.6114797824488898) -- (-1.2640163528473731,1.918062998166016) -- cycle;  
\end{scope}
\begin{scope}
    \clip (0,0) circle (2 and 1);
    \fill[orange!15] (1.2640163528473731,-1.918062998166016) -- (4.648928010093716,-1.6114797824488898) -- (3.2229595648977796,2.324464005046858) -- (-0.16195209234856356,2.0178807893297317) -- cycle; 
\end{scope}
\draw  (0,0) ellipse (2 and 1);
\draw (2*0.18,-0.98) node[below right]{$\mathcal P^\prime$};
\draw (-1.3,1.15) node{$2K\delta r$};
\draw (2*0.7,0) node[right]{$\Omega^\prime$};
\draw (-1.2,0) node[left]{$\R^2\setminus\overline {\Omega^\prime}$};
%--
\clip(-3.793211017391042,-3.0945705637773004) rectangle (4.638462463849253,3.281882506409062);
\draw (0,0) circle (1);
\draw plot(\x,{(-0--0.17539161384480992*\x)/-0.06354331270920366});
\draw [dashed] plot(\x,{(-0.0499088955818578--0.17539161384480992*\x)/-0.06354331270920366});
\draw [dashed] plot(\x,{(--0.0499088955818578--0.17539161384480992*\x)/-0.06354331270920366});
\draw [densely dotted] plot(\x,{(--0.0998177911637156--0.17539161384480992*\x)/-0.06354331270920366});
\draw [densely dotted] plot(\x,{(-0.0998177911637156--0.17539161384480992*\x)/-0.06354331270920366});
\draw (1.2640163528473731,-1.918062998166016)-- (4.648928010093716,-1.6114797824488898);
\draw [line width=1.2, red!50 ] plot(\x,{(--0.10446660156234223--0.17539161384480992*\x)/-0.06354331270920366});
\draw [line width=1.2, red!50 ] plot(\x,{(-0.10446660156234226--0.17539161384480992*\x)/-0.06354331270920366});
\draw [dash dot] plot(\x,{(--0.10446660156234223--0.17539161384480992*\x)/-0.06354331270920366});
\draw [dash dot] plot(\x,{(-0.10446660156234226--0.17539161384480992*\x)/-0.06354331270920366});
%--
\draw (0,0) node[circle,fill,inner sep=1]{} node[above]{$\xi^\prime\in\mathcal P^\prime\cap\partial\Omega^\prime$};
%--
\draw (-1.206089918170865,2.9268808399279753)--(-1.2021407857827802,2.9147925659972023)--(-1.1981916533946952,2.9027042920664297)--(-1.1953100608830227,2.8804988091439157)--(-1.19242846837135,2.8582933262214016)--(-1.1906144157360898,2.8259706343071467)--(-1.1872326706234622,2.7973877456343645)--(-1.1838509255108343,2.768804856961582)--(-1.1789014879208384,2.74396177153027)--(-1.172974892129554,2.7218487566409664)--(-1.1670482963382696,2.6997357417516628)--(-1.1601445423456966,2.6803527974043657)--(-1.152939156595395,2.6620772174524756)--(-1.1457337708450934,2.6438016375005855)--(-1.1382267533370634,2.626633421944102)--(-1.1308529198890418,2.6112313615025777)--(-1.12347908644102,2.5958293010610536)--(-1.116238437053007,2.582193395734489)--(-1.1088254895587908,2.5664778273147455)--(-1.1014125420645748,2.5507622588950025)--(-1.0938272964641562,2.532967027382082)--(-1.0862571984925353,2.5139479657536565)--(-1.0786871005209144,2.4949289041252305)--(-1.0711321501780906,2.4746860123812984)--(-1.0647301507826405,2.4564456767141407)--(-1.0583281513871903,2.438205341046984)--(-1.0530791029391136,2.4219675614566025)--(-1.049110345248034,2.4066401730965157)--(-1.0451415875569545,2.391312784736429)--(-1.0424531206228724,2.3768957876066383)--(-1.0411059123692192,2.362750704102552)--(-1.0397587041155658,2.348605620598465)--(-1.039752754542341,2.3347324507200815)--(-1.0410911491168224,2.323320583807128)--(-1.0424295436913038,2.311908716894174)--(-1.0451122824134913,2.30295815294665)--(-1.0449112732287633,2.2887261252829316)--(-1.0447102640440353,2.274494097619213)--(-1.0416255069523912,2.2549806062392985)--(-1.0384256646397223,2.237585168356536)--(-1.0352258223270536,2.220189730473773)--(-1.0319108947933602,2.204912346088162)--(-1.0283187380601293,2.190735678028266)--(-1.0247265813268982,2.1765590099683703)--(-1.0208571953941292,2.163483058234189)--(-1.0155987144110037,2.1481879636921066)--(-1.0103402334278784,2.1328928691500244)--(-1.003692657394396,2.1153786318000383)--(-0.9971208000170365,2.0992510937354747)--(-0.990548942639677,2.083123555670911)--(-0.9840528039184401,2.06838271689177)--(-0.9774044458411728,2.0545931197493768)--(-0.9707560877639054,2.0408035226069843)--(-0.9639555103306079,2.0279651671013394)--(-0.9555855683337351,2.010328833545725)--(-0.9472156263368623,1.9926924999901094)--(-0.9372763197764148,1.9702581883845247)--(-0.9284630856033558,1.9481503416563655)--(-0.9196498514302968,1.9260424949282067)--(-0.9119626896446263,1.9042611130774736)--(-0.9057909279573558,1.886788297965508)--(-0.8996191662700854,1.8693154828535425)--(-0.8949628046812143,1.8561512344803432)--(-0.8904375336192389,1.8374668018980265)--(-0.8859122625572635,1.81878236931571)--(-0.8815180820221835,1.7945777525242745)--(-0.8766889660754812,1.7705200266503445)--(-0.8718598501287786,1.746462300776414)--(-0.8665957987704535,1.722551465819988)--(-0.8607908992225375,1.7007191149973808)--(-0.8549859996746215,1.6788867641747733)--(-0.8486402519371142,1.6591328974859825)--(-0.8417790070633544,1.6398385680785834)--(-0.8349177621895945,1.6205442386711844)--(-0.8275410201795816,1.6017094465451767)--(-0.8202688458411569,1.5704443358542208)--(-0.8129966715027321,1.5391792251632648)--(-0.8058290648358954,1.4954837959073606)--(-0.7972471976235902,1.4621620566775397)--(-0.7886653304112851,1.428840317447719)--(-0.7786692026535109,1.405892268243981)--(-0.7752247147668013,1.3960663061699194)--(-0.7717802268800916,1.386240344095858)--(-0.7748873788644461,1.3895364691514729)--(-0.7511051793615062,1.353738065500147) -- (-0.6495118005216055,1.2025715870119151)--(-0.6223721721748797,1.1620964445332773)--(-0.618764466471663,1.15603909193427)--(-0.6131472647023193,1.1403155327898398)--(-0.6075300629329754,1.1245919736454097)--(-0.5999033650975042,1.099202207955556)--(-0.5950960151004695,1.0782300749908276)--(-0.5902886651034349,1.0572579420260992)--(-0.5883006629448367,1.0407034417864958)--(-0.5841402225472112,1.0249533114740514)--(-0.5799797821495859,1.009203181161607)--(-0.5736469035129335,0.9942574207763215)--(-0.5666968067853562,0.9784634300717531)--(-0.559746710057779,0.9626694393671845)--(-0.5521793952392767,0.9460272183433321)--(-0.5418548860495096,0.9256821754491269)--(-0.5315303768597426,0.9053371325549218)--(-0.5184486732987101,0.881289267790363)--(-0.505077362519328,0.8584532928037435)--(-0.49170605173994575,0.835617317817124)--(-0.47804513374221336,0.8139932326084434)--(-0.46486962448129326,0.7893907570429841)--(-0.45169411522037317,0.7647882814775246)--(-0.4390040146962654,0.7372074155552866)--(-0.4295693015872227,0.7087555319236009)--(-0.42013458847817997,0.6803036482919151)--(-0.4139552627842022,0.650980746950781)--(-0.40379139654369967,0.6227354060056985)--(-0.3936275303031971,0.5944900650606161)--(-0.37947912351616986,0.567322284511585)--(-0.36938297542806975,0.5498243877585065)--(-0.3592868273399696,0.5323264910054281)--(-0.35324293795079653,0.5244984780483022)--(-0.3438762709335921,0.5165954932102697)--(-0.3345096039163876,0.5086925083722369)--(-0.3218201592711518,0.5007145516532976)--(-0.30417716116468585,0.48616630341282213)--(-0.2865341630582198,0.4716180551723466)--(-0.2639376114905236,0.4504995154103349)--(-0.2517772940792156,0.4390518447815668)--(-0.23961697666790766,0.4276041741527987)--(-0.2378928934129878,0.425827372657274)--(-0.2297434646319155,0.41551870543272895)--(-0.2215940358508432,0.40521003820818385)--(-0.20701926154361833,0.3863695052546184)--(-0.19454697865933746,0.36408180964968684)--(-0.18207469577505656,0.3417941140447553)--(-0.17170490431371946,0.31605925578845734)--(-0.16297380548199025,0.2904668322828608)--(-0.15424270665026107,0.2648744087772643)--(-0.14715030044813965,0.23942442002236908)--(-0.13501888900151532,0.21508102678804358)--(-0.12288747755489096,0.19073763355371812)--(-0.10571706086376369,0.16750083583996234)--(-0.09322100850063358,0.15341871694930578)--(-0.08072495613750348,0.13933659805864917)--(-0.0729032681023705,0.13440915799109168)--(-0.0645147398819582,0.1266484053337397)--(-0.056126211661545924,0.11888765267638768)--(-0.047170843255854335,0.1082935874292411)--(-0.03963299765622802,0.09685990694267137)--(-0.032095152056601704,0.08542622645610164)--(-0.025974829263040674,0.07315293073010873)--(-0.017955703509112884,0.052885646895535765)--(-0.00993657775518509,0.032618363060962795)--(-0.000018649040890542394,0.004357091117809775)--(0.007733575544082556,-0.019061004266560773)--(0.015485800129055656,-0.04247909965093133)--(0.021072320584707312,-0.06105401847651943)--(0.02672917285332887,-0.08015034710060176)--(0.03238602512195043,-0.0992466757246841)--(0.03811320920354185,-0.1188644141472606)--(0.050145752491541275,-0.14095817033367974)--(0.0621782957795407,-0.16305192652009887)--(0.0805161982739481,-0.18762170047036064)--(0.09466496737402691,-0.20667284133962308)--(0.1088137364741057,-0.22572398220888545)--(0.11877337217985592,-0.23925648999714866)--(0.12893778525434252,-0.2521059116395057)--(0.13910219832882909,-0.26495533328186277)--(0.149471388772052,-0.2771216687783138)--(0.15791203224116612,-0.28625043116880405)--(0.16635267571028028,-0.29537919355929426)--(0.17286477220528573,-0.3014703828438239)--(0.17870736266785836,-0.30659285298016103)--(0.184549953130431,-0.3117153231164982)--(0.18972303756057082,-0.31586907410464277)--(0.19702293988810765,-0.3233894956093376)--(0.20432284221564445,-0.33090991711403245)--(0.21374956244057836,-0.3417970091352776)--(0.2250841819337242,-0.3589175407561555)--(0.23641880142687005,-0.3760380723770334)--(0.24966132018822787,-0.39939204359754377)--(0.2613186623179533,-0.4225268571242885)--(0.2729760044476787,-0.4456616706510332)--(0.2830481699457715,-0.4685773264840116)--(0.29178297715640034,-0.4903478745209332)--(0.30051778436702914,-0.5121184225578548)--(0.3079152332901937,-0.5327438627987191)--(0.31561824709681313,-0.5553444618980942)--(0.3233212609034326,-0.5779450609974692)--(0.3313298395935068,-0.6025208189553553)--(0.3323063932981326,-0.6251347897985335)--(0.3332829470027584,-0.6477487606417118)--(0.3272274757219357,-0.6684009443701822)--(0.32431836284277893,-0.6827290243915993)--(0.32140924996362213,-0.6970571044130163)--(0.3216464954861312,-0.7050610807273799)--(0.32294308336791866,-0.7163195133370062)--(0.32659560149077194,-0.7420908348515215)--(0.32465829419921427,-0.7863074389127157) --(0.311335443024155,-0.9609718976704773) -- (0.2994034641562679,-1.0986094467786842) -- (0.2951335377190738,-1.21373359985481)--(0.2948767375297443,-1.2420667190060737)--(0.29567873688484225,-1.2539171852722393)--(0.3029508763975676,-1.2758313313331533)--(0.310223015910293,-1.2977454773940673)--(0.3239652955806458,-1.3297233032497289)--(0.3329882302045789,-1.3478069325151383)--(0.342011164828512,-1.3658905617805475)--(0.34631475440602555,-1.370079994455704)--(0.35619800960880044,-1.3772092570872163) --(0.3905475355428846,-1.4002163603048299)--(0.39955088544615736,-1.406025108303063)--(0.4020946644246661,-1.4075735117116612)--(0.4114638659130318,-1.418850644020651)--(0.4208330674013975,-1.4301277763296407)--(0.43702769139961994,-1.4511336375390211)--(0.4483418590157893,-1.4721169754770056)--(0.4596560266319587,-1.4931003134149898)--(0.46608973786607444,-1.5140611280815768)--(0.4703590586902038,-1.5314041641671792)--(0.474628379514333,-1.5487472002527813)--(0.47673330992847507,-1.5624724577573974)--(0.478819811342312,-1.579215773931963)--(0.480906312756149,-1.5959590901065284)--(0.48297438516968105,-1.615720464951043)--(0.48617347327551064,-1.6385099331303445)--(0.48937256138134005,-1.6612994013096463)--(0.4937026651794671,-1.6871169628237341)--(0.49831850768944685,-1.7089556523809215)--(0.5029343501994266,-1.7307943419381087)--(0.5078359314212588,-1.7486541595383946)--(0.5126477623903898,-1.7645691838631865)--(0.5174595933595207,-1.7804842081879781)--(0.5221816740759505,-1.794454439237276)--(0.5276165698862119,-1.8082343751085965)--(0.5330514656964733,-1.8220143109799174)--(0.5391991766005664,-1.835603951673261)--(0.5453301197767512,-1.849618416140483)--(0.551461062952936,-1.863632880607705)--(0.5575752384012124,-1.878072168848805)--(0.5653290507275431,-1.8983158948724086)--(0.5730828630538739,-1.9185596208960123)--(0.5824763122582586,-1.9446077847021186)--(0.5892094894964947,-1.9669873054812683)--(0.5959426667347307,-1.989366826260418)--(0.6000155720068185,-2.008077704012612)--(0.6068043896942552,-2.028504331626212)--(0.6135932073816919,-2.0489309592398124)--(0.6230979374844768,-2.071073336714817)--(0.6304203231790084,-2.088866247049523)--(0.63774270887354,-2.1066591573842293)--(0.6428827501598177,-2.120102600578636)--(0.647959596808378,-2.1340262891918695)--(0.6530364434569382,-2.1479499778051037)--(0.6580500954677813,-2.162353911837166)--(0.6643145858664501,-2.1773626742190917)--(0.6705790762651188,-2.1923714366010163)--(0.6780944050516129,-2.207985027332802)--(0.6867824765154319,-2.22563327575301)--(0.695470547979251,-2.243281524173217)--(0.7053313621203949,-2.2629644302818455)--(0.715240234528444,-2.284946885491233)--(0.7251491069364933,-2.30692934070062)--(0.7351060376114473,-2.331211345010765)--(0.7443255766455441,-2.354569654457447)--(0.7535451156796407,-2.3779279639041286)--(0.7620272630728796,-2.4003625784873464)--(0.7706264586082512,-2.423834856853426)--(0.7792256541436228,-2.4473071352195053)--(0.7879418978211259,-2.471817077368444)--(0.7969032847577125,-2.4943217676800304)--(0.8058646716942992,-2.5168264579916175)--(0.8150712018899692,-2.5373258964658527)--(0.8211688930316767,-2.549949976286579)--(0.8272665841733843,-2.5625740561073043)--(0.8302554362611294,-2.5673227772745206)--(0.845618781083719,-2.594077590186728) --(0.9038750333595256,-2.6965262346289776)--(0.9190301032553087,-2.723459161501821)--(0.9216024054894403,-2.7285641105903093)--(0.9270151437149374,-2.7397562049716093)--(0.9324278819404346,-2.7509482993529097)--(0.9406810561572962,-2.7682275390270203)--(0.9492981028456132,-2.788246811331348)--(0.95791514953393,-2.8082660836356754)--(0.9668960686937018,-2.8310253885702195)--(0.9728190474803293,-2.846076643533358) --(0.984472103093924,-2.8758143084799648);
%--
\end{tikzpicture}
%Geogebra with K=2, alpha=1.75, delta=0.14
  \caption{$D^{-1}$ of \cref{Rf deformation 1} and $\{x:\dist(x,\mathcal P^\prime)=2K\delta r\}$ in red.}
  \label{Rf deformation 2}
\end{subfigure}%
\caption{Almost the worst situation: $\mathcal P\perp \{y=0\}$.}
\label{Rf deformation both}
\end{figure}
        
        \begin{claim}\label{relation distances under linear transformation}
            For any $x,y\in \R^2$, $K^{-1} \dist(x,y) \leq \dist (Dx, Dy) \leq \dist (x,y)$.
            \begin{proof}
                Indeed, since the minimum (resp. maximum) eigenvalue of $D$ is $K^{-1}$ (resp. 1), we have
                $$
                \frac{\dist(Dx, Dy)}{\dist(x,y)} = \frac{\|Dx-Dy\|}{\|x-y\|} \in [K^{-1}, 1] .
                $$
            \end{proof}
        \end{claim}
    
        \begin{claim}\label{condition 2 dilation reifenberg flat}
            One component of
            \begin{equation}\label{condition 2 dilation reifenberg flat 2}
                B_r (\xi') \cap \{ x\in \R^2 : \dist(x, \mathcal P ' ) \geq 2K \delta r \}
            \end{equation}
            is contained in $\Omega^\prime$ and the other is contained in $\R^2 \setminus \Omega^\prime$.
            \begin{proof}
                By the previous claim we get
                \begin{equation}\label{one component dilation reifenberg flat}
                    D^{-1} \left( \{x\in \R^2 : \dist (x,\mathcal P) \geq 2\delta r\} \right) 
                    \supset \{x\in \R^2 : \dist (x,\mathcal P^\prime) \geq 2 K \delta r\} ,
                \end{equation}
                and from the definition of $D^{-1}$ we have $B_r (\xi' ) \subset D^{-1} \left( B_r (\xi) \right)$. By \rf{one component dilation reifenberg flat}, and since $2K\delta < 2\sqrt{15}/16 < 1$, in particular $B_r (\xi^\prime ) \cap D^{-1} \left( \{x : \dist (x,\mathcal P) \geq 2\delta r\} \right) \not = \emptyset$, and the claim follows from the Reifenberg flat condition of $\Omega$.
            \end{proof}
        \end{claim}
    
        First we check
        \begin{equation}\label{condition 1.1 dilation reifenberg flat}
            \sup_{y \in  \partial \Omega'\cap B_r (\xi')}\dist(y, \mathcal P' \cap B_r (\xi ')) \leq K \delta r .
        \end{equation}
        From the Reifenberg flatness of $\Omega$, see \cref{def:Reifenberg flat}\rf{def:Reifenberg flat 1}, we have
        $$
        \sup_{x\in \partial \Omega \cap B_r (\xi)} \dist (x, \mathcal P \cap B(\xi, r)) \leq \delta r.
        $$
        Given $y \in{\partial \Omega}^\prime \cap B_r (\xi^\prime)$, take $z \in \mathcal P \cap B_r (\xi)$ with $\dist(D y, z)\leq \delta r$. Since $\mathcal P^\prime$ is a line through $\xi$, the distance $\dist(y,\mathcal P^\prime)$ is attained at $B_r (\xi^\prime)$, that is,
        $$
        \dist (y, \mathcal P^\prime \cap B_r (\xi^\prime)) 
        = \dist (y, \mathcal P^\prime) 
        \leq \dist (y, D^{-1} z).
        $$
        By \cref{relation distances under linear transformation},
        $$
        \dist (y, \mathcal P^\prime \cap B_r (\xi^\prime))
        \leq K \dist (D y, z) \leq K\delta r,
        $$
        and \rf{condition 1.1 dilation reifenberg flat} follows.
        
        Now we turn to prove
        \begin{equation}\label{condition 1.2 dilation reifenberg flat}
            \sup_{y \in \mathcal P'\cap B_r (\xi')}\dist(y, \partial \Omega' \cap B_r (\xi ')) \leq \frac{2\sqrt{2} K \delta}{\sqrt{1+\sqrt{1-4K^2 \delta^2}}} r.
        \end{equation}
%beginning of picture of {condition 1.2 dilation reifenberg flat}
		\begin{figure}%[h]
			\centering
			\begin{tikzpicture}[scale=3]
                    \def \dosKDelta {0.7}
                    \def \ellipseDilatation {2.5} %ellipse 1 and '3'
                    \draw [line width=1.5, red!50 ] (0.6,0) -- (0.6,-0.8*\dosKDelta);
                    \draw [line width=1.5, red!50 ] (0,0) -- (0.6,0);
                    \filldraw[color=black, fill=red!10 ] (-{(1-\dosKDelta^2)^0.5},0) -- (-{(1-\dosKDelta^2)^0.5},\dosKDelta) -- node[left]{\scriptsize$\frac{2\sqrt{2} K \delta}{\sqrt{1+\sqrt{1-4K^2 \delta^2}}} r$} (-1,0) -- node[below] {\scriptsize$r(1-\sqrt{1-4K^2\delta^2})$} cycle;
                    % \draw (-1,\dosKDelta*0.5) node[left]{\tiny$\frac{2\sqrt{2} K \delta}{\sqrt{1+\sqrt{1-4K^2 \delta^2}}} r$};
                    % \draw ({0.5*(-1-(1-\dosKDelta^2)^0.5)},0) node[below] {\tiny$r(1-\sqrt{1-4K^2\delta^2})$};
                    \def \pointPx {{0.7*(-1)+(1-0.7)*(-(1-\dosKDelta^2)^0.5)}}
                    \draw (\pointPx,0) node[circle,fill,inner sep=1]{} node[above right]{$y$};
                    \draw[dashed] (\pointPx,0) -- (-{(1-\dosKDelta^2)^0.5},\dosKDelta);
                    \clip(-1.8,-1.02) rectangle (1.8,1.02); %picture region
                    \begin{scope}
                        \clip (0,0) ellipse (1 and \ellipseDilatation);
                        \fill[blue!15] (-10,-10) rectangle (10,-\dosKDelta);  
                    \end{scope}
                    \begin{scope}
                        \clip (0,0) ellipse (1 and \ellipseDilatation); 
                        \fill[orange!15] (-10,\dosKDelta) rectangle (10,10);  
                    \end{scope}
                    %default value of line width is 0.4
				\draw[line width=0.2] (0,0) circle (1);
                    \draw[dotted] (0,0) ellipse (1 and \ellipseDilatation);
				\draw (-0.5,-0.8) node[below left] {$B_r(\xi^\prime)$};
				\draw (0.6,-\dosKDelta) node[below right] {$\R^2 \setminus \overline\Omega$};
				\draw (0.7,\dosKDelta) node[above right] {$\Omega$};
                    \draw plot(\x,{0});
                    \draw[densely dotted] plot(\x,{\dosKDelta});
                    \draw[densely dotted] plot(\x,{-\dosKDelta});
                    \draw (0.6,\dosKDelta*0.55) node[right] {$\ell$};
                    \draw (0.6,-10) -- (0.6,10);
                    \draw (0,0) node[circle,fill,inner sep=1]{} node[above]{$\xi^\prime\in\mathcal P^\prime\cap\partial\Omega^\prime$};
                    \draw (0.6/2,0) node[below]{$\leq r$};
                    \draw (0.6,0) node[circle,fill,inner sep=1]{} node[above right]{$y$};
                    \draw (0.6,-0.8*\dosKDelta) node[circle,fill,inner sep=1]{} node[right]{$p\in\partial\Omega^\prime$};
                    \draw (0.6,-0.5*0.8*\dosKDelta) node[right]{$\leq 2K\delta r$};
                    \draw (1.7,0) node[above] {$\mathcal P^\prime$};
                    \draw [<->] (1.2,0) -- node[right]{$2K\delta r$}(1.2,\dosKDelta);
			\end{tikzpicture}
			\caption{Setting of \rf{condition 1.2 dilation reifenberg flat}.}\label{visual proof}
		\end{figure}%end of picture of {condition 1.2 dilation reifenberg flat}
        By \cref{condition 2 dilation reifenberg flat} we have that for each line $\ell$ orthogonal to $\mathcal P^\prime$ such that
        $$
        \ell \cap B_r (\xi') \cap \{ x\in \R^2 : \dist(x, \mathcal P ' ) \geq 2K \delta r \} \not = \emptyset ,
        $$
        there is a point $p\in \ell \cap \partial \Omega'  \cap \{x \in  B_r (\xi' ) : \dist (x, \mathcal P ') < 2K\delta r\}$, since $\partial \Omega'$ must separate each component in \rf{condition 2 dilation reifenberg flat 2}. Using this fact and $\xi^\prime \in \mathcal P^\prime \cap \partial \Omega^\prime$ (the center of the ball), for every $y\in \mathcal P^\prime \cap B_r (\xi^\prime)$ we obtain
        \begin{equation*}
        \dist (y, \partial \Omega' \cap B_r (\xi')) 
        \leq  \min \left\{r, \dist \left(y, B_r (\xi') \cap \{ x : \dist(x, \mathcal P ' ) \geq 2K \delta r \} \right) \right\}.
        \end{equation*}
        Using basic trigonometric computations, see \Cref{visual proof}, we have that for every $y\in \mathcal P^\prime \cap B_r (\xi^\prime)$,
        $$
        \dist \left(y, B_r (\xi') \cap \{ x : \dist(x, \mathcal P ' ) \geq 2K \delta r \} \right) \leq \frac{2\sqrt{2} K \delta}{\sqrt{1+\sqrt{1-4K^2 \delta^2}}} r.
        $$
        Since $K\delta < \sqrt{15}/16 < \sqrt{3}/4$, the previous value is less than $r$, and hence \rf{condition 1.2 dilation reifenberg flat} is proved.
        
        Notice that $K\delta r < 2K\delta r < \frac{2\sqrt{2} K \delta}{\sqrt{1+\sqrt{1-4K^2 \delta^2}}} r$. From this, \rf{condition 1.1 dilation reifenberg flat} and \rf{condition 1.2 dilation reifenberg flat} we get
        $$
        \dist_\HH (\partial \Omega ' \cap B_r (\xi') , \mathcal P ' \cap B_r (\xi')) \leq \frac{2\sqrt{2} K \delta}{\sqrt{1+\sqrt{1-4K^2 \delta^2}}} r ,
        $$
        and \cref{def:Reifenberg flat}\rf{def:Reifenberg flat 1} is verified.
        
        Since the last term is strictly larger than $2K\delta r$, and $2 \frac{2\sqrt{2} K \delta}{\sqrt{1+\sqrt{1-4K^2 \delta^2}}} <1$ when $\delta < \frac{\sqrt{15}}{16K}$, we get
        $$
        B_r (\xi^\prime) \cap\{x : \dist (x,\mathcal P^\prime) \geq 2 K \delta r\}
        \supset
        B_r (\xi^\prime) \cap \left\{x : \dist (x,\mathcal P^\prime) \geq 2 \frac{2\sqrt{2} K \delta}{\sqrt{1+\sqrt{1-4K^2 \delta^2}}} r \right\} \not = \emptyset,
        $$
        and the second condition of Reifenberg flat, \cref{def:Reifenberg flat}\rf{def:Reifenberg flat 2}, is achieved by \cref{condition 2 dilation reifenberg flat}.
        
        In conclusion, if $\delta < \sqrt{15}/(16K)$ then $\Omega' = D^{-1} (\Omega)$ is $\left( \frac{2\sqrt{2} K \delta}{\sqrt{1+\sqrt{1-4K^2 \delta^2}}}, r_0 \right)$-Reifenberg flat. In particular $\Omega' = D^{-1} (\Omega)$ is $\left( \frac{8\sqrt{15}}{15} K\delta, r_0 \right)$-Reifenberg flat.
    \end{proof}
\end{lemma}

\subsection{Reduction 1: Lipschitz diagonalization of the symmetric part}\label{sec:reduction diagonalization}

Let us see how \cref{covering elliptic measure} implies \cref{only lambda dependence covering elliptic measure}. We will do this in two steps. First we show how to get rid of the assumption on the decomposition of $A_0$, using the linear transformation in \cref{symmetric to RBR}.

\begin{claim}\label{covering elliptic measure sym}
	In \cref{covering elliptic measure}, the hypothesis $A_0=R^T B R$ is unnecessary.
	\begin{proof}
	For every square matrix $X$ we define its symmetric part $X_0\coloneqq\frac{X+X^T}{2}$.
	
	By \cref{symmetric to RBR} there exists a constant diagonal matrix $D=\begin{psmallmatrix}
		1/K & 0 \\
		0 & 1
	\end{psmallmatrix}$ with $K=K(\lambda) \geq 1$ such that $\widetilde{(A_0)} (\cdot) = D^{-1} A_0(D\cdot) D^{-1}$ can be written as $ \widetilde{(A_0)} = R^T B R \in C^{0,1} (D^{-1} (U_{r_0} (\partial \Omega)))$ with $R\in C^{0,1} (D^{-1} (U_{r_0} (\partial \Omega)))$ a rotation, and $B\in C^{0,1} (D^{-1} (U_{r_0} (\partial \Omega)))$ diagonal.

	Setting $\widetilde A (\cdot) \coloneqq D^{-1}A(D\cdot) D^{-1}$, we have that the symmetric part of the matrix $\widetilde A$ is 
	$$
	\widetilde A_0 \coloneqq
	\frac{D^{-1}A(D\cdot) D^{-1} + \left(D^{-1}A(D\cdot) D^{-1}\right)^T}{2} 
	= \frac{D^{-1}A(D\cdot) D^{-1} + D^{-1}A^T(D\cdot) D^{-1}}{2} =\widetilde{(A_0)},
	$$
	and hence $\widetilde A_0 = R^T B R$ as before.
	
	Note that 
	$$
	U_{r_0} (\partial D^{-1}(\Omega)) \coloneqq \{x : \dist(x,\partial D^{-1}\Omega)<r_0\}
	\subset D^{-1} (U_{r_0}(\partial \Omega)),
	$$
	and so these matrices are Lipschitz in $U_{r_0} (\partial D^{-1} (\Omega))$.

	Denoting $ \widetilde \omega \coloneqq \omega_{D^{-1} \Omega, \widetilde A}$ the elliptic measure in $D^{-1} \Omega$ with matrix $\widetilde A$, by \cref{elliptic measure deformation} we have $\omega^x (\cdot) = \widetilde \omega^{D^{-1}x} (D^{-1}\cdot)$ for any $x\in \Omega$. By \cref{reifenberg flat dilation} we have that $D^{-1}\Omega$ is $(8\sqrt{15} K\delta / 15, r_0)$-Reifenberg flat. Set $\widetilde p \coloneqq D^{-1} p \in D^{-1} (\Omega)$. Since $\dist (p, \partial \Omega)>r_0$, $\dist (\widetilde p, \partial D^{-1}\Omega) > r_0$ and we are in position to apply \cref{covering elliptic measure} with this pole $\widetilde p$.
	
	First, we need to compute the ellipticity constant, the Lipschitz seminorm, and the $L^\infty$ norm of the matrix $\widetilde A$. Recall $\lambda \geq 1$ is the ellipticity constant of $A$. For $\xi, \eta \in \R^2$,
	\begin{align*}
	\langle \widetilde A \xi, \eta \rangle &= \langle A(D\cdot) D^{-1} \xi, D^{-1} \eta \rangle \leq \lambda |D^{-1} \xi| |D^{-1} \eta| \leq \lambda K^2 |\xi| |\eta|, \\
	\langle \widetilde A \xi, \xi \rangle &= \langle A(D\cdot) D^{-1} \xi, D^{-1} \xi \rangle \geq \lambda^{-1} |D^{-1} \xi|^2 \geq \lambda^{-1} |\xi|^2 \geq \left(\lambda K^2\right)^{-1} |\xi|^2,
	\end{align*}
	i.e., the ellipticity constant of $\widetilde A$ is $\lambda K^2$. If one seeks optimal constants, choosing $\widetilde A/K$ the ellipticity constant becomes $\lambda K$. The Lipschitz seminorm of $\widetilde A$ in $D^{-1} (U_{r_0} (\partial \Omega))$ is at most $K^2 [A]_{C^{0,1} (U_{r_0} (\partial \Omega))}$. Indeed, for any two points $x,y \in D^{-1} (U_{r_0} (\partial \Omega))$ with $x\not = y$,
	\begin{align*}
	\frac{|D^{-1}A(Dx)D^{-1}-D^{-1}A(Dy)D^{-1}| }{|x-y|} &\leq K^2 \frac{|A(Dx)-A(Dy)| }{|x-y|}  \\
	&\leq  K^2 [A]_{C^{0,1}(U_{r_0} (\partial \Omega))} \frac{|Dx-Dy| }{|x-y|} \\
	&\leq K^2 [A]_{C^{0,1}(U_{r_0} (\partial \Omega))}.
	\end{align*}
	The $L^\infty$ norm is $\|\widetilde A\|_{L^\infty (\R^2)} \leq K^2 \|A\|_{L^\infty (\R^2)}$.
	
	By \cref{covering elliptic measure} there exists
	$$
	\delta_0 = \delta_0 (K(\lambda)^2 \lambda, K(\lambda)^2 [A]_{C^{0,1}(U_{r_0} (\partial \Omega))} \cdot K(\lambda)^2 \|A\|_{L^\infty (\R^2)}) >0
	$$
	such that if $0<8\sqrt{15} K(\lambda) \delta /15\leq \delta_0$, i.e., $0<\delta\leq 15 \delta_0 / (8\sqrt{15}K(\lambda))$, and taking $M$ big enough such that \cref{covering elliptic measure} holds, then there is a set $\widetilde F \subset \partial D^{-1}(\Omega)$ such that $\widetilde \omega^{\widetilde p} (\widetilde F) \geq C^{-1} \tau$ and with a covering $\widetilde F \subset \bigcup_i B(\widetilde z_i, r_i)$ where
	\begin{enumerate}[label=({\arabic*}*)]
		\item $\sum_i r_i \leq C M^\tau$,
		\item $\sum_{\{i \, : \, r_i > \rho\}} r_i \leq C M^{-1}$,
	\end{enumerate}
	with universal constant $C$.
	
	Defining $F \coloneqq D (\widetilde F)$ we have $\omega^p (F) = \widetilde \omega^{\widetilde p} (\widetilde F) \geq C^{-1} \tau$, and $F\subset \bigcup_i D(B(\widetilde z_i, r_i))$. Finally, as $D(B(\widetilde z_i, r_i)) \subset B(D\widetilde z_i, r_i)$, then $\{B(D\widetilde z_i, r_i)\}_i$ is a covering of $F$ satisfying the same properties.
	\end{proof}
\end{claim}

\subsection{Reduction 2: The dependence on the Lipschitz seminorm}\label{sec:reduction lipschitz seminorm}

This is the second step to show that \cref{covering elliptic measure} implies \cref{only lambda dependence covering elliptic measure}. Note that the flatness constant of the Reifenberg flat domain on \cref{covering elliptic measure} (hence also on \cref{covering elliptic measure sym}) depends also on the Lipschitz seminorm $\kappa$ of the matrix. Below, we see that in fact these results imply \cref{only lambda dependence covering elliptic measure} by a rescaling argument. Here the flatness constant is determined solely by the ellipticity of the matrix.

\begin{proof}[Proof of \cref{only lambda dependence covering elliptic measure} assuming \cref{covering elliptic measure}]
		Fix $r\in (0,1]$ be such that $r\kappa\|A\|_{L^\infty} \leq 1$, and let $\widetilde A (\cdot) = A(r\cdot)$, $\widetilde \Omega = \Omega /r$ and $\widetilde \omega \coloneqq \omega_{\widetilde \Omega, \widetilde A}$ be the elliptic measure with respect to the matrix $\widetilde A$ in $\widetilde \Omega$.
		
		The matrix $\widetilde A$ is Lipschitz in $\{x/r : \dist(x, \partial \Omega) < r_0\} = \{x : \dist(x, \partial \widetilde \Omega) < r_0/r\}$, and $\widetilde \Omega$ is $(\delta, r_0 / r)$-Reifenberg flat since $\Omega$ is $(\delta, r_0)$-Reifenberg flat. By the uniqueness of the elliptic measure we have $\widetilde \omega^{z/r} (\cdot /r) = \omega^z (\cdot)$ for any $z\in \Omega$.
		
		With this ``zoom'' the ellipticity constant of $\widetilde A$ becomes the same as the one of $A$, $\|\widetilde A\|_{L^\infty} = \| A\|_{L^\infty} $ and $[\widetilde A]_{C^{0,1}} = r [ A]_{C^{0,1}} $, which implies $[\widetilde A]_{C^{0,1} } \cdot \|\widetilde A\|_{L^\infty}  = r \kappa \| A\|_{L^\infty} \leq 1$.
		
		This allows us to invoke \cref{covering elliptic measure sym} for the elliptic measure $\widetilde \omega^{p/r}$ since $\dist (p/r, \Omega / r) > r_0 / r$, as $\dist (p, \Omega ) > r_0$. Note that now we don't have the dependence on the Lipschitz seminorm and $L^\infty$ norm of the matrix since we are in the case $[\widetilde A]_{C^{0,1}} \|\widetilde A\|_{L^\infty} \leq 1$. Hence, there exists $\delta_0 = \delta_0 (\lambda) >0$ such that if $0<\delta \leq \delta_0$, then for $M$ big enough (to satisfy \cref{covering elliptic measure sym}) and setting $\rho$ such that $0<\rho/r<1/M$, we can find a set $F' \subset \partial \Omega / r$ such that $\widetilde \omega^{p/r} (F') \geq C^{-1} \tau$ and with a covering $F' \subset \bigcup_i B(z'_i, r_i)$ such that
		\begin{enumerate}[label=({\arabic*}*)]
			\item $\sum_i r_i \leq CM^\tau$,
			\item $\sum_{\{i \, : \, r_i > \rho /r\}} r_i \leq CM^{-1}$.
		\end{enumerate}
		Set $F=rF'$. Then $\omega^p (F) = \widetilde \omega^{p/r} (F') \geq C^{-1} \tau$, and $F=rF' \subset \bigcup_i B(z_i, r r_i)$, which implies
		\begin{enumerate}
			\item $\sum_i r r_i \leq rCM^\tau \leq CM^\tau$,
			\item $\sum_{\{i \, : \, r r_i > \rho\}} r r_i = \sum_{\{i \, : \, r_i > \rho /r\}} r r_i \leq r CM^{-1} \leq CM^{-1}$,
		\end{enumerate}
		as claimed.
\end{proof}

\section{Elliptic measures in CDC domains}

In this section we collect the key properties of elliptic measures in CDC domains for the proof of \cref{only lambda dependence covering elliptic measure}. The first one, frequently called Bourgain's lemma (see \cite[Lemma 1]{Bourgain1987} for the harmonic case), is the following lemma.

\begin{lemma}[{\cite[Lemma 11.21]{Heinonen2006}}]\label{Bourgain lemma}
	Let $\Omega \subset \R^{n+1}$, $n\geq 1$, be a bounded CDC domain with constant $c_0$ and radius $s_0$, and let $A$ be a real uniformly elliptic (not necessarily symmetric) matrix. Then there exists a constant $\tau \in (0,1)$, depending only on $n$, $c_0$ and the ellipticity constant of the matrix $A$, such that for $E\subset \partial \Omega$, $x_0 \in \partial \Omega$ and $0<r\leq s_0$, we have the following:
	\begin{enumerate}
		\item \label{Bourgain lemma upper bound}if $B(x_0, 2r) \cap E = \emptyset$, then $\omega^p_{\Omega, A} (E) \leq 1-\tau < 1$, and
		\item \label{Bourgain lemma lower bound}if $B(x_0, 2r) \cap \partial \Omega \subset E$, then $\omega^p_{\Omega, A}  (E) \geq \tau >0 $,
	\end{enumerate}
	for any point $p\in B(x_0, r) \cap \Omega$.
\end{lemma}

The proof of \cref{only lambda dependence covering elliptic measure} is based on a modification of the domain, without losing the initial information. In the following lemma we obtain the first step in that modification. This is the analogue of \cite[Lemma 1.1]{Wolff1993}.

\begin{lemma}\label{inverse max pple}
	Let $\Omega \subset \R^{n+1}$, $n\geq 1$, be a bounded CDC domain with constants $c_0$ and radius $s_0$, let $x\in \partial \Omega$ and $0<r\leq s_0$, and let $A$ be a real uniformly elliptic (not necessarily symmetric) matrix. Then for any $k > 2$, there exists a constant $C$ depending only on $n$, $c_0$ and the ellipticity constant of the matrix $A$ such that
	$$
	\omega_{\Omega\setminus \overline{B(x,r)},A}^p \left(\overline{B(x,r)} \right) \leq C \omega_{\Omega,A}^p ( B(x,k r)), \text{ for all } p \in \Omega \setminus \overline{B(x,r)}.
	$$
	
	\begin{proof}
		The proof follows the same argument as \cite[Lemma 1.1]{Wolff1993}.
		
		Since $k>2$, by \cref{Bourgain lemma}\rf{Bourgain lemma lower bound} there exists $C > 1$, depending on $n$, $c_0$ and the ellipticity constant of the matrix, such that $\omega^p_{\Omega,A} (B(x,k r)) \geq \omega^p_{\Omega,A} (\overline{B(x,2 r)})\geq C^{-1} $ for any $p \in \overline{B(x,r)} \cap \Omega$. The lemma follows by the maximum principle in $\Omega \setminus \overline{B(x,r)}$ by standard techniques.
	\end{proof}
\end{lemma}

Later on we will need to have some control on the Radon-Nikodym derivative of the elliptic measure of the modified domain with respect to its surface measure, see \cref{density new elliptic measure} below. First we compute the CDC constants of an annulus, which will be used later to control this density in a modified domain.

\begin{lemma}\label{annulus >3}
	Let $k >3$. The annulus $\mathcal A_k = B(0, k^2 r) \setminus \overline{B(0,r)} \subset \R^{n+1}$, $n\geq 1$, satisfies the CDC with constant $c_0=c_0(k)$ and radius $s_0\coloneqq\left(\frac{k+1}{2}\right)^2 r$, and moreover
	\begin{enumerate}
	\item\label{inclusion cdc annulus} $\overline{B(0, k r)} \subset B(x_0, s_0)$ and
	\item\label{does not intersect other component} $\overline{B(x_0, 2s_0)} \subset B(0, k^2 r)$,
	\end{enumerate}
	for any $x_0 \in \partial B(0,r)$, i.e., the inner circle.
\begin{proof}
	In the following computations we find the radius $s_0$ to have the CDC on the annulus $\mathcal A_k$ with the conditions \rf{inclusion cdc annulus} and \rf{does not intersect other component}. From the first condition we get $s_0 > (k+1)r$, and from the second $2s_0 + r < k^2 r$, i.e., $2s_0 < (k^2 -1)r$. In order to have existence in $s_0$ we need $2(k +1)r < (k^2 -1)r$, whence we need $k > 3$. Let $s_0$ be the middle point in $\left( (k+1)r, \frac{k^2-1}{2} r \right)$,
	$$
	s_0 \coloneqq \frac{(k +1)r + \frac{k^2-1}{2} r }{2} = \frac{k^2+2k+1}{4} r \eqqcolon C(k) \cdot r .
	$$
	
	Now we want to see that for $k > 3$, the annulus $\mathcal A_k$ satisfies the capacity density condition with $s_0 =\frac{k^2+2k+1}{4} r = C(k ) \cdot r$. By definition of $C(k)$, given a point $x_0\in \partial \mathcal A_k$, the ball $B(x_0, C(k)\cdot r)$ does not intersect the other component of $\partial \mathcal A_k$, by condition \rf{does not intersect other component}.
	
		We want to see that there exists $c_0 = c_0 (k)$ such that 
		$$
		\Capacity \left( \overline{B (x_0,s)} \cap \mathcal A_k^c, B(x_0, 2s) \right) \geq c_0 (k) \cdot s^{n-1},
		$$
		for all $x_0 \in \partial \mathcal A_k$ and $0<s\leq s_0$.
		
		\textbf{Case 1.} Suppose $x_0 \in \partial B(0,r) \subset \partial \mathcal A_k$. Let $0<s\leq s_0 = C(k) \cdot r$. By the choice of $C(k)$ we have $B(x_0,s) \cap \mathcal A_k^c = B (x_0,s) \cap \overline{B(0,r)}$.
		
		Set $\xi = x_0 - x_0 \frac{s}{2 s_0}$. So $|\xi - x_0| = \left| x_0 \frac{s}{2s_0} \right| = \frac{rs}{2s_0} = \frac{s}{2 C(k)}$, and note that $|\xi - x_0| \leq \frac{r}{2}$. In particular $B\left( \xi, \frac{s}{2 C(k)} \right) \subset B(x_0, s) \cap \mathcal A_k^c$. Also, $B(x_0, 2s) \subset B\left( \xi, 2s + |x_0 - \xi| \right) = B\left( \xi, 2s+\frac{s}{2C(k)} \right)$. From these two inclusions, the monotonicity of the capacity and \cite[(2.13)]{Heinonen2006}, we have
		$$
		\Capacity \left(\overline{B (x_0,s)} \cap \mathcal A_k^c, B(x_0, 2s) \right)
		\geq \Capacity \left(\overline{B \left(\xi, \frac{s}{2C(k)}\right)}, B\left(\xi, 2s + \frac{s}{2C(k)}\right) \right) \approx_k s^{n-1}.
		$$
		
		\textbf{Case 2.} Suppose $x_0 \in \partial B(0,k^2 r) \subset \partial \mathcal A_k$. Let $0<s\leq s_0 = C(k) \cdot r$. Define $\xi = x_0 + \frac{x_0}{|x_0|}\frac{s}{2}$. Hence $B\left( \xi, \frac{s}{2}\right) \subset B(x_0, s) \cap \mathcal A_k^c$ and $B(x_0, 2s) \subset B\left( \xi, \frac{5}{2}s \right)$. Arguing as before we get
		$$
		\Capacity \left(\overline{B (x_0,s)} \cap \mathcal A_k^c, B(x_0, 2s) \right) \gtrsim s^{n-1},
		$$
        as claimed.
\end{proof}
\end{lemma}

Now we study the density of the elliptic measure in an annulus. For a Hölder matrix $A$, here we use that $L_A$-harmonic functions are Hölder continuous up to the boundary, see \cref{C^a to the boundary}.

\begin{lemma}\label{L infty elliptic measure annulus}
	Let $k > 3$, $0<r\leq 1$ and $\mathcal A_k = B(x,k^2 r) \setminus \overline{B(x,r)} \subset \R^{n+1}$, $n\geq 1$. Let $A$ be a real uniformly elliptic (not necessarily symmetric) matrix. Suppose also that $A\in C^\alpha \left( B(x, 2k^2 r) \right)$ with $0<\alpha < 1$. Then the elliptic measure in the annulus $\mathcal A_k$ (arising from the matrix $A$) satisfies
	$$
	\omega_{\mathcal A_k,A}^z (Y) \lesssim \frac{\sigma (Y)}{r^n}, \text{ for any } z\in \partial B(x,k r) \text{ and any } Y\subset \partial B(x,r) ,
	$$
	with constant depending only on $k$, $[A]_{C^\alpha}$ and the ellipticity of $A$.
	\begin{proof}
		Suppose without loss of generality that the annulus is centered at the origin and denote $B_t \coloneqq B (0,t)$. We can also assume that $r=1$. Indeed, denote $\omega \coloneqq \omega_{\mathcal A_k,A}$ and $\widetilde \omega \coloneqq \omega_{\widetilde{\mathcal A_k},\widetilde A}$ the elliptic measure associated to the matrix $\widetilde A (\cdot) \coloneqq A(r\cdot)$, where the rescaled annulus is $\widetilde{\mathcal A_k} = B_{k^2} \setminus \overline{B_1}$. After rescaling and by the uniqueness of the elliptic measure,
		$$
		\omega^z (Y) = \widetilde \omega^{z'} ( Y') \text{ where } z'=z/r \text{ and }  Y'=Y/r ,
		$$
		see \cref{elliptic measure deformation}. The matrix $\widetilde A$ has the same ellipticity constant as $A$, and the Hölder seminorm is improved because $[\widetilde A]_{C^\alpha} = [A]_{C^\alpha} \cdot r^\alpha$ whenever $0<r<1$. If the lemma were true with $r=1$ then writing $p=z/r$ we would get
		$$
		\omega^{z} (Y)
		=\widetilde \omega ^p (Y') 
		\leq C_k \sigma(Y')
		= C_k \frac{\sigma (Y)}{r^n},
		$$
		as claimed.
		
		Let $p \in \partial B_k$ and let $g^T_p$ the Green function of the annulus $\mathcal A_k$ with pole at $p$. Then, using \rf{elliptic measure density},
		\begin{equation}\label{elliptic measure and gradient green function}
		 \omega ^p (Y) = -\int_{Y} \langle { A}^T (\xi) \nabla g^T_p (\xi), \nu (\xi) \rangle \, d\sigma (\xi) 
		\lesssim \int_{Y} |  \nabla g^T_p (\xi) | \, d\sigma (\xi) .
		\end{equation}
		We would be done if we can bound $|\nabla g^T_p|\leq C_k$. To obtain this we apply \cref{C^a to the boundary}. In the next paragraphs we check its hypothesis.
		
		The function $g^T_p$ is $L_{ A^T}$-harmonic in $B_2 \setminus B_1$ since we are in the case $k > 3$. Moreover $g^T_p \equiv 0$ in $\partial B_1$. We need to verify that $g^T_p \in W^{1,2} \left( B_2 \setminus B_1 \right)$. Recall that the Green function is constructed in \rf{def green function} as $g^T_p = -\EE^T_p + h$ where $h$ is a $L_{\widetilde A^T}$-harmonic function with $h=\EE^T_p$ on $\partial {\mathcal A_k}$, and $\EE^T_p$ is the fundamental solution with pole at $p$. Hence,
		\begin{equation}\label{L2 norm green function}
		\begin{aligned}
			\|g^T_p\|_{L^2 \left( B_{2} \setminus B_1 \right)}  \leq &  \|g^T_p\|_{L^2 \left( B_{5/2} \setminus B_1 \right)}
			\leq \max_{y\in \partial B_{5/2} } g^T_p (y) 
			 \leq \max_{y\in \partial B_{5/2} } |\EE^T_p (y) | 
			+ \max_{y\in \partial B_{5/2}  } |h (y)| \\
			\leq & \max_{y\in \partial B_{5/2}  } |\EE^T_p (y) | 
			+ \max_{y\in \partial B_{k^2} \cup \partial B_1 } |\EE^T_p (y)| .
		\end{aligned}
 		\end{equation}
 		In the planar case, we have $|\EE^T_p (y)| \lesssim 1+ \left|\log |y-p|\right|$ by \cref{pointwise_bound_log}, and in particular, from \rf{L2 norm green function} we obtain
 		\begin{equation*}
 				\|g^T_p\|_{L^2 \left( B_{2} \setminus B_1 \right)}  \lesssim \max_{y\in \partial B_{5/2} } \left[1+\left|\log |y-p|\right| \right]
 				+ \max_{y\in \partial B_{k^2}  \cup \partial B_1 } \left[ 1+ \left|\log |y-p|\right| \right] \leq C_k.
 		\end{equation*}
 		In higher dimensions, $n\geq 2$, the fundamental solution is bounded by $|\EE^T_p (y)|\lesssim |p-y|^{1-n}$ (see \cite[Theorem 3.1 (3.55)]{Hofmann2007}). From this bound and \rf{L2 norm green function} we get
 		\begin{equation*}
 			\|g^T_p\|_{L^2 \left( B_{2} \setminus B_1 \right)} \lesssim  \max_{y\in \partial B_{5/2}  } | y-p |^{n-1} 
 			+ \max_{y\in \partial B_{k^2}  \cup \partial B_1 }  | y-p |^{n-1}  \leq C_k .
 		\end{equation*}
 		In fact, we have obtained $\|g^T_p\|_{L^2 \left( B_{5/2} \setminus B_1 \right)} \lesssim C_k$, and hence by Caccioppoli's inequality in the annulus $B_{2} \setminus B_1$ we also obtain the upper bound for the gradient,
 		$$
 		\|\nabla g^T_p\|_{L^2 \left( B_2 \setminus B_1 \right)} 
 		\lesssim 
 		\| g^T_p\|_{L^2 \left( B_{5/2} \setminus B_1 \right)} 
 		 \lesssim C_k ,
 		$$
 		implying $g^T_p \in W^{1,2} (B(0,2)\setminus B(0,1))$ depending only on $k$ and the ellipticity constant.
 		
 		Consider $\Omega = B_2 \setminus \overline{B}_1$, $T=\partial B_1$ and $\Omega' =  B_{3/2} \setminus \overline{B}_1 $ in \cref{C^a to the boundary}. Note that $\Omega^\prime \subset \Omega$, $T\subset \partial \Omega'$ and $\dist \left( \left( B_{3/2} \setminus \overline{B}_1 \right) \cup \partial B_1, \partial B_2 \right) = 1/2$. Then we get $ g^T_p \in C^{1,\alpha} \left( \left( B_{3/2} \setminus \overline{B}_1 \right) \cup \partial B_1  \right)$ with
 		$$
 		\max_{y\in \partial B_1 } |\nabla g^T_p (y)| \leq \|g^T_p\|_{1; \Omega'} \overset{\text{Thm~\ref{C^a to the boundary}}}{\lesssim} \max_{y\in B_2 \setminus \overline{B}_1} g^T_p (y) \leq \max_{y\in \partial B_{5/2} } g^T_p (y) \leq C_k ,
 		$$
 		and the lemma follows.
	\end{proof}
\end{lemma}

\Cref{inverse max pple} relates the elliptic measure on the initial domain with the elliptic measure on the domain minus a fixed ball. Next in \cref{density new elliptic measure}, which is the analogue of {\cite[Lemma 1.2]{Wolff1993}}, we study the elliptic measure on this last setting. Combining \cref{inverse max pple,density new elliptic measure} we will obtain density properties of the elliptic measure on a modified domain.

\begin{lemma}\label{density new elliptic measure}
	Set $k >3$ and $0<r\leq 1$. Let $\Omega \subset \R^{n+1}$, $n\geq 1$, be a bounded Wiener regular domain, $x\in \partial \Omega$ and let $A$ be a real uniformly elliptic (not necessarily symmetric) matrix. Suppose also that $A\in C^\alpha (\{y\in \R^{n+1} : \dist(y, \partial \Omega) < 2k^2 r\})$ with $0<\alpha < 1$.
	
	Set $\widetilde \Omega = \Omega \setminus \overline B$ where $B=B(x,r)$. Then $\omega_{\widetilde \Omega,A}^p \big|_{\partial B}$ is absolutely continuous with respect to $\sigma$ for any $p \in \Omega \setminus k\overline{B}$, and for $z\in \partial B$,
	$$
	\frac{d\omega_{\widetilde \Omega,A}^p}{d\sigma} (z) \leq \frac{C}{r^n} \omega_{\Omega \setminus k\overline{B},A}^p \left( k\overline{B} \right), \text{ for any } p \in \Omega \setminus k\overline{B},
	$$
	with constant $C$ depending only on $n$, $k$, $[A]_{C^\alpha}$ and the ellipticity of $A$.
\end{lemma}

Following the scheme of the proof of \cite[Lemma 1.2]{Wolff1993}, to obtain \cref{density new elliptic measure} we study the elliptic measure of the annulus $\mathcal A_k \coloneqq  B(x,k^2 r) \setminus \overline{B(x,r)}$ when $k > 3$ (in order to apply \cref{annulus >3}) and $0<r\leq 1$ (to have a control on the Hölder seminorm of the matrix $A\in C^\alpha$). However, some technicalities are needed due to the variability of the coefficients of the matrix.

\begin{proof}[Proof of \cref{density new elliptic measure}]
    During the proof we write $\omega_{\cdot}$ instead of $\omega_{\cdot,A}$.
    
    To obtain the result it suffices to prove $\omega^p_{\widetilde \Omega} (Y) \lesssim \frac{\sigma(Y)}{r^n} \omega^p_{\Omega \setminus k \overline{B}} (k \overline{B})$ for all $p \in \Omega \setminus k \overline{ B }$ and every $Y\subset \partial B$, and in fact, it is enough to assume that $Y$ is open. Indeed, fixed $p\in \Omega\setminus k\overline B$, for $\varepsilon >0$ let $U\supset Y$ be an open set (relative to $\partial B$) such that $\sigma(U)\leq \sigma(Y)+\varepsilon$. If the lemma were true for open sets, then $\omega_{\widetilde\Omega}^p (Y)\leq \omega_{\widetilde\Omega}^p (U) \lesssim \frac{\sigma(U)}{r^n} \omega^p_{\Omega \setminus k \overline{B}} (k \overline{B})=\frac{\sigma(Y)+\varepsilon}{r^n} \omega^p_{\Omega \setminus k \overline{B}} (k \overline{B})$ and the general case would follow taking $\varepsilon\to 0$.
	
	Let us assume $Y$ is an open set relative to $\partial B$, and fix $p\in \Omega\setminus k\overline{B}$. Again, for $\varepsilon>0$ let $V\supset k\overline B$ be an open set such that $\omega_{\Omega\setminus k\overline B}^p (V)\leq \omega_{\Omega\setminus k\overline B}^p (k\overline B)+\varepsilon$, $\psi\in C_c (V)$ such that $\characteristic_{k\overline{B}} \leq \psi \leq \characteristic_V$, and let $v_{\Omega \setminus k \overline B}$ denote the $L_A$-harmonic extension of $\psi$ in $\Omega\setminus k\overline B$. 
 
    Let $N\coloneqq \max_{x\in \partial k B} \omega^x_{\widetilde \Omega}(Y)$ and let $x_0 \in \partial k B$ such that $\omega^{x_0}_{\widetilde \Omega}(Y) = N$. Define the annulus $\mathcal A_k \coloneqq  k^2 B \setminus \overline B$. By \cref{difference of two elliptic measures} below we obtain
	\begin{equation}\label{split integral modified domain step 1}
		N-\omega^{x_0}_{\widetilde \Omega \cap \mathcal A_k} (Y) 
		= \omega^{x_0}_{\widetilde \Omega} (Y) - \omega^{x_0}_{\widetilde \Omega \cap \mathcal A_k} (Y)
            = \int_{\Omega \cap \partial k^2 B} \omega^\xi_{\widetilde \Omega} (Y) \, d \omega_{\widetilde \Omega \cap \mathcal A_k}^{x_0} (\xi) .
	\end{equation}
	By the maximum principle in $\widetilde \Omega \setminus k \overline B$ we have that $\omega^\xi_{\widetilde \Omega} (Y) \leq \omega^{x_0}_{\widetilde \Omega}(Y) = N$ for $\xi \in \partial k^2 B$, and hence
	\begin{equation}\label{split integral modified domain step 2}
	\int_{\Omega \cap \partial k^2 B} \omega^\xi_{\widetilde \Omega} (Y) \, d\omega_{\widetilde \Omega \cap \mathcal A_k}^{x_0} (\xi) \leq  N \omega_{\widetilde \Omega\cap \mathcal A_k}^{x_0} (\Omega \cap \partial k^2 B) .
	\end{equation}
    All in all, from \rf{split integral modified domain step 1} and \rf{split integral modified domain step 2} then
	\begin{equation}\label{reduced N and annulus}
	N - \omega^{x_0}_{\widetilde \Omega \cap \mathcal A_k} (Y) \leq N \omega_{\widetilde \Omega\cap \mathcal A_k}^{x_0} (\Omega \cap \partial k^2 B) .
	\end{equation}
	
	Also, by the maximum principle and the fact that the annulus $\mathcal A_k = k^2 B \setminus \overline B$ when $k > 3$ satisfies the CDC with the precise conditions in \cref{annulus >3}, we have that the right-hand side of \rf{reduced N and annulus} is controlled by
	\begin{equation}\label{measure of big circle in annulus}
	\omega_{\widetilde \Omega \cap \mathcal A_k}^{x_0} (\Omega \cap \partial k^2 B) \leq \omega_{\mathcal A_k}^{x_0} (\partial k^2 B) \leq 1-\tau < 1,
	\end{equation}
	for $\tau \in (0,1)$, depending also on $k$. Indeed, this last step follows by applying \cref{Bourgain lemma}\rf{Bourgain lemma upper bound} to any $y \in \partial B$ with the choice of $s_0$ in \cref{annulus >3}, because in that case $x_0 \in \partial k B \subset B(y, s_0)$ and $B(y, 2 s_0) \cap \partial k^2 B =\emptyset$.
	
	From \rf{reduced N and annulus} and \rf{measure of big circle in annulus} we obtain $N-\omega^{x_0}_{\widetilde \Omega \cap \mathcal A_k } (Y) \leq N \omega_{\widetilde \Omega \cap \mathcal A_k}^{p_0} (\partial k^2 B) \leq (1-\tau) N$, equivalently $\omega^{x_0}_{\widetilde \Omega \cap \mathcal A_k} (Y) \geq \tau N = \tau \omega^{x_0}_{\widetilde \Omega} (Y)$. From this and the maximum principle,
	$$
	\tau \omega^{x_0}_{\widetilde \Omega} (Y) \leq \omega^{x_0}_{\widetilde \Omega \cap \mathcal A_k} (Y) \leq \omega^{x_0}_{\mathcal A_k} (Y).
	$$
	
	We have reduced to the elliptic measure in the annulus $\mathcal A_k$. By \cref{L infty elliptic measure annulus} we have $\omega_{\mathcal A_k}^{x_0} (Y) \lesssim \sigma(Y)/r^n$, and so $\omega^{x_0}_{\widetilde \Omega} (Y) \lesssim \omega_{\mathcal A_k}^{x_0} (Y) \lesssim \sigma (Y)/r^n$. Since $x_0 \in \partial kB$ was chosen to achieve the maximum of $\omega^\cdot_{\widetilde \Omega} (Y)$ in $\partial kB$, we obtain that for any $x \in \widetilde\Omega\cap\partial k B$,
	$$
	\omega^x_{\widetilde \Omega} (Y) \leq \omega^{x_0}_{\widetilde \Omega} (Y) \lesssim  \frac{\sigma (Y)}{r^n} = 
	\frac{\sigma (Y)}{r^n} \cdot v_{\Omega \setminus k \overline B} (x),
	$$
	where the last equality is just because $v_{\Omega \setminus k \overline B} (\xi) = \psi=1$ for $\xi\in \widetilde\Omega\cap\partial k B$. Moreover, $\omega^x_{\widetilde\Omega} (Y)=0$ if $x\in \partial\widetilde\Omega \setminus kB$. The same inequality follows for $x \in \widetilde \Omega \setminus k \overline{B}$ by the maximum principle in $\widetilde \Omega \setminus k \overline{B}$. Evaluating at the fixed pole $p\in \Omega\setminus k\overline B$, by the choice of $V$ we have
     $$
     \omega_{\widetilde\Omega}^p (Y)
     \lesssim \frac{\sigma(Y)}{r^n}\cdot v_{\Omega \setminus k \overline B}(p)
     \leq \frac{\sigma(Y)}{r^n}\cdot \omega_{\Omega \setminus k \overline B}^p (V)
     \leq \frac{\sigma(Y)}{r^n}\cdot (\omega_{\Omega \setminus k \overline B}^p (k\overline B)+\varepsilon),
     $$
     and the lemma follows by taking $\varepsilon\to 0$.
\end{proof}

For the sake of completeness here we provide a proof of \rf{split integral modified domain step 1}.

\begin{lemma}\label{difference of two elliptic measures}
    Let $\Omega\subset \R^{n+1}$, $n\geq 1$, be a Wiener regular domain and $A$ be a real uniformly elliptic (not necessarily symmetric) matrix. Let $\widetilde\Omega \subset \Omega$ be a Wiener regular domain. For any Borel set $E\subset\partial\Omega\cap\partial\widetilde\Omega$ and $p\in\widetilde\Omega$, there holds
    $$
    \omega_{\Omega,A}^p (E) - \omega_{\widetilde\Omega,A}^p (E) = \int_{\partial\widetilde\Omega\setminus \partial\Omega} \omega_{\Omega,A}^\xi (E)\, d\omega_{\widetilde\Omega,A}^p (\xi). 
    $$
    \begin{proof}
        During the proof we write $\omega=\omega_{\Omega,A}$ and $\widetilde\omega = \omega_{\widetilde\Omega,A}$.
        
        Fixed $p\in \widetilde\Omega$, for $m\geq 1$, by the inner (just for Borel sets) and outer regularity of Radon measures, let $K_m\subset E$ be a compact set and $U_m \supset E$ be an open set such that
        \begin{equation}\label{regularity interior and exterior Radon}
        \begin{aligned}
            \omega^p (U_m)-1/m &\leq \omega^p (E) \leq \omega^p (K_m)+1/m,\text{ and}\\
            \widetilde\omega^p (U_m)-1/m &\leq \widetilde\omega^p (E) \leq \widetilde\omega^p (K_m)+1/m.
        \end{aligned}
        \end{equation}
        Moreover, we take $K_m\subset K_{m+1}$ and $U_m \supset U_{m+1}$ by redefining the sequences suitably. Finally, take $K\coloneqq\bigcup_{m\geq 1} K_m$ and $U\coloneqq\bigcap_{m\geq 1} U_m$.
        
        Let $\varphi = \varphi_m \in C_c (U_m)$ be such that $\characteristic_{K_m} \leq \varphi \leq \characteristic_{U_m}$, and let $u=u_m$ and $\widetilde u = \widetilde u_m$ denote the $L_A$-harmonic extension of $\varphi$ in $\Omega$ and $\widetilde\Omega$ respectively. By the monotonicity of the integral,
        \begin{subequations}
        \begin{align}
            \label{control elliptic extension by elliptic measure}&\omega^\xi (K_m) \leq u(\xi) \leq \omega^\xi (U_m),\text{ for all }\xi\in \Omega\text{, and}\\
            &\widetilde\omega^\xi (K_m) \leq \widetilde u(\xi) \leq \widetilde\omega^\xi (U_m), \text{ for all } \xi \in \widetilde \Omega.
        \end{align}
        \end{subequations}
        Using \rf{regularity interior and exterior Radon}, we get
        \begin{equation}\label{comparison both elliptic measures and elliptic extensions}
        |\omega^p (E)-\widetilde\omega^p (E) - (u(p)-\widetilde u(p))|\leq 2/m.
        \end{equation}
        On the other hand, note that $u(\xi)-\widetilde u (\xi) = (u(\xi)-\varphi(\xi))\cdot\characteristic_{\partial\widetilde\Omega\setminus \partial\Omega} (\xi)$ if $\xi\in \partial\widetilde\Omega$, and in particular, writing $u-\widetilde u$ as the $L_A$-harmonic extension in $\widetilde\Omega$ of its boundary values in $\partial\widetilde\Omega$ and evaluating at the point $p\in \Omega$ we get
        $$
        u(p)-\widetilde u(p) = \int_{\partial\widetilde\Omega\setminus\partial\Omega} (u(\xi)-\varphi(\xi)) \, d\widetilde\omega^p (\xi).
        $$
        From this, \rf{control elliptic extension by elliptic measure} and $0\leq\int_{\partial\widetilde\Omega\setminus\partial\Omega} \varphi (\xi)\, d\widetilde\omega^p(\xi)\leq \widetilde\omega^p (U_m\setminus K_m)\leq 2/m$, we therefore obtain
        $$
        \int_{\partial\widetilde\Omega\setminus\partial\Omega} \omega^\xi (K_m) \, d\widetilde\omega^p (\xi)-\frac{2}{m}
        \leq u(p)-\widetilde u(p)\leq \int_{\partial\widetilde\Omega\setminus\partial\Omega} \omega^\xi (U_m) \, d\widetilde\omega^p (\xi),
        $$
        which together with \rf{comparison both elliptic measures and elliptic extensions} gives
        $$ 
        \int_{\partial\widetilde\Omega\setminus\partial\Omega} \omega^\xi (K_m) \, d\widetilde\omega^p (\xi) - \frac{4}{m}
        \leq \omega^p (E)-\widetilde\omega^p (E)\leq \int_{\partial\widetilde\Omega\setminus\partial\Omega} \omega^\xi (U_m) \, d\widetilde\omega^p (\xi) + \frac{2}{m}.
        $$
        
        Note that the set $K=\bigcup_{m\geq 1} K_m \subset E$ satisfies $\omega^p (E) = \omega^p (K)$, and since $K$ is measurable, then $\omega^p (E\setminus K)=0$. Since elliptic measures are mutually absolutely continuous for any two different poles, this implies that $\omega^\xi (E\setminus K)=0$ for any $\xi\in \Omega$. Again, since $K$ is measurable we obtain that $\omega^\xi (E) = \omega^\xi (K)$ for any $\xi\in\Omega$. By the same argument, using now that $E$ is Borel and so measurable for every $\omega^\xi$ with $\xi\in \Omega$, the set $U=\bigcap_{m\geq 1} U_m$ satisfies $\omega^\xi (E)=\omega^\xi (U)$ for any $\xi\in \Omega$.
        
        Taking $m\to \infty$, by the monotone convergence theorem, the equation above becomes
        $$
        \int_{\partial\widetilde\Omega\setminus\partial\Omega} \lim_{m\to \infty}\omega^\xi (K_m) \, d\widetilde\omega^p (\xi)
        \leq \omega^p (E)-\widetilde\omega^p (E)\leq \int_{\partial\widetilde\Omega\setminus\partial\Omega} \lim_{m\to \infty}\omega^\xi (U_m) \, d\widetilde\omega^p (\xi),
        $$
        and the lemma follows since $\lim_{m\to \infty}\omega^\xi (K_m) = \omega^\xi (K)=\omega^\xi (E)=\omega^\xi (U)=\lim_{m\to \infty} \omega^\xi (U_m)$.       
    \end{proof}
\end{lemma}

\section{Proof of the weak version of the Main Lemma}\label{proof of weak version of main lemma}

According to the previous reductions in \cref{weakening of the main lemma}, to obtain the \cref{only lambda dependence covering elliptic measure} it suffices to prove \cref{covering elliptic measure}, which we intend to do in this section modulo the proof of \rf{key integral log bound} which is deferred to \cref{sec:log integral}.

In this section we work with bounded $(\delta, r_0)$-Reifenberg flat domains $\Omega \subset \R^2$. Recall that for $\delta >0$ small enough we have that $\Omega$ is an NTA domain (see \cite[Section 3]{Kenig1997}), and hence it satisfies the capacity density condition. See \cref{rem:NTA are CDC}.

\begin{proof}[Proof of \cref{covering elliptic measure}]
	Let $M>0$ be big enough and $0<\rho < 1/M$. Denote $\omega  \coloneqq \omega_{\Omega,A}^p$.
	
	For $x\in \partial \Omega$ define the `high density value' as
	\begin{equation}\label{high density function}
	h(x)\coloneqq \sup \{r\geq \rho : \omega (B(x,r))\geq Mr \},
	\end{equation}
	and $h(x)=\rho$ if the supremum runs over an empty set. Note that $\rho \leq h(x) \leq 1/M$ for every $x\in \partial \Omega$, because $\omega$ is a probability measure.
	
	\begin{definition}[Good balls]
		For $x\in \partial \Omega$, we say that the ball $B(x,r)$ is good, $B(x,r) \in \good$, if $r> h(x)$, i.e., for any $s \geq r$ we have $\omega (B(x,s)) < Ms$.
	\end{definition}
	
	For $x\in \R^2$ define
	$$
	d(x)\coloneqq \inf_{B\in \good} \left[ |x-c(B)|+r(B)  \right] ,
	$$
	where $c(B)$ is the center of the ball and $r(B)$ its radius. For any $B\in \good$ we have that $r(B) \geq h(c(B)) \geq \rho$, which implies $d(x) \geq \rho$.
	
	Note that for every ball $B(\xi,r) \in \good$ we have $r\geq h(\xi)$. Therefore
	$$
	d(x) = \inf_{\xi\in\partial \Omega} \left[ |x-\xi| + h(\xi) \right].
	$$
	This function is $1$-Lipschitz as it is the infimum of $1$-Lipschitz functions.
	
	\begin{rem}
		For $x\in \partial \Omega$ it follows $\rho\leq d(x)\leq h(x)$, and hence $h(x)=\rho$ implies $d(x)= h(x)$. Moreover, if $d(x)=\rho$ then $x\in \partial \Omega$.
	\end{rem}
	
	Let $0<\varepsilon<1$ be small enough (to be fixed in \rf{lower bound kernel and log gradient} below) depending on the ellipticity constant $\lambda$ and the product $\kappa \|A\|_{L^\infty (\R^2)}$. Let $\mathcal{I} = \mathcal{I}_{\varepsilon^{-2}}$ be the family of maximal dyadic cubes $Q \in \mathcal D (\R^2)$ such that $Q\cap \partial \Omega \not = \emptyset$ and $\ell (Q) \leq \varepsilon^2 d(x)$ for all $x\in Q$.
	
	\begin{lemma}\label{properties family I}
		Let $\mathcal I$ be the family defined above. Then:
		\begin{enumerate}
			\item \label{prop1 family I} Every $Q\in \mathcal I$ satisfy $\frac{\varepsilon^2 d(x)}{3} < \ell (Q) \leq 2 \varepsilon^2 d(x)$ for all $x\in \frac{\varepsilon^{-2}}{4} Q$.
			\item \label{prop2 family I} If $Q_1, Q_2 \in \mathcal{I}$ and $\frac{\varepsilon^{-2}}{4} Q_1 \cap \frac{\varepsilon^{-2}}{4} Q_2 \not=\emptyset$, then $\frac{\ell(Q_1)}{6} < \ell (Q_2) < 6\ell (Q_1)$.
			\item \label{prop3 family I}$\{\frac{\varepsilon^{-2}}{4} Q\}_{Q\in \mathcal I}$ has finite superposition, with superposition number $N=N_\varepsilon$ depending on $\varepsilon$ only.
		\end{enumerate}
	\begin{proof}
		Let $x\in \frac{\varepsilon^{-2}}{4} Q$. We start by proving $\ell (Q) \leq 2 \varepsilon^2 d(x)$ in \rf{prop1 family I}. Take any $y\in Q$. By the election of $y$ and since $d$ is $1$-Lipschitz,
		$$
		d(x)=d(y)+d(x)-d(y) \geq \varepsilon^{-2} \ell(Q) -|x-y| \geq \varepsilon^{-2} \ell(Q) - \diam \frac{\varepsilon^{-2}}{4} Q \geq \frac{\ell (Q)}{2\varepsilon^2}.
		$$
		For the other inequality in \rf{prop1 family I}, let $\widehat{Q}$ be the dyadic father of $Q$, i.e., the unique $\widehat Q \in \DD(\R^2)$ such that $Q \subset \widehat Q$ and $\ell (\widehat Q) = 2\ell (Q)$. Since $Q$ is maximal, there exists $y\in \widehat Q$ such that $2\ell (Q) = \ell (\widehat Q) > \varepsilon^2 d(y)$, and hence
		$$
		d(x)=d(x)-d(y)+d(y) < |x-y| + \frac{2\ell(Q)}{\varepsilon^2} \leq \diam \frac{\varepsilon^{-2}}{4} Q + \frac{2\ell(Q)}{\varepsilon^2} \leq \frac{3 \ell (Q)}{\varepsilon^2},
		$$
		and with this we conclude the proof of \rf{prop1 family I}.
		
		Let $Q_1, Q_2 \in \mathcal I$ such that $\frac{\varepsilon^{-2}}{4} Q_1 \cap \frac{\varepsilon^{-2}}{4} Q_2 \not=\emptyset$. Take $x\in \frac{\varepsilon^{-2}}{4} Q_1 \cap \frac{\varepsilon^{-2}}{4} Q_2$, and then \rf{prop2 family I} follows from \rf{prop1 family I} by
		$$
		\ell(Q_1) \leq 2\varepsilon^2 d(x) < 6\ell(Q_2). 
		$$
		
		Given $Q\in \mathcal I$, there is only a finite number $N=N_\varepsilon$ of cubes $P\in \DD (\R^2)$ such that $\ell(Q)/6 < \ell (P) < 6 \ell (Q)$ and $\frac{\varepsilon^{-2}}{4}Q \cap \frac{\varepsilon^{-2}}{4}P \not=\emptyset$, which gives \rf{prop3 family I}.
	\end{proof}
	\end{lemma}

	\begin{lemma}\label{good ball existence}
		There exists $\gamma \geq 1$ such that for every $Q\in \mathcal{I}$ there exists a ball $G_Q \in \good$ with $r(G_Q) \approx \varepsilon^{-2} \ell (Q)$, satisfying the inclusions $\varepsilon^{-2} Q \subset \gamma G_Q$ and $G_Q \subset \gamma \varepsilon^{-2} Q$.
		\begin{proof}
			Given $Q\in \mathcal I$, fix any $x\in Q \cap \partial \Omega$, and take $\xi\in \partial \Omega$ such that 
			\begin{equation}\label{comparability d and x+h}
			d(x)\leq |x-\xi| + h(\xi) \leq 1.1 d(x) .
			\end{equation}
			Define the ball $G_Q \coloneqq B\left( \xi, 2(|x-\xi|+ h(\xi)) \right)$, and hence $G_Q \in \good$, since $r(G_Q) \geq 2h(\xi)\geq h(\xi)+\rho > h(\xi)$.
			
			We claim that $r (G_Q) \approx \varepsilon^{-2}\ell(Q)$. Indeed, from \rf{comparability d and x+h} and \rf{prop1 family I} in \cref{properties family I},
			$$
			\frac{1}{2} r(G_Q) = |x-\xi| + h(\xi) \overset{\text{\rf{comparability d and x+h}}}{\approx} d(x) \overset{\text{\rf{prop1 family I}}}{\approx} \varepsilon^{-2} \ell (Q) .
			$$
			Also, using the previous comparability, the distance between $x$ and $\xi$ is controlled above by
			$$
			|x-\xi| \leq  |x-\xi| + h(\xi) = \frac{1}{2} r(G_Q) \approx \varepsilon^{-2} \ell (Q),
			$$
			which implies that there exists a universal constant $\gamma \geq 1$ such that $\varepsilon^{-2} Q \subset \gamma G_Q$ and $G_Q \subset \gamma \varepsilon^{-2} Q$.
		\end{proof}
	\end{lemma}

	For each cube $Q\in \mathcal I$, fix a point $z_Q \in Q \cap\partial\Omega$ and define $B_Q \coloneqq B(z_Q, r_Q)$ with $r_Q \coloneqq \varepsilon d(z_Q)$. Next, we modify the domain as in \cite{Wolff1993}, but using the family $\{B_Q\}_{Q\in \mathcal I}$. To do so, define
	$$
	\widetilde\Omega\coloneqq\Omega\setminus\bigcup_{Q\in \mathcal I} B_Q,
	$$
	and denote $\widetilde \omega \coloneqq \omega_{\widetilde \Omega,A}^p$ its elliptic measure with pole $p$. From \rf{prop1 family I} in \cref{properties family I} and $d(\cdot) \geq \rho$ on $\partial \Omega$ we have that $\mathcal I$ is finite. In particular, the family $\{B_Q\}_{Q\in \mathcal I}$ is finite, and $\partial \widetilde{\Omega}$ is smooth except at finitely many points.
	
	Recall $\mathcal P$ is the approximating hyperplane in \cref{def:Reifenberg flat}. Since the function $d(\cdot)$ is $1$-Lipschitz, and so $\varepsilon d(\cdot)$ is $\varepsilon$-Lipschitz, the same proof as in \cite[Lemma 2.2]{Azzam2017a} applies to obtain the following lemma.
	
	\begin{lemma}
	Let $r_0 \in (0,\infty]$ and let $\varepsilon >0$ be small enough. Then:
	\begin{enumerate}
		\item (Analogue of {\cite[Lemma 2.2]{Azzam2017a}}) There exists $\delta_0 = \delta_0 (\varepsilon)>0$ such that if $\Omega \subset\R^2$ is $(\delta, r_0)$-Reifenberg flat with $0<\delta<\delta_0$, then the modified domain $\widetilde \Omega$ as above is $(c\varepsilon^{1/2}$, $r_0 /2)$-Reifenberg flat.
		\item (See {\cite[Lemma 2.3(c)]{Azzam2017a}}) For every $Q\in \mathcal I$, there exists a Lipschitz function $f_Q: \mathcal P (z_Q, 30r(B_Q))\cap 10B_Q \to \mathcal P (z_Q, 30r(B_Q))^\perp$ with Lipschitz constant at most $c\varepsilon^{1/2}$.
	\end{enumerate}
	\end{lemma}
	
	For any $Q\in \mathcal I$, by the maximum principle and \cref{inverse max pple} (with $k=10$) respectively,
	\begin{equation}\label{new old elliptic measure}
	\widetilde \omega \left(\overline{B_Q}\right) \leq \omega_{\Omega \setminus \overline{B_Q}} \left( \overline{B_Q} \right) \lesssim \omega (10 B_Q) .
	\end{equation}
	Moreover, for all $z\in \partial \widetilde \Omega \cap \partial B_Q$, by \cref{density new elliptic measure} (with $k =\sqrt {10}$), the maximum principle and \cref{inverse max pple} (with $k = \sqrt {10}$) respectively,
	\begin{equation}\label{density elliptic measure}
	\frac{d\widetilde \omega}{d\sigma} (z) 
	\lesssim \frac{\omega_{\widetilde \Omega \setminus \overline{\sqrt {10} B_Q}} \left(\overline{ \sqrt {10} B_Q} \right) }{r(B_Q)}
	\leq \frac{\omega_{ \Omega \setminus \overline{\sqrt {10} B_Q}} \left( \overline{\sqrt {10} B_Q} \right) }{r(B_Q)} \lesssim \frac{\omega \left( 10 B_Q \right) }{r(B_Q)}.
	\end{equation}
	
	By the existence of a good ball $G_Q \in \good$ with $r(G_Q) \approx \varepsilon^{-2} \ell (Q)$ and $\varepsilon^{-2} Q \subset \gamma G_Q$ (see \cref{good ball existence}), the ball $10 B_Q$ has bounded density with respect to the initial elliptic measure:
	\begin{multline*}
	\omega \left(10 B_Q \right) 
	= \omega\left(B(z_Q, 10 \varepsilon d (z_Q))\right) 
	\overset{\text{L.\ref{properties family I}\rf{prop1 family I}}}{\leq} \omega (3\cdot 3 \cdot 10 \varepsilon^{-1} Q) 
	\overset{\varepsilon \ll 1}{\leq} 
	\omega (\varepsilon^{-2} Q)
	\leq \omega (\gamma G_Q) \\
	\overset{\gamma G_Q \in \good}{\leq} \gamma r(G_Q ) M 
	\approx \gamma \varepsilon^{-2} \ell(Q)M 
	\overset{\text{L.\ref{properties family I}\rf{prop1 family I}}}{\approx} \gamma d(z_Q) M
	= \gamma \varepsilon^{-1} r(B_Q) M .
	\end{multline*}
	In particular, combined with \rf{density elliptic measure} this implies
	\begin{equation}\label{not excessive density}
	\frac{d\widetilde \omega}{d\sigma} (z) \lesssim M \text{ for all }z\in \partial \widetilde \Omega,
	\end{equation}
	where the involved constant depends on $\gamma$ and $\varepsilon$.
	
	Let us smooth out the domain $\widetilde \Omega$. Recall that $\widetilde \Omega$ is smooth except at finitely many points $\{\xi_j\}_{j\in J}$, with $\# J<\infty$ depending on $M$ and $\rho$. Let $0<s<\min_{Q\in\mathcal I} r(B_Q)/1000$ be a small enough parameter to be fixed later. For each point $\xi_j$, $j\in J$, let $B_1, B_2 \in \{B_Q : Q\in \mathcal I\}$ be the two intersecting balls such that $\xi_j \in \partial B_1 \cap \partial B_2$, and let $c_1$, $c_2$ be their centers and $r_1$, $r_2$ be their radii respectively. Take the unique point $q_j\in \widetilde \Omega \cap \partial B(c_1, r_1+s) \cap \partial B(c_2, r_2+s)$, and denote $B_j \coloneqq B(q_j,s)$. The ball $B_j$ is tangent to $B_1$ and $B_2$, and define $T_j$ to be the bounded open region enclosed between $B_1$, $B_2$ and $B_j$. We define the new smooth domain
	$$
	\doublewidetilde \Omega = \doublewidetilde\Omega_s \coloneqq \widetilde \Omega \setminus \bigcup_{j\in J} \overline{T_j},
	$$
	taking small enough $s$ such that
	\begin{enumerate}
		\item for each $Q\in \mathcal I$, $\sigma (\partial B_Q \cap \partial \doublewidetilde \Omega) \geq 0.9 \sigma (\partial B_Q \cap \partial \widetilde \Omega)$, and
		\item $\doublewidetilde \Omega$ is a Lipschitz domain with the same Lipschitz character as $\widetilde \Omega$.
	\end{enumerate}
	Note that $\sigma(\partial\doublewidetilde\Omega\setminus\partial\widetilde\Omega) \leq \# J \cdot 2\pi s$, and so we can take this value to be as small as needed.

	Denote $\doublewidetilde \omega \coloneqq \omega_{\doublewidetilde \Omega,A}^p$. By the maximum principle, $\doublewidetilde \omega (E) \leq \widetilde \omega (E)$ for any $E\subset \partial \doublewidetilde \Omega \cap \partial \widetilde \Omega$. Consequently,
	\begin{equation}\label{newnew old elliptic measure}
		\doublewidetilde \omega (\overline{B_Q}) \leq \widetilde \omega (\overline{B_Q}) \overset{\text{\rf{new old elliptic measure}}}{\lesssim} \omega (10B_Q) \text{ for any } Q\in \mathcal I,
	\end{equation}
	and
	\begin{equation}\label{newnew not excessive density}
		\frac{d\doublewidetilde \omega}{d\sigma} (z) \leq \frac{d\widetilde \omega}{d\sigma} (z) \overset{\text{\rf{not excessive density}}}{\lesssim} M \text{ for all } z \in \partial \doublewidetilde \Omega \cap \partial \widetilde \Omega .
	\end{equation}
	
	Let $K(\cdot) \coloneqq \frac{d\doublewidetilde{\omega}}{d\sigma} (\cdot)$ be the Radon-Nikodym derivative. By \rf{elliptic measure density} we have
	$$
	K(\xi) = -\langle A^T(\xi) \nabla g^T(\xi), \nu(\xi) \rangle \text{ for } \xi \in \partial \doublewidetilde \Omega,
	$$ where $g^T$ is the Green function in $\doublewidetilde \Omega$ with respect to the matrix $A^T$. By \rf{key integral log bound} (proved in \cref{sec:log integral}), if $\varepsilon>0$ is small enough depending on the ellipticity constant $\lambda$ and the product $\kappa\|A\|_{L^\infty (\R^2)}$, then there is a constant $\textrm{const}_{\Omega, A} > 0$ such that 
	\begin{equation}\label{lower bound kernel and log gradient}
	-\infty < -\textrm{const}_{\Omega, A} 
	\leq \int_{\partial \doublewidetilde \Omega} \log |S(\xi) \nabla g^T (\xi)| \, d\doublewidetilde\omega (\xi) 
	= \int_{\partial \doublewidetilde \Omega} K(\xi) \log |S(\xi) \nabla g^T (\xi)| \, d \sigma (\xi) ,
	\end{equation}
	where $A_0 = \frac{A+A^T}{2}$ and $S = A_0^{1/2}$, i.e., $S^T S=A_0$.
	
	For every $\xi \in \partial \doublewidetilde \Omega$ we can write
	$$
	\langle A^T(\xi) \nabla g^T(\xi), \nu(\xi) \rangle 
	= \langle  \nabla g^T(\xi), A(\xi) \nu(\xi) \rangle
	= \langle  \nabla g^T(\xi), c_1 (\xi) \nu(\xi) \rangle + \langle  \nabla g^T(\xi), c_2 (\xi) t(\xi) \rangle,
	$$
	where $\nu$ (resp. $t$) is the outward normal (resp. tangential) vector of $\partial \doublewidetilde \Omega$, and $c_1 (\xi)$ (resp. $c_2 (\xi)$) is the projection of $A(\xi) \nu (\xi)$ into $\nu (\xi)$ (resp. $t(\xi)$). In particular $c_1 (\xi) = \langle A(\xi) \nu (\xi) , \nu (\xi) \rangle \approx 1$, and since $\partial \doublewidetilde \Omega$ is smooth, we have that $\langle  \nabla g^T (\xi), c_2 (\xi) t(\xi) \rangle = c_2 (\xi) \partial_t g^T (\xi) = 0$ in $\partial \doublewidetilde \Omega$ as $g^T |_{\partial \Omega} = 0$. Hence $-\langle A^T(\xi) \nabla g^T(\xi), \nu(\xi) \rangle  \approx |\nabla g^T (\xi)|$, and since $| \nabla g^T(\xi) |^2 \approx \langle A_0 \nabla g^T(\xi), \nabla g^T(x) \rangle = |S(\xi) \nabla g^T (\xi)|^2$, we obtain
	\begin{equation}\label{relation gradient green and kernel}
	|S(\xi) \nabla g^T (\xi)| \leq -C \langle A^T(\xi) \nabla g^T(\xi), \nu(\xi) \rangle = C K(\xi) .
	\end{equation}
	
	From \rf{lower bound kernel and log gradient} and \rf{relation gradient green and kernel},
	\begin{align*}
	-\infty < -\textrm{const}_{\Omega, A}
	& \leq \int_{\partial \doublewidetilde \Omega} K(\xi) \log |S(\xi) \nabla g^T (\xi)| \, d \sigma (\xi) 
	\leq \int_{\partial \doublewidetilde \Omega} K(\xi) \log C K(\xi) \, d \sigma (\xi)  \\
	& = \int_{\partial \doublewidetilde \Omega} K(\xi) \log C \, d \sigma (\xi)  + \int_{\partial \doublewidetilde \Omega} K(\xi) \log K(\xi) \, d \sigma (\xi) \\
	& = \log C +  \int_{\partial \doublewidetilde \Omega} K(\xi) \log K(\xi) \, d \sigma (\xi),
	\end{align*}
	whence we obtain
	\begin{equation}\label{integral log kernel bounded}
		-\infty <  -C_{\Omega, A} \leq \int_{\partial \doublewidetilde \Omega} K(\xi) \log K(\xi) \, d \sigma (\xi) .
	\end{equation}
	
	In view of \rf{integral log kernel bounded} and the fact that $K(\cdot) \lesssim M$ on $\partial \doublewidetilde \Omega\cap\partial \widetilde \Omega$ by \rf{newnew not excessive density}, if $M$ is big enough (provided $s$ is small enough depending on $M$ and $\rho$), then the set of points $\xi\in\partial \doublewidetilde \Omega \cap\partial \widetilde \Omega$ with density $K(\xi)< M^{-\tau}$ indeed has elliptic measure at most $1-\tau/4$, uniformly in $\rho$ and $M$. More precisely, we will obtain
	\begin{equation}\label{points with very small density has small elliptic measure}
	\doublewidetilde\omega\left( \left\{\xi \in \partial \doublewidetilde \Omega \cap\partial \widetilde \Omega : K(\xi) \geq M^{-\tau}\right\} \right)
	=\int_{\left\{\xi \in \partial \doublewidetilde \Omega \cap \partial \widetilde \Omega : K(\xi) \geq M^{-\tau}\right\}} K(\xi) \, d\sigma(\xi) \geq \frac{\tau}{4}
	\end{equation}
	whenever $M$ is big enough depending on $\tau$ and the constant $\textrm{const}_{\Omega, A}$, taking small enough $s$ depending on $M$ and $\rho$. Note that if we had $K(\cdot)\lesssim M$ in the whole boundary $\partial\doublewidetilde\Omega$, using Tchebyshoff's inequality as in \cite[p.~170]{Wolff1993} we would directly get \rf{points with very small density has small elliptic measure} from \rf{integral log kernel bounded}. However, we only have the bound of $K$ on $\partial \doublewidetilde \Omega\cap\partial \widetilde \Omega$ and we will need to estimate the elliptic measure (and the $K\log^+ K\, d\sigma$ measure) of the ``bad'' set $\partial \doublewidetilde \Omega\setminus\partial \widetilde \Omega$.
	\begin{proof}[Proof of \rf{points with very small density has small elliptic measure}]
        Let $\log^+ (\cdot) \coloneqq \max \{0, \log (\cdot)\}$ and $\log^- (\cdot) \coloneqq -\min \{0, \log (\cdot)\}$, and recall that $\doublewidetilde \Omega = \doublewidetilde\Omega_s$. We first establish that
        \begin{equation}\label{eq:KlogK elliptic measure of small set is small}
        \int_{\partial\doublewidetilde\Omega\setminus\partial\widetilde\Omega} K\log^+ K\, d\sigma \to 0 \text{ as } s\to 0,
        \end{equation}
        which in particular implies
        \begin{equation}\label{eq:elliptic measure of small set is small}
        \doublewidetilde\omega(\partial\doublewidetilde\Omega\setminus\partial\widetilde\Omega)=\int_{\partial\doublewidetilde\Omega\setminus\partial\widetilde\Omega} K\, d\sigma
        \leq e\sigma(\partial\doublewidetilde\Omega\setminus\partial\widetilde\Omega)
        +\int_{\partial\doublewidetilde\Omega\setminus\partial\widetilde\Omega} K\log^+ K\, d\sigma \to 0 \text{ as }s\to 0.
        \end{equation}
        As previously noted, for sufficiently small $s>0$, the Lipschitz character of $\doublewidetilde\Omega$ is controlled by the one of $\widetilde\Omega$, and therefore we have that it is NTA with uniform constants. Thus, the involved constants in the subsequent computations will be independent of $s$.

        We claim that there exists $\alpha>0$ such that
        \begin{equation}\label{eq:better decayment doublewidetilde omega}
        \doublewidetilde\omega(B(x,r))\lesssim \left(\frac{r}{r_0}\right)^\alpha \text{ for all } x\in\partial\Omega \text{ and } 0<r\leq r_0/1000,
        \end{equation}
        where $r_0$ is the scale where Reifenberg flatness of $\Omega$ is granted. To prove this, note that there exists $\beta\in (0,1)$ such that $\doublewidetilde\omega(B\setminus B/2)\geq \beta\doublewidetilde\omega(B/2)$ for any ball $B$ centered at $\partial\doublewidetilde\Omega$ with radius $r_B\leq r_0/100$, as a consequence of the doubling property of elliptic measure in NTA domains (see \cite[(1.3.7)]{Kenig1994}) and the fact that $B\setminus B/2\not=\emptyset$ by Reifenberg flatness. It follows that $\doublewidetilde\omega(B/2) \leq \doublewidetilde\omega(B)/(1+\beta)$. Iterating this inequality for a fixed $x\in\partial\Omega$ and $B_0=B(x,r_0/200)$, we have $\doublewidetilde\omega (2^{-k}B_0)\leq \doublewidetilde\omega (B_0)/(1+\beta)^k \leq 1/(1+\beta)^k$ for all $k\geq 0$ integer. This readily implies the existence of $\alpha=\alpha(\beta)>0$ such that \rf{eq:better decayment doublewidetilde omega} holds.
        
        Let us fix $j\in J$ momentarily. For $B_j$, let $q_j=c(B_j)$ denote its center. By the change of pole formula in NTA domains (see \cite[Corollary 1.3.8]{Kenig1994}), we have
        \begin{equation}\label{eq:K is comparable to local Kj}
        K(z)
        \approx \frac{d\doublewidetilde\omega^{q_j}}{d\sigma}(z) \doublewidetilde\omega (\overline{B_j}) \text{ for all }z\in\partial\doublewidetilde\Omega\cap\partial B_j.
        \end{equation}
        Let $\doublewidetilde g_{q_j}^T$ be the Green function of $\doublewidetilde\Omega$ with pole at $q_j$, and $u (\cdot)= \doublewidetilde g_{q_j}^T (s\cdot+q_j)$, which is the Green function of $\{(z-q_j)/s : z\in\doublewidetilde\Omega\}$ with pole at $0$ for the matrix $A(s\cdot+q_j)^T$. Arguing as in \cref{L infty elliptic measure annulus} (using \cref{C^a to the boundary} with $T=\{(z-q_j)/s : z\in \partial\doublewidetilde\Omega\cap \partial B_j\}$) we obtain $s |\nabla \doublewidetilde g_{q_j}^T (z)| = |\nabla u ((z-q_j)/s))| \lesssim \max_{\partial B_{1/2} (0)} u \approx 1$ for all $z\in \partial\doublewidetilde\Omega\cap\partial B_j$, where we used \cite[Lemma 5.4]{Guillen-quasiconformal} in the last step\footnote{In fact, it would suffice to show that $|\nabla \doublewidetilde g_{q_j}^T (z)|\lesssim \frac{|{\log s}|}{s}\leq s^{-\alpha/2-1}$ for small enough $s=s(\alpha,\diam\, \doublewidetilde\Omega)$, which follows from the pointwise estimate in \cref{pointwise_bound_log}.}. Consequently,
        $$
        \frac{d\doublewidetilde\omega^{q_j}}{d\sigma}(z) \lesssim \frac{1}{s} \text{ for all }z\in\partial\doublewidetilde\Omega\cap\partial B_j.
        $$
        
        Using both this and \rf{eq:better decayment doublewidetilde omega} in \rf{eq:K is comparable to local Kj}, since $j\in J$ was arbitrary, we get
        $$
        K(z)
        \lesssim \frac{s^{\alpha-1}}{r_0^\alpha} \text{ for all }z\in\partial\doublewidetilde\Omega\cap\partial\widetilde\Omega.
        $$
        This bound suffices to establish \rf{eq:KlogK elliptic measure of small set is small}, as we can now estimate
        $$
        \int_{\partial\doublewidetilde\Omega\setminus\partial\widetilde\Omega} K\log^+ K\, d\sigma 
        \lesssim \frac{1}{r_0^\alpha}\int_{\partial\doublewidetilde\Omega\setminus\partial\widetilde\Omega} s^{\alpha-1} \log^+ (s^{\alpha-1})\, d\sigma 
        \lesssim \frac{\# J s^\alpha \log^+ (s^{\alpha-1})}{r_0^\alpha},
        $$
        which goes to zero as $s\to 0$.

        Taking $s>0$ small enough in \rf{eq:elliptic measure of small set is small}, to prove \rf{points with very small density has small elliptic measure} it suffices to see
		\begin{equation}\label{points with very small density has small elliptic measure reduction}
			\int_{\left\{\xi \in \partial \doublewidetilde \Omega : K(\xi) \geq M^{-\tau}\right\}} K \, d\sigma \geq \frac{\tau}{2}.
		\end{equation}
		We now prove \rf{points with very small density has small elliptic measure reduction}. First we bound
		\begin{equation}\label{log+ term}
		\int_{\partial \doublewidetilde \Omega} K\log^+ K \, d\sigma
		= \int_{\partial \doublewidetilde \Omega \cap\partial\widetilde\Omega} K\log^+ K \, d\sigma
		+ \int_{\partial \doublewidetilde \Omega\setminus \partial\widetilde\Omega} K\log^+ K \, d\sigma
        \overset{\text{\rf{eq:KlogK elliptic measure of small set is small}}}{\leq} \int_{\partial \doublewidetilde \Omega \cap\partial\widetilde\Omega} K\log^+ K \, d\sigma
		+ 1,
		\end{equation}
        where in the last step we took $s>0$ small enough in \rf{eq:KlogK elliptic measure of small set is small}. Writing $K(\cdot) \leq e^C M$ on $\partial \doublewidetilde \Omega\cap\partial\widetilde\Omega$, see \rf{newnew not excessive density},
		\begin{equation}\label{log+ term in big region}
		\int_{\partial \doublewidetilde \Omega \cap\partial\widetilde\Omega} K\log^+ K \, d\sigma
		\leq (\log M + C )\int_{\left\{\xi \in \partial \doublewidetilde \Omega : K(\xi) \geq 1 \right\}} K \, d\sigma
		\leq \log M\int_{\left\{\xi \in \partial \doublewidetilde \Omega : K(\xi) \geq 1 \right\}} K \, d\sigma + C.
		\end{equation}
		By \rf{log+ term} and \rf{log+ term in big region} we can control the term with $K\log^+K$ as
		\begin{equation*}
		\int_{\partial \doublewidetilde \Omega} K\log^+ K \, d\sigma
		\leq  \log M\int_{\left\{\xi \in \partial \doublewidetilde \Omega : K(\xi) \geq 1 \right\}} K \, d\sigma + C +1.
		\end{equation*}
	
		From this and \rf{integral log kernel bounded},
		\begin{align*}
			\int_{\partial \doublewidetilde \Omega} K \log^- K \, d\sigma 
			&= \int_{\partial \doublewidetilde \Omega}  K \log^+ K \, d\sigma - \int_{\partial \doublewidetilde \Omega}  K \log K \, d\sigma
			\leq \log M\int_{\left\{\xi \in \partial \doublewidetilde \Omega : K(\xi) \geq 1 \right\}} K \, d\sigma + C_{\Omega,A} \\
			&\leq \log M\int_{\left\{\xi \in \partial \doublewidetilde \Omega : K(\xi) \geq M^{-\tau} \right\}} K \, d\sigma + C_{\Omega,A}.
		\end{align*}
		Using this in the last inequality of the following computations, we obtain
		\begin{align*}
			\tau \log M &-  \tau\log M \int_{\left\{\xi \in \partial \doublewidetilde \Omega : K(\xi) \geq M^{-\tau}\right\}} K\, d\sigma
			= \log M^\tau \int_{\partial \doublewidetilde \Omega} K \, d\sigma  -  \log M^\tau \int_{\left\{\xi \in \partial \doublewidetilde \Omega : K(\xi) \geq M^{-\tau}\right\}} K\, d\sigma \\
			&= \log M^\tau \int_{\left\{\xi \in \partial \doublewidetilde \Omega : K(\xi) < M^{-\tau}\right\}} K\, d\sigma
			\leq  \int_{\left\{\xi \in \partial \doublewidetilde \Omega : K(\xi) < M^{-\tau}\right\}} K \log^- K \, d\sigma 
			\leq  \int_{\partial \doublewidetilde \Omega} K \log^- K \, d\sigma \\
			&\leq \log M\int_{\left\{\xi \in \partial \doublewidetilde \Omega : K(\xi) \geq M^{-\tau} \right\}} K \, d\sigma + C_{\Omega,A},
		\end{align*}
		which gives
		$$
		\int_{\left\{\xi \in \partial \doublewidetilde \Omega : K(\xi) \geq M^{-\tau}\right\}} K\, d\sigma \geq \frac{\tau \log M - C_{\Omega,A}}{(1+\tau)\log M}
		= \frac{\tau}{1+\tau} - \frac{C_{\Omega,A}}{(1+\tau)\log M} \geq \frac{\tau}{2},
		$$
		as claimed in \rf{points with very small density has small elliptic measure reduction}, by letting $M>0$ be big enough depending on $\tau$ and $C_{\Omega,A}$.
	\end{proof}
	
	We are now in position to find the final set $F$ with the claimed properties. Let $\{B_{Q_n}\}_n \subset \{B_Q\}_{Q\in \mathcal I} $ be the subfamily satisfying either $(\HD)$ or $(\TS)$:
	\begin{enumerate}
		\item[($\HD$)] $h(z_{Q_n}) > \rho$, where $z_{Q_n}=c(B_{Q_n})$ is the center of the ball $B_{Q_n}$.
		\item[($\TS$)] $h(z_{Q_n}) = \rho$ and $\partial B_{Q_n} \cap \left \{\xi \in \partial\doublewidetilde\Omega\cap\partial \widetilde \Omega : K(\xi) \geq M^{-\tau}\right \} \not = \emptyset$.
	\end{enumerate}
	Here $\HD$ stands for `high density' and $\TS$ for `touching set'.
	
	\begin{notation}
		We write $B_Q \in \HD$ if $B_Q$ satisfies ($\HD$), and $B_Q \in \TS$ if $B_Q$ satisfies ($\TS$).
	\end{notation}

	With this choice,
	\begin{equation*}
		\frac{\tau}{4} \overset{\text{\rf{points with very small density has small elliptic measure}}}{\leq} \doublewidetilde \omega \left( \left\{\xi\in\partial\doublewidetilde\Omega\cap\partial\widetilde\Omega : K(\xi) \geq M^{-\tau}\right\} \right)
		\leq \doublewidetilde \omega \left(\overline{\bigcup_n B_{Q_n}}\right)
		\leq \sum_n \doublewidetilde \omega \left( \overline{B_{Q_n}} \right),
	\end{equation*}
	and by \rf{newnew old elliptic measure},
	$$
	\frac{\tau}{4} \leq \sum_n \doublewidetilde \omega \left( \overline{B_{Q_n}} \right) \lesssim \sum_n \omega (10 B_{Q_n}).
	$$
	If $\varepsilon$ is small enough, then $10 B_Q = B\left( z_Q, 10 \varepsilon d(z_Q) \right) \subset \frac{\varepsilon^{-2}}{4} Q$ for each $Q\in \mathcal I$. Moreover, since $\{\frac{\varepsilon^{-2}}{4} Q\}_{Q\in \mathcal I}$ has finite overlapping by \cref{properties family I}\rf{prop3 family I}, $\{10 B_Q\}_{Q\in \mathcal I}$ has also finite overlapping with constant depending on $\varepsilon$ only. From this we obtain,
	\begin{equation}\label{first candidate of final set lemma}
	\tau \lesssim \sum_n \omega (10 B_{Q_n}) \lesssim \omega \left( \bigcup_n 10 B_{Q_n} \right).
	\end{equation}
	
	At this point we have found a subset of $\partial \Omega$ (covered by balls) with elliptic measure bounded uniformly from below. Moreover, the radii $\varepsilon d(z_{Q_n})$ of these balls $B_{Q_n}$ are smaller than the `high density value' $h(z_{Q_n})$, which will allow us to have a control on the sum of the radii.
	
	First we need to define the set $F \subset \partial \Omega$ and its covering. For each $B_{Q_n}\in\HD$, we have $10 B_{Q_n} = B\left( z_{Q_n}, 10 \varepsilon d(z_{Q_n}) \right)$, and since $d(\cdot) \leq h(\cdot)$ on $\partial \Omega$, $10 B_{Q_n}  \subset B\left( z_{Q_n}, 10 \varepsilon h(z_{Q_n}) \right) \subset B\left( z_{Q_n}, 10 h(z_{Q_n}) \right)$. Since the family $\left\{ B\left( z_{Q_n}, 10 h(z_{Q_n}) \right) \right\}_{B_{Q_n} \in\HD}$ is finite, by means of the $3R$-covering theorem consider a disjoint subfamily
	$$
	\left\{ B_m \right\}_m
	\subset 
	\left\{ B\left( z_{Q_n}, 10  h(z_{Q_n}) \right) \right\}_{B_{Q_n} \in\HD}
	$$
	such that
	$$
	\bigcup_{B_{Q_n} \in\HD} B\left( z_{Q_n}, 10 h(z_{Q_n}) \right) \subset \bigcup_m 3B_m .
	$$
	
	Let us define
	$$
	F \coloneqq \left(\bigcup_m  3B_m \cup \bigcup _{B_{Q_n} \in\TS} 10 B_{Q_n}\right) \cap \partial \Omega .
	$$
	Note that $\bigcup_n 10 B_{Q_n} \subset \bigcup_m  3B_m \cup \bigcup _{B_{Q_n} \in \TS} 10 B_{Q_n}$ implies $\omega (F) \gtrsim \tau$ by \rf{first candidate of final set lemma}. Next we show that the covering
	$$
	\left\{ 3B_m \right\}_m \cup 
	\left\{ 10 B_{Q_n}  \right\}_{B_{Q_n} \in \TS}
	$$
	satisfies the properties \rf{covering elliptic measure cond1} and \rf{covering elliptic measure cond2} in \cref{covering elliptic measure}.
	
	We can control the radius of the balls with high density
	$$
	\sum_m r(3B_m)
	= 30  \sum_m h\left( c(B_m) \right) 
	\overset{\text{\rf{high density function}}}{\leq} \frac{30 }{M} \sum_m \omega \left( B_m \right)
	= \frac{30 }{M} \omega \left( \bigcup_m B_m \right)
	\lesssim \frac{1}{M},
	$$
	by the definition of $h(\cdot)$ in \rf{high density function}, and the fact that the balls $\{B_m\}_m$ are pairwise disjoint. Also, $r\left( 3B_m \right) = 30  h\left( c(B_m) \right) > 30 \rho$. We have shown the second property of the covering.
	
	Recall that the balls in $\TS$ intersect the set $\left \{\xi \in \partial\doublewidetilde\Omega \cap \partial \widetilde \Omega : K(\xi) \geq M^{-\tau}\right\}$. For each $B_{Q_n} \in \TS$ consider a point $x_{Q_n} \in \partial B_{Q_n} \cap \left \{\xi \in \partial\doublewidetilde\Omega \cap \partial \widetilde \Omega : K(\xi) \geq M^{-\tau}\right\}$. Then, by \rf{newnew not excessive density} and \rf{density elliptic measure} we have
	$$
	M^{-\tau} \leq \frac{d\doublewidetilde \omega}{d\sigma} (x_{Q_n})
	\overset{\text{\rf{newnew not excessive density}}}{\leq} \frac{d\widetilde \omega}{d\sigma} (x_{Q_n})
	\overset{\text{\rf{density elliptic measure}}}{\lesssim} \frac{\omega \left( 10 B_{Q_n} \right) }{r(B_{Q_n})},
	$$
	which implies
	\begin{equation*}
	\sum_{B_{Q_n} \in\TS} r\left( 10 B_{Q_n} \right) 
	\lesssim M^{\tau} \sum_{B_{Q_n} \in \TS} \omega \left( 10 B_{Q_n} \right)
	\lesssim M^{\tau}  \omega \left( \bigcup_{{Q_n} \in \TS} 10 B_{Q_n} \right) 
	\leq M^\tau ,
	\end{equation*}
	obtaining the first property of the covering. Note that for $B_{Q_n} \in \TS$ we have $h(z_{Q_n}) = \rho$, and in particular $d(z_{Q_n}) = \rho$, which implies $r(10 B_{Q_n}) = 10 \varepsilon d (z_Q) =  10 \varepsilon \rho$.
\end{proof} %proof of theorem covering

\section{Proof of \texorpdfstring{\cref{dimension elliptic measure}}{Theorem~\ref{dimension elliptic measure}}}\label{sec:proof of thm}
In this section we will prove \cref{dimension elliptic measure}. First we make the reduction to the case of bounded domains (\cref{bounded implies unbounded theorem}) and then we prove the theorem using \cref{only lambda dependence covering elliptic measure}.

\subsection{Reduction to bounded domains}\label{sec:proof of thm-reduction to bounded domains}

First, we state some lemmas.

Let $\Omega \subset \R^{n+1}$ with $n \geq 1$, and let $B$ be a ball centered at $\partial \Omega$. By the maximum principle \cite[p.~46]{Gilbarg2001} we have 
\begin{equation}\label{easy inequality localized elliptic measure}
	\omega^z_{\Omega \cap B,A} (E) \leq \omega^z_{\Omega,A} (E) \text{ for any } E\subset \partial \Omega \cap B \text{ and } z\in \Omega \cap B .
\end{equation}

The converse inequality may fail. However, the following weaker relation holds.

\begin{lemma}\label{hard inequality localized elliptic measure}
	Let $\Omega \subset \R^{2}$ be a (possibly unbounded) Wiener regular domain, $A$ be a real uniformly elliptic matrix and $B$ be a ball centered at $\partial \Omega$ with $\Capacity (B \cap \partial \Omega, 4B) \not =0$. Let $E\subset \partial \Omega \cap B$ be a Borel set and $z_E \in \partial 1.5 B$ such that $\omega^{z_E}_{\Omega,A} (E) = \max_{z\in \partial 1.5 B} \omega^z_{\Omega,A} (E)$. Then
    $$
    \omega^{z_E}_{\Omega,A} (E)
    \lesssim \frac{\Capacity (2B, 4B)}{\Capacity (B \cap \partial \Omega, 4B)}\omega^{z_E}_{\Omega \cap 4B,A} (E),
    $$
    where the constant involved depends on the ellipticity constant of the matrix $A$ and the dimension. The same also holds for bounded Wiener regular domains $\Omega \subset \R^{n+1}$ when $n\geq 2$.
	\begin{proof}
		During this proof we write $\omega_{\cdot}$ instead of $\omega_{\cdot,A}$.
		
		By \cref{difference of two elliptic measures} we have
		\begin{equation}\label{two terms in elliptic measure}
            \omega^{z_E}_\Omega (E) - \omega^{z_E}_{\Omega \cap 4B} (E)
		= 
		\int_{\partial 4B \cap \Omega} \omega^\xi_\Omega (E) \, d\omega^{z_E}_{\Omega \cap 4B} (\xi) .
		\end{equation}
            By the maximum principle\footnote{In the planar case the maximum principle on $\Omega \setminus 1.5\overline B$ holds even if $\Omega$ is not bounded, since $\omega^\cdot_\Omega (E) \in [0,1]$.} in $\Omega \setminus 1.5\overline B$, $\omega^\xi_\Omega (E) \leq \omega^{z_E}_\Omega (E)$ for $\xi \in 4\partial B \cap \Omega$. From this and \rf{two terms in elliptic measure}, we get
		\begin{equation}\label{relation localized elliptic measure}
                \omega^{z_E}_\Omega (E) \leq  \omega^{z_E}_{\Omega \cap 4B} (4 \partial B \cap \Omega) \cdot \omega^{z_E}_\Omega (E) + \omega^{z_E}_{\Omega\cap 4B} (E).
		\end{equation}
		
		It remains to bound $\omega^{z_E}_{\Omega \cap 4B} (4 \partial B \cap \Omega)$. By the maximum principle we have 
		$$
		\omega^{z_E}_{\Omega \cap 4B} (4 \partial B \cap \Omega) \leq \omega^{z_E}_{4B \setminus (B \cap \partial \Omega)} (4 \partial B) .
		$$
		By \cite[Lemma 6.21]{Heinonen2006}, since $z_E \in 1.5 \partial B$, we have
		\begin{equation}\label{constant relation localized elliptic measure}
			1- \omega^{z_E}_{4B \setminus (B \cap \partial \Omega)} (4 \partial B) \geq \widetilde c \cdot  \Capacity (B \cap \partial \Omega, 4B) / \Capacity (2B, 4B) \eqqcolon c \in (0,1).
		\end{equation}
		In particular
		\begin{equation}\label{reduced domain does not see big ball}
		\omega^{z_E}_{\Omega \cap 4B} (4 \partial B \cap \Omega) \leq \omega^{z_E}_{4B \setminus (B \cap \partial \Omega)} (4 \partial B) \leq 1-c .
		\end{equation}
		
		From \rf{relation localized elliptic measure} and \rf{reduced domain does not see big ball} we get
		\begin{equation*}
			\omega^{z_E}_\Omega (E)  \leq  \omega^{z_E}_{\Omega \cap 4B} (4 \partial B \cap \Omega) \cdot \omega^{z_E}_\Omega (E) + \omega^{z_E}_{\Omega\cap 4B} (E) 
			 \leq  (1-c)  \cdot \omega^{z_E}_\Omega (E) + \omega^{z_E}_{\Omega\cap 4B} (E),
		\end{equation*}
		obtaining
		$$
		c \cdot \omega^{z_E}_\Omega (E) \leq \omega^{z_E}_{\Omega\cap 4B} (E) ,
		$$
		as claimed, with $c$ as in step \rf{constant relation localized elliptic measure}.
	\end{proof}
\end{lemma}

We remark that $\Capacity (B \cap \partial \Omega, 4B)=0$ would imply $\omega_D (B\cap \partial \Omega)=0$ for any domain $D$ with $B\cap\partial\Omega\subset\partial D$, see \cite[Theorems 10.1 and 11.14]{Heinonen2006}. For $(\delta, r_0)$-Reifenberg flat domains, $\Capacity (B \cap \partial \Omega, 4B) / \Capacity (2B, 4B)\approx 1$ whenever $r_B\leq r_0$. In fact it is only needed the exterior Corkscrew condition, see \cref{rem:NTA are CDC}.

\begin{lemma}\label{localization argument}
	Let $\Omega \subset \R^2$ be a (possibly unbounded) Wiener regular domain and $A$ be a real uniformly elliptic matrix. Let $\{B_i\}_i$ be a pairwise disjoint collection of balls centered at $\partial \Omega$ with $\omega_{\Omega,A}\left(\partial \Omega \setminus \bigcup_i B_i \right)=0$, and $F_i \subset \partial (\Omega \cap 4 B_i)$ with $\omega_{\Omega \cap 4 B_i ,A} (F_i)=1$. Then $F\coloneqq \partial \Omega \cap \bigcup_i F_i$ satisfies $\omega_{\Omega,A} (F)=1$. The same also holds for bounded Wiener regular domains $\Omega \subset \R^{n+1}$ when $n\geq 2$.
	\begin{proof}
		In this proof we denote $\omega \coloneqq \omega_{\Omega,A}$ and $\omega_{\Omega \cap 4 B_i} \coloneqq \omega_{\Omega \cap 4 B_i ,A}$.
		
		Let $p\in \Omega\setminus \{\infty\}$. Abusing notation, we write $F^c=\partial\Omega\setminus F$. Since $\omega^p (F^c \setminus \bigcup_i B_i)  \leq \omega^p (\partial \Omega \setminus \bigcup_i B_i) =0$ and the balls are pairwise disjoint, we can conclude
		\begin{equation}\label{split the set in each ball}
			\omega^p (F^c) = \omega^p \left( \bigcup_i F^c \cap B_i \right)=\sum_i \omega^p \left( F^c \cap B_i \right).
		\end{equation}
		We claim that each term in the right-hand side is zero. Indeed, for each $i$ fix a pole $p_i\in \Omega\cap 4B_i$, by the Borel regularity of $\omega_{\Omega\cap 4B_i}^{p_i}$ let $E_i \supset F^c \cap B_i$ be a Borel set with $\omega_{\Omega\cap 4B_i}^{p_i} (E_i)=\omega_{\Omega\cap 4B_i}^{p_i} (F^c\cap B_i)$, and finally let $z_i \in 1.5\partial B_i$ such that $\omega^{z_i} (E_i) = \max_{z\in 1.5\partial B_i} \omega^z (E_i)$. With this choice of $z_i$, by \cref{hard inequality localized elliptic measure} we have
	\begin{equation}\label{finding the set for each ball}
		\omega^{z_i} (E_i) \lesssim \omega^{z_i}_{\Omega \cap 4B_i} (E_i). 
	\end{equation}
		Since the balls are pairwise disjoint we have that $F^c \cap B_i \subset (F_i^c \cap B_i) \cap \partial \Omega \subset F_i^c$, so
        $$
        \omega_{\Omega\cap 4B_i}^{p_i} (E_i)=\omega^{p_i}_{\Omega \cap 4B_i} (F^c \cap B_i) = 0.
        $$
        Hence by the Harnack inequality (denoting its use with $H$), \rf{finding the set for each ball} and $\omega_{\Omega\cap 4B_i}^{p_i} (E_i) = 0$, for every index $i$ there holds
		\begin{equation}\label{the set has zero measure in each ball}
			\omega^p (F^c \cap B_i) \leq \omega^p (E_i) 
            \overset{(H)}{\lesssim} \omega^{z_i} (E_i)
            \overset{\text{\rf{finding the set for each ball}}}{\lesssim} \omega_{\Omega\cap4B_i}^{z_i} (E_i)
            \overset{(H)}{\lesssim} \omega_{\Omega\cap 4B_i}^{p_i} (E_i)
            = 0 .
		\end{equation}
		Notice that for each $i$ the constants involved in the use of Harnack inequality and \cref{hard inequality localized elliptic measure} in \rf{finding the set for each ball} and \rf{the set has zero measure in each ball} depend on $i$, but the right-hand side in \rf{finding the set for each ball} is zero.
		
		By \rf{the set has zero measure in each ball} we have that the sum in the right-hand side of \rf{split the set in each ball} is zero as claimed. Therefore the set $F\coloneqq\partial \Omega \cap \bigcup_i F_i$ satisfies $\omega^p (F)=1$.
	\end{proof}
\end{lemma}

\begin{lemma}\label{bounded reifenberg flat}
	Let $r_0\in (0,\infty]$ and let $\varepsilon >0$ be small enough. There exists $\delta_0 = \delta_0 (\varepsilon) >0$ such that if $\Omega \in \R^2$ is $(\delta,r_0)$-Reifenberg flat with $\delta \in (0,\delta_0)$ and $B$ is a ball centered at $\partial \Omega$ with radius $r_B \leq r_0 / 100$, then there exists a bounded $(\varepsilon, r)$-Reifenberg flat domain $D \subset \{z : \dist (z,\partial \Omega) < r_0 / 2\}$ (for some $r\in (0,r_0 / 2)$) with $\Omega \cap 4B = D \cap 4B$.
\end{lemma}

Note that we are not interested in the precise dependence of $r$ with respect to $r_0$, because we seek for a qualitative result in \cref{dimension elliptic measure}. It is quite likely that with some care the previous result could be made quantitative.

\begin{proof}
	The proof uses the construction in \cite[Definition 2.1 and Lemma 2.2]{Azzam2017a}. Set $0<\varepsilon<1/100$ and $E\coloneqq \partial \Omega \cap 5\overline B$. Let $\WW_{\varepsilon^{-2}} (E^c)$ be the set of maximal dyadic cubes $Q\subset E^c$ such that $\diam (\varepsilon^{-2} Q) \leq r_0$ and $\varepsilon^{-2} Q \cap E = \emptyset$.
	
	Denote $\mathcal I$ the family of cubes $Q\in \WW_{\varepsilon^{-2}} (E^c)$ such that $Q\cap \partial \Omega\not = \emptyset$. For each cube $Q\in \mathcal I$ fix a point $z_Q \in Q \cap \partial \Omega$, and set $r_Q \coloneqq \varepsilon \min\{r_0, \dist (z_Q,E)\}$ and $B_Q \coloneqq B(z_Q, r_Q)$. Consider the enlarged domain
	$$
	\Omega_\varepsilon^+ \coloneqq \Omega \cup \bigcup_{Q\in \mathcal I} B_Q \supset \Omega.
	$$
	By \cite[Lemma 2.2]{Azzam2017a}, this new domain is $(c\varepsilon^{1/2},r_0 /2)$-Reifenberg flat, where the constant $c$ depends only on the dimension, provided the initial domain is $(\delta, r_0)$-Reifenberg flat with $\delta \in (0,\delta_0)$ and $\delta_0$ is small enough depending on $\varepsilon$.
	
	Consider the domain $D_0 \coloneqq \Omega_\varepsilon^+ \cap 10B$. Clearly $D_0 \subset \{z : \dist(z,\partial\Omega) < r_0 / 2\}$ as $10r_B \leq r_0/10$. Let us smooth the corners of $D_0$ out, where the $(c\varepsilon^{1/2},s)$-Reifenberg flat condition fails for all $s>0$. Note that this may only happen in a finite number of points $\{\xi_j\}_{j\in J} \in \partial D_0 \cap \partial 10B$ because $\Omega_\varepsilon^+$ is contructed as a countable union of balls. Fix a small parameter $\tau$. For each $\xi_j$ of these points in $\partial D_0\cap \partial 10B$, let $B_j \in \{B_Q : Q\in \mathcal I\}$ such that $\xi_j \in \partial 10B \cap \partial B_j$, and let $c_j$ and $r_j$ denote its centers and radii respectively. Consider now the unique point $p_j\in D_0 \cap (\partial B (c_B, 10r_B - \tau) \cap \partial B(c_j, r_j-\tau))$. In particular, the ball $B (p_j,\tau)$ is tangent to $\partial 10B$ and $\partial B_j$. Let $T_j = 10B \cap B_j \cap B(p_j, \tau)^c$ the bounded open region enclosed between the previous balls. Taking $\tau$ to be small enough, the final domain
	$$
	D \coloneqq D_0 \setminus \bigcup_{j\in J} T_j \subset D_0
	$$
	satisfies $D\cap 4B = \Omega_\varepsilon^+ \cap 4B=\Omega\cap 4B$ and is $(c\varepsilon^{1/2},r)$-Reifenberg flat for some $r>0$ depending on $\tau$.
	\end{proof}

\begin{claim}\label{bounded implies unbounded theorem}
	If \cref{dimension elliptic measure} holds for bounded $(\delta_0, r_0)$-Reifenberg flat domains, then there exists $\delta_1 = \delta_1 (\delta_0) >0$ such that \cref{dimension elliptic measure} holds for unbounded $(\delta_1, r_0)$-Reifenberg flat domains.
	\begin{proof}
		First we want to remark that if \cref{dimension elliptic measure} holds for $(\delta, r_0)$-Reifenberg flat domains for a fixed $r_0>0$, then by means of a dilation it holds for $(\delta, r)$-Reifenberg flat domains for any $r>0$.
		
		Let $\varepsilon >0$ given by \cref{bounded reifenberg flat} and let $\varepsilon^\prime \coloneqq \min \{\varepsilon, \delta_0 /2\}$ be small enough. Let $\delta_1 = \delta_1 (\varepsilon^\prime)$ given by \cref{bounded reifenberg flat}. Let $\Omega \subset \R^2$ be an unbounded $(\delta_1, r_0)$-Reifenberg flat domain. Let $\{B_i\}_i \subset \{B(\xi,r)\}_{\xi\in \partial \Omega, 0<r<r_0 /100}$ be a disjoint family with $\omega_{\Omega,A} (\partial \Omega \setminus \bigcup_i B_i) = 0$, by Vitali's covering theorem. For each ball $B_i$ let $D_i$ be the bounded $(\varepsilon, r_i)$-Reifenberg flat domain from \cref{bounded reifenberg flat}, for some $r_i \in (0, r_0 /2)$. As $\varepsilon^\prime < \delta_0$, in particular each $D_i$ is a bounded $(\delta_0, r_i)$-Reifenberg flat domain.
		
		As we are assuming that \cref{dimension elliptic measure} holds for bounded $(\delta_0, r)$-Reifenberg flat domains for any $r>0$, for each $i$ take $F_i \subset \partial D_i$ with $\omega_{D_i,A} (F_i)=1$ and $\sigma$-finite one-dimensional Hausdorff measure. From $\Omega \cap 4B_i= D_i \cap 4B_i$, the maximum principle and $\omega_{D_i} (F_i)=1$ we get
		$$
		\omega_{\Omega \cap 4B_i,A} ((\partial \Omega \setminus F_i)\cap 4B_i) = \omega_{D_i \cap 4B_i,A} ((\partial D_i \setminus F_i)\cap 4B_i) \leq \omega_{D_i,A} ((\partial D_i \setminus F_i)\cap 4B_i) = 0 .
		$$
		In particular $\omega_{\Omega \cap 4B_i,A} ((F_i \cap 4B_i)\cup \partial 4B_i)=1$. As $F_i$ has $\sigma$-finite length, so does $\widetilde F_i \coloneqq (F_i \cap 4B_i)\cup \partial 4B_i$. By \cref{localization argument}, the set $F=\partial \Omega \cap \bigcup_i \widetilde F_i$ satisfies $\omega_\Omega (F)=1$, and clearly has $\sigma$-finite length.
	\end{proof}
\end{claim}
 
\subsection{Proof for bounded domains}
 
\Cref{dimension elliptic measure} follows from \cref{only lambda dependence covering elliptic measure} as it is done in \cite[Proof of Theorem 1]{Wolff1993}, with some small modifications. For the sake of completeness we give the detailed proof.

	\begin{proof}[Proof of \cref{dimension elliptic measure}]
		By \cref{bounded implies unbounded theorem} we can assume without loss of generality that $\Omega$ is bounded. We denote $\omega\coloneqq \omega_{\Omega,A}$.
		
		Let $\phi: [0,\infty) \to [0,\infty)$ be any increasing function with $\lim_{t\to 0} \phi(t)/t = 0$, and consider the $\phi$-Hausdorff content
		$$
		h_\phi (E) =\inf\left\{ \sum_i \phi(r_i) : E\subset \bigcup_i B(z_i, r_i) \right\} .
		$$
		
		Now we claim that there exists  $F_\phi \subset \partial \Omega$ with $\omega (F_\phi) = 1$ and $h_\phi (F_\phi)=0$. Indeed, suppose that the pole $z\in \Omega$ is such that $\dist(z,\partial \Omega)>r_0$. Set $\tau = 1/2$ and fix $0<r\leq 1$ to be small enough as in \cref{only lambda dependence covering elliptic measure}. Let $\varepsilon >0$ and fix $M>\varepsilon^{-1}$ satisfying the hypothesis in \cref{only lambda dependence covering elliptic measure}. Take $\rho$ small enough such that $0<\rho /r < 1/M$ and $\phi(\gamma) < \varepsilon M^{-1/2} \gamma$ for all $\gamma \leq \rho$. Then, by \cref{only lambda dependence covering elliptic measure} with these parameters, we obtain a set $F_\varepsilon \subset \partial \Omega$ such that $\omega^z (F_\varepsilon)\geq C^{-1}$ and with a covering $F_\varepsilon \subset \bigcup_i B(z_i, r_i)$ with $\sum_i r_i \leq C M^{1/2}$ and $\sum_{\{i \, : \, r_i >\rho\}} r_i \leq CM^{-1}$. This covering satisfies
		$$
		h_\phi (F_\varepsilon) \leq \sum_{r_i\leq \rho} \phi (r_i) + \sum_{r_i > \rho} \phi (r_i) 
		\leq \varepsilon M^{-1/2} \sum_{r_i \leq \rho} r_i + \sum_{r_i > \rho}  r_i 
		\leq \varepsilon M^{-1/2} CM^{1/2} + CM^{-1} \leq C\varepsilon .
		$$
		
		Define now $F_\infty \subset \partial \Omega$ as
		$$
		F_\infty \coloneqq \limsup_{j\to \infty} F_{1/j^2} = \bigcap_{j\geq 1} \bigcup_{k\geq j} F_{1/k^2}.
		$$
		With this choice we have
		$$
		\omega^z (F_\infty ) 
		= \lim_{j\to \infty} \omega^z \left(\bigcup_{k\geq j} F_{1/k^2}\right)
		 \geq \limsup_{j\to \infty} \omega^z (F_{1/j^2}) \geq C^{-1} ,
		$$
		and as $F_\infty \subseteq \bigcup_{k\geq j} F_{1/k^2}$ for any $j\geq 1$, then
		$$
		h_\phi (F_\infty) \leq h_\phi \left( \bigcup_{k\geq j} F_{1/k^2} \right)
		\leq \sum_{k\geq j} h_\phi (F_{1/k^2}) \leq C \sum_{k\geq j} \frac{1}{k^2} ,
		$$
		which gives $h_\phi (F_\infty)=0$ letting $j\to \infty$.
		
		Let $\{z_k\}_{k=1}^\infty$ be a countable dense subset of $\Omega$. Fix $z_k$ and set $d_k = \dist (z_k, \partial \Omega)$. As $\Omega$ is $(\delta, r_0)$-Reifenberg flat and the matrix $A$ is Lipschitz in $\{x : \dist (x,\partial \Omega) < r_0\}$, then $\Omega$ is $(\delta, r_k)$-Reifenberg flat and $A$ is Lipschitz in $\{x : \dist (x,\partial \Omega) < r_k\}$ with $r_k = \min \{ r_0, d_k /2 \}$. By the choice of $r_k$ we are in the situation $\dist(z_k , \partial \Omega) > r_k$. By the same argument done in the previous paragraphs we get a set $F_k \subset \partial \Omega$ (relative to $z_k$) such that
		$$
		\omega^{z_k} (F_k) \geq C^{-1} \text{ and } h_\phi (F_k) = 0.
		$$
		
		Define $F_\phi \coloneqq \bigcup_{k=1}^\infty F_k$. The condition $h_\phi (F_k) = 0$ for every $k\geq 1$ means that for every $\epsilon >0$ there exists a covering $F_k \subset \bigcup_i B(z_i^{k,\epsilon}, r_i^{k,\epsilon})$ such that $\sum_i \phi (r_i^{k,\epsilon}) \leq \epsilon / 2^k$. Hence, $F_\phi = \bigcup_{k=1}^\infty F_k \subset \bigcup_{k=1}^\infty \bigcup_i B(z_i^{k,\epsilon}, r_i^{k,\epsilon})$ which gives $h_\phi (F_\phi) =0$ because
		$$
		h_\phi (F_\phi) \leq \sum_{k=1}^\infty \sum_i \phi(r_i^{k,\epsilon}) \leq \sum_{k=1}^\infty \frac{\epsilon}{2^k} = \epsilon .
		$$
		Moreover $\omega^{z_k} (F_\phi) \geq \omega^{z_k} (F_k) \geq C^{-1}$ for any $z_k$. As $\{z_k\}_{k=1}^\infty \subset \Omega$ is dense and $\omega^z (F_\phi)$ is $L_A$-harmonic with respect to $z$ (in particular continuous), then $\omega^p (F_\phi) \geq C^{-1}$ for any $p\in \Omega$. By \cite[Lemma 11.16]{Heinonen2006} we conclude $\omega^p (F_\phi) = 1$ for any $p\in \Omega$.
		
		Finally let
		$$
		F=\left\{ \xi \in \partial \Omega : \limsup_{r \to 0} \frac{\omega^p (B(\xi, r))}{r} > 0 \right\}.
		$$
		We claim that this set has $\sigma$-finite length and $\omega (F)=1$. We start by proving $\omega (F) = 1$, and later we will show that it has $\sigma$-finite length.
		
		Suppose to get a contradiction that $\omega^p (F)\not = 1$, that is,
		$$
		\omega^p (F^c) = \omega^p \left( \left\{ \xi\in \partial \Omega : \lim_{r\to 0} \frac{\omega^p (B(\xi, r))}{r}  = 0 \right\} \right) > 0 .
		$$
		Egorov's theorem ensures that for every $s > 0$ there exists a measurable set $V=V_s \subset F^c$ such that $\omega^p (V)<s$ and $\omega^p (B(\xi, r))/r \to 0$ as $r\to 0$ uniformly on $F^c \setminus V$. Since $0<\omega^p (F^c) = \omega^p (V) + \omega^p (F^c \setminus V) <s +\omega^p (F^c \setminus V)$, if we take $s>0$ small enough, say $s = \omega^p (F^c)/2$, then $\omega^p (F^c \setminus V) > 0$. The set $Y\coloneqq F^c \setminus V$ has non zero elliptic measure and the limit
		$$
		\lim_{r\to 0} \left\| \frac{\omega^p (B(\cdot, r))}{r} \right\|_{L^\infty (Y)} = 0 .
		$$
		
		Note that the function $\phi(r)\coloneqq \sup_{\xi \in Y} \omega^p (B(\xi, r))$ satisfies the conditions of the rate function in the beginning of this proof, and by definition it satisfies $\omega^p (B(\xi,r)) \leq \phi (r)$ for all $\xi \in Y$ and all $r>0$. All in all,
		\begin{enumerate}
			\item\label{phi prop 1} $\omega^p (B(\xi,r)) \leq \phi (r)$ for all $\xi \in Y$ and all $r>0$,
			\item $\phi$ is increasing, and
			\item $\phi (r) / r \to 0$ as $r\to 0$.
		\end{enumerate}
		For this particular function $\phi$, let $F_\phi$ be the set constructed in the beginning of this proof. That is, a set $F_\phi$ with $h_\phi (F_\phi)=0$ and $\omega^p (F_\phi)=1$ for every $p\in \Omega$. Hence $\omega^p (Y\cap F_\phi) = \omega^p (Y) > 0$, and moreover $h_\phi (Y\cap F_\phi) \leq h_\phi (F_\phi) = 0$. Consequently, we can cover $Y\cap F_\phi$ with balls $B(\xi_i,r_i)$ centered at $Y\cap F_\phi$ such that $\sum_i \phi (r_i) < \omega^p (Y\cap F_\phi) / 2$. With this we get
		$$
		\omega^p (Y\cap F_\phi ) \leq \sum_i \omega^p (B(\xi_i, r_i)) \overset{\text{\rf{phi prop 1}}}{\leq} \sum_i \phi (r_i) < \frac{\omega^p (Y\cap F_\phi )}{2}.
		$$
		This is a contradiction, and hence $\omega^p (F) = 1$.
		
		It remains to prove that $F$ has $\sigma$-finite one-dimensional Hausdorff measure. The set $F$ can be written as 
		$$
		F=\bigcup_{j\geq 1} F^{1/j} \text{ with } F^{1/j} \coloneqq \left\{ \xi \in \partial \Omega : \limsup_{r \to 0} \frac{\omega^p (B(\xi, r))}{r} > 1/j \right\} .
		$$
		Therefore, it suffices to see that every $F^{1/j}$ has finite one-dimensional Hausdorff measure. Each point $\xi \in F^{1/j}$ has arbitrarily small neighborhoods $B(\xi,r)$ such that $\omega^p (B(\xi,r)) \geq  r/j$. Given $\varepsilon>0$ small, consider the family of these balls centered at $F^{1/j}$ with radius at most $\varepsilon$, i.e.,
		$$
		\mathcal B_\varepsilon \coloneqq \{B(\xi, r) : \xi\in F^{1/j}, r<\varepsilon, \text{ and } \omega^p (B(\xi,r)) \geq  r/j\}.
		$$
		By the Besicovitch covering theorem there is a countable subfamily $\{B(\xi_i, r_i)\}_i \subset \mathcal B_\varepsilon$ such that no point belongs to more than a fixed finite number $C$ (it depends on the dimension only) of these balls. Hence,
		$$
		\sum_i r_i \leq j \sum_i \omega^p (B(\xi_i, r_i)) \leq Cj,
		$$
		and letting $\varepsilon \to 0$ we obtain $\HH^1 (F^{1/j}) \leq Cj$, as claimed.
	\end{proof}

Hence \cref{dimension elliptic measure} is proved under the assumption that \rf{key integral log bound} holds.

\section{\texorpdfstring{$L\log L(d\sigma)$}{LlogL} type estimate for small densities: Proof of \rf{key integral log bound}}\label{sec:log integral}

The purpose of this section is to prove \rf{key integral log bound}, under the hypothesis of \cref{covering elliptic measure}. More specifically, we prove the following result.

\begin{lemma}\label{key integral log bound absolute lemma}
	Let $\Omega \subset \R^2$ be a bounded $(\delta, r_0)$-Reifenberg flat domain, $p \in \Omega$ with $\dist (p, \partial \Omega)>r_0$, and $A$ be a real uniformly elliptic (not necessarily symmetric) matrix with ellipticity constant $\lambda$. Suppose also that $A$ is $\kappa$-Lipschitz in $U_{r_0} (\partial \Omega) \coloneqq \{x\in \R^2 : \dist (x, \partial \Omega) < r_0\}$ and that its symmetric part $A_0 = \frac{A+A^T}{2}$ is of the form $A_0=R^TBR$, with $R \in C^{0,1} (U_{r_0} (\partial \Omega))$ a rotation and $B \in C^{0,1} (U_{r_0} (\partial \Omega))$ diagonal.
	
	Then there exists $\delta_0 = \delta_0 (\lambda, \kappa \|A\|_{L^\infty (\R^2)}) >0$ and $C = C (\lambda, \kappa, r_0, \diam\, \partial\Omega) \in (0,\infty)$ such that if $\delta \leq \delta_0$, then for any $(\delta, r_0)$-Reifenberg flat domain $\widetilde \Omega \subset \Omega$ with smooth boundary $\partial \widetilde \Omega$ and small enough $\dist_\HH (\partial \Omega, \partial \widetilde \Omega)$, we have
	\begin{equation}\label{key integral log bound absolute}
		\left| \int_{\partial \widetilde \Omega} \log |S\nabla g_p^T (\xi)|^2 \, d\widetilde \omega^p (\xi) \right| \leq C < +\infty,
	\end{equation}
	where $\widetilde \omega^p=\omega^p_{\widetilde\Omega,A}$ is the elliptic measure in $\widetilde \Omega$ with respect to the matrix $A$, $g_p^T$ is the Green function in $\widetilde \Omega$ with respect to the matrix $A^T$, and $S$ is the square root matrix of the symmetric part $A_0=(A+A^T)/2$, i.e., $S^TS=A_0$.
\end{lemma}

\begin{rem}
    As argued from \rf{lower bound kernel and log gradient} to \rf{integral log kernel bounded}, this implies that the Radon-Nykodym derivative $\frac{d\widetilde\omega^p}{d\sigma}$ satisfies the following $L\log L (d\sigma)$ type estimate
    $$
    -\infty < C^\prime (\lambda,\kappa,\diam\, \partial\Omega) \leq \int_{\partial\widetilde \Omega} \frac{d\widetilde\omega^p}{d\sigma}(\xi)\log \frac{d\widetilde \omega^p}{d\sigma}(\xi) \, d\sigma(\xi).
    $$
\end{rem}

By symmetry, estimate \rf{key integral log bound absolute} is equivalent to the existence of a constant $C=C(\lambda, \kappa,\diam\partial\Omega)$ depending only on the ellipticity constant, the Lipschitz seminorm of the matrix $A$, and the diameter $\diam \, \partial \Omega$ (but not on $\widetilde\Omega$) such that
\begin{equation}\label{key integral log bound absolute symmetric}
\left| \int_{\partial \widetilde \Omega} \log |S\nabla g (\xi)|^2 \, d\widetilde\omega_T^p (\xi) \right| \leq C < \infty ,
\end{equation}
where $S=A_0^{1/2}$, i.e., $S^T S=A_0$, $g=g_p$ is the Green function in $\widetilde \Omega$ with respect to the matrix $A$ with pole $p$, and $\widetilde \omega_T = \omega_{\widetilde\Omega,A^T}^p$ denotes the elliptic measure in $\widetilde\Omega$ with respect to the matrix $A^T$ and a pole $p\in \widetilde \Omega$ such that $\dist(p,\partial \Omega)>r_0$. The existence of the matrix $S$ is granted by the fact that the symmetric matrix $A_0$ is uniformly elliptic with the same ellipticity constant as $A$, and hence positive definite. In fact, $S=R^T\sqrt{B}R$, where $\sqrt{B}=(\delta_{i=j} \sqrt{b_{ij}})_{1\leq i,j\leq2}$.

Throughout all this section, when dealing with terms within solid integrals we will use the following notation:
\begin{itemize}
	\item We write $h(x) = \OO(f(x))$ to denote ``$|h(x)| \leq C f(x)$ almost everywhere with respect to the Lebesgue measure''.
	
	\item Given a Lipschitz function $h$, we write $|\nabla h(x)| \leq C$ instead of ``$|\nabla h(x)| \leq C$ almost everywhere with respect to the Lebesgue measure''. Recall that Lipschitz functions are differentiable almost everywhere by Rademacher's theorem, see \cite[Theorem 7.3]{Mattila1995}.
\end{itemize}

\subsection{Directional derivatives and the dual space}\label{sec:set up and notation}

We introduce some extra notation regarding the directional derivatives.

\begin{notation}
	Definition of directional derivatives.
	\begin{itemize}
		\item The (vertical) vector $e_i \coloneqq (
		0,  \ldots , 0 , \overset{\ithposition}{1} , 0 , \ldots , 0 )^T$ has $1$ in position $i$ and $0$'s otherwise. In $\R^2$ there are only two such vectors, namely $e_1 = (1,0)^T$ and $e_2=(0,1)^T$.
		
		\item $R_i \coloneqq e_i^T \cdot R$ corresponds to the $i$-th row of $R$.
		
		\item The $\partial^i$-directional derivative $\partial^i \coloneqq \partial_{R^T_i}$ (superscript) is defined as
		$$
		\partial^i  f \coloneqq \partial_{R^T_i} f = \langle \nabla f , R^T_i \rangle = R_i \cdot \nabla f .
		$$
		Here $\langle \cdot , \cdot \rangle$ denotes the standard scalar product in $\R^{n+1}$. The derivative $\partial_i = \partial_{e_i}$ (subscript) is the usual one in the direction $e_i$.
		
		\item The $R$-directional gradient $\nabla_R$ is defined as
		$$
		\nabla_R f \coloneqq R\cdot \nabla f=\left( \partial^i f \right)_{i=1}^{n+1} .
		$$
	\end{itemize}
\end{notation}	

\begin{rem}
	The directional derivative $\partial^i$ preserves the usual properties in sums, products, and logarithms, i.e., $\partial^i (f+g) = \partial^i f + \partial^i g$, $\partial^i (fg)= g\partial^i f  + f \partial^i g$, and $\partial^i \log f = \frac{\partial^i f}{f}$.
\end{rem}

Note that
\begin{equation}\label{two directional derivatives}
	\left| \partial^{ij} f \right| = 
	\left| \partial^{i} (R_j \cdot \nabla f) \right| 
	= \left| \partial^i R_j \cdot \nabla f + R_j \partial^i \nabla f \right|
	\lesssim |\nabla f| + \left| \nabla^2 f \right|,
\end{equation}
where we denote $|\nabla^2 f|^2 \coloneqq \sum_{i,j=1}^{n+1} (\partial_{i,j} f)^2$.

Since the directions $R_i$ change depending on the point, the usual integration by parts formula does not apply. Instead we have the following formula.

\begin{claim}[Integration by parts formula]
	Let $U\subset \R^{n+1}$ be an open set. For $f, h \in W^{1,2}(U)$, if $fh \in W_0^{1,2} (U)$ then
	\begin{equation}\label{integration parts special}
		\int_{U} \partial^i f(x) h(x) \, dx = -\int_{U} f(x)h(x)\divv R^T_i (x) \, dx- \int_{U} f(x) \partial^i h(x) \, dx .
	\end{equation}
	\begin{proof}
		The product rule gives $\partial^i f(x) h(x)  = \partial^i (f(x) h(x))  - f(x) \partial^i h(x)$. By definition of the directional derivative $\partial^i$ we can write
		\begin{align*}
			\partial^i (f(x) h(x)) &=  R_i (x) \cdot \nabla (f(x) h(x))  = \sum_{j=1}^{n+1} \partial_j (f(x) h(x)) R_{i,j} (x) \\
			& =\sum_{j=1}^{n+1} \partial_j \left( f(x) h(x) R_{i,j}  (x) \right) 
			-f(x) h(x) \sum_{j=1}^{n+1} \partial_j R_{i,j} (x) \\
			&=\divv (fhR^T_i)(x) - f(x)h(x) \divv R^T_i (x) ,
		\end{align*}
		which imply
		$$
		\int \partial^i f(x) h(x) \, dx = \int  \divv (fhR^T_i)(x) \, dx - \int f(x)h(x) \divv R^T_i (x)  \, dx - \int f(x) \partial^i h(x) \, dx .
		$$
		Using that $C^{\infty}_c(U)$ is dense in $W_0^{1,2} (U)$ and the divergence theorem, the first element in the right-hand side is $\int  \divv (fhR^T_i)(x) \, dx = 0$.
	\end{proof}
\end{claim}

Again, since the directions $R$ depend on the point, the directional derivatives do not commute. Instead, the following formula relating $\partial^{\alpha, \beta}$ and $\partial^{\beta,\alpha}$ is available, where $\partial^{\alpha, \beta} \coloneqq \partial^\alpha \partial^\beta$. We refer to this as the `almost'-commutative property.

\begin{claim}\label{directional derivatives does not commute}
Let $U\subset \R^{n+1}$ be an open set and $f\in W^{2,2} (U)$. For $\alpha, \beta \in \{1,\ldots, n+1\}$,
\begin{equation*}
\partial^{\alpha , \beta} f - \partial^{\beta,\alpha} f = (\partial^\alpha R_\beta - \partial^\beta R_\alpha)  \cdot \nabla f  .
\end{equation*}
\begin{proof}
By definition we have $\partial^\beta f =  R_\beta \cdot \nabla f = \sum_{i=1}^{n+1} \partial_i f R_{\beta,i}$. Applying $\partial^\alpha$ gives
$$
\partial^{\alpha, \beta}  f =  R_\alpha \cdot \nabla \left(  \sum_{i=1}^{n+1} \partial_i f R_{\beta,i}   \right) 
= \sum_{j=1}^{n+1} \partial_j \left( \sum_{i=1}^{n+1} \partial_i f R_{\beta,i} \right) R_{\alpha,j} .
$$
Expanding the derivative of this last expression and arranging we obtain
$$
\partial^{\alpha, \beta}  f = \sum_{i,j = 1}^{n+1} \partial_{j,i} f R_{\beta,i} R_{\alpha,j} + \sum_{i=1}^{n+1} \partial_i f \sum_{j=1}^{n+1} \partial_j R_{\beta,i} R_{\alpha,j}.
$$
The sum inside the second term in the right-hand side is precisely $\sum_{j=1}^{n+1} \partial_j R_{\beta,i} R_{\alpha,j} = R_\alpha \cdot \nabla R_{\beta,i} = \partial^\alpha R_{\beta,i}$. Hence
$$
\partial^{\alpha, \beta}  f = \sum_{i,j = 1}^{n+1} \partial_{j,i} f R_{\beta,i} R_{\alpha,j} 
+ \sum_{i=1}^{n+1} \partial_i f \partial^\alpha R_{\beta,i}
 =  \sum_{i,j = 1}^{n+1} \partial_{j,i} f R_{\beta,i} R_{\alpha,j}  
 + \partial^\alpha R_\beta \cdot \nabla f .
$$

By symmetry we obtain $\partial^{\beta, \alpha}  f = \sum_{i,j = 1}^{n+1} \partial_{j,i} f R_{\alpha,i} R_{\beta,j} +  \partial^\beta R_\alpha \cdot \nabla f$, and subtracting $\partial^{\alpha, \beta} f - \partial^{\beta, \alpha} f$ we get \cref{directional derivatives does not commute}.
\end{proof}
\end{claim}

In particular, setting $\alpha=1$ and $\beta=2$ in the planar case we have
\begin{equation*}
	\partial^{1,2} f - \partial^{2,1} f =  (\partial^1 R_2 - \partial^2 R_1) \cdot \nabla f ,
\end{equation*}
and applying it to the Green function $g$, we get
\begin{equation}\label{derivative_commutator}
	\partial^{1,2} g - \partial^{2,1} g = (\partial^1 R_2 - \partial^2 R_1) \cdot  \nabla g = \OO(|\nabla g|).
\end{equation}

\subsubsection{Definition of the dual space}

Let $U\subset \widetilde\Omega$ and $u \in L^1_{\loc} (U)$. Define the linear functional 
$$
T_u (\psi) \coloneqq \int_{U} u\psi,
$$
whenever it makes sense. In particular, $T_u\in L^{\infty}_c (U)^\prime$. To simplify the notation, we will also write $u(\psi)$ to denote $T_u (\psi)$.

The $\partial_i$-derivative functional is defined as
$$
T_{\partial_i u}(\psi) \coloneqq -\int_U u\partial_i \psi,
$$
whenever it makes sense. Note that $T_{\partial_i u} \in W^{1,\infty}_c (U)^\prime$ for $u\in L^1_{\loc}(U)$. Moreover, we also have $T_{\partial_i u} \in W^{1,2}_c (U)^\prime$ whenever $u\in L^2_{\loc} (U)$. As before, we will also write $D u(\psi)$ to denote $T_{D u} (\psi)$ for any differential operator $D$. In general, the functional
$$
T_{fD(u)}(\psi) \coloneqq T_{D(u)} (f\psi),
$$
is in $W^{1,\infty}_c (U)^\prime$ as long as either $f\in W^{1,2}_{\loc}(U)$ and $u\in L^2_{\loc}(U)$, or $f\in W^{1,\infty}_{\loc} (U)$ and $u\in L^1_{\loc} (U)$.

We claim that whenever $u\in L^2_{\loc}(U)$ and $f\in W^{1,2}_{\loc}(U)$, in the dual space $W^{1,\infty}_c (U)^\prime$, the $\partial_i$-derivative functional satisfies the product rule
\begin{equation}\label{product rule distributional final}
	T_{\partial_i (f u)} = T_{\partial_i f u} + T_{f \partial_i u}.
\end{equation}
Indeed, for any $\psi \in W_c^{1,\infty}(U)$, the previous equality reads as
$$
-\int_U fu\partial_i \psi = \int_U \partial_i f u \psi - \int_U u \partial_i (f\psi),
$$
which holds by the Leibniz's rule a.e.\ for $W^{1,2}$-functions.

Using the definitions of $T_{\partial_i}$ and the directional derivative $\partial^i (\cdot) = R_i \nabla (\cdot)$, see the beginning of \cref{sec:set up and notation}, the directional $\partial^i$-derivative functional $T_{\partial^i u} : W^{1,\infty}_c (U) \to \R$, for $u\in L^1_{\loc}(U)$, is defined as
\begin{equation}\label{def:distributional directional derivative}
	T_{\partial^i u} \coloneqq T_{R_i \nabla u}= \sum_j T_{R_{i,j}\partial_j u} .
\end{equation}
Equivalently, for $\psi\in W^{1,\infty}_c (U)$,
\begin{equation}\label{def2:distributional directional derivative}
	\begin{aligned}
		T_{\partial^i u}(\psi) &= \sum_j T_{\partial_j u} (R_{i,j} \psi)
		= - \sum_j \int_{U} u\partial_j (R_{i,j} \psi)
		= -\sum_j \int_{U} u\psi \partial_j R_{i,j} - \sum_j \int_{U} u R_{i,j} \partial_j \psi \\
		&= - \int_{U} u\psi \divv R_i^T - \int_{U} u \partial^i \psi .
	\end{aligned}
\end{equation}
This agrees with the integration by parts formula in \rf{integration parts special} when $u\in W^{1,2} (U)$.

Next we claim that, in the dual space $W_c^{1,\infty} (U)^\prime$, for $u\in L^2_{\loc}(U)$ and $f\in W^{1,2}_{\loc}(U)$, the directional $\partial^i$-derivative functional satisfies the product rule
\begin{equation}\label{product rule distributional directional derivative}
	T_{\partial^i (f u)} = T_{\partial^i f u} + T_{f \partial^i u}.
\end{equation}
Indeed, $R_{i,j}\psi \in W_c^{1,\infty}(U)$ when $\psi\in W^{1,\infty}_c (U)$, and hence, for $f \in W^{1,2}_{\loc} (U)$,
$$
T_{\partial^i (fu)} (\psi) \overset{\text{\rf{def:distributional directional derivative}}}{=} \sum_j T_{\partial_j (fu)} (R_{i,j} \psi)
\overset{\text{\rf{product rule distributional final}}}{=}
\sum_j T_{\partial_j f u} (R_{i,j} \psi) + T_{f\partial_j u} (R_{i,j} \psi) 
\overset{\text{\rf{def:distributional directional derivative}}}{=}
T_{\partial^i fu} (\psi) + T_{f\partial^i u} (\psi) ,
$$
as claimed.

Another important property is the following `almost'-commutative property of the directional derivatives (compare to \cref{directional derivatives does not commute}). For $u\in W^{1,1}_{\loc} (U)$, in the dual space $W^{1,\infty}_c (U)^\prime$, we claim that
\begin{equation}\label{directional derivatives does not commute distribution}
	T_{\partial^{1} (\partial^2 u)} - T_{\partial^2 (\partial^1 u)} = T_{(\partial^1 R_2 - \partial^2 R_1) \cdot \nabla u}.
\end{equation}
Indeed, for $\psi\in W^{1,\infty}_c(U)$ and $\alpha,\beta\in \{1,2\}$,
$$
T_{\partial^\alpha (\partial^\beta u)} (\psi) \overset{\text{\rf{def2:distributional directional derivative}}}{=} -\int_U \partial^\beta u \psi \divv R_\alpha^T - \int_U \partial^\beta u \partial^\alpha \psi .
$$
Given $\{u_k\}_{k\geq 1} \subset C_c^\infty (U)$ with $\lim_{k\to \infty}\|u_k - u\|_{W^{1,1}(\supp \psi)} =0$, this equals
\begin{equation}\label{distributional alpha-beta derivative}
\begin{aligned}
T_{\partial^\alpha (\partial^\beta u)} (\psi)
=&-\int_U \partial^\beta u_k \psi \divv R_\alpha^T - \int_U \partial^\beta u_k \partial^\alpha \psi \\
&- \int_U \partial^\beta (u-u_k) \psi \divv R_\alpha^T - \int_U \partial^\beta (u - u_k) \partial^\alpha \psi .
\end{aligned}
\end{equation}
By the integration by parts in \rf{integration parts special} applied to $\partial^\beta u_k$, we have that the first row in the right-hand side is precisely $\int_U \partial^\alpha (\partial^\beta u_k) \psi$. Applying now the `almost'-commutative property (\cref{directional derivatives does not commute}) here we obtain
\begin{align*}
	T_{\partial^\alpha (\partial^\beta u)} (\psi)
	=&\int_U \partial^\beta (\partial^\alpha u_k) \psi 
	+ \int_U ((\partial^\alpha R_\beta - \partial^\beta R_\alpha) \cdot \nabla u_k )\psi\\
	&- \int_U \partial^\beta (u-u_k) \psi \divv R_\alpha^T - \int_U \partial^\beta (u - u_k) \partial^\alpha \psi .
\end{align*}
Again, by the integration by parts in \rf{integration parts special}, we can replace the first term in the right-hand side to obtain
\begin{equation*}
\begin{aligned}
	T_{\partial^\alpha (\partial^\beta u)} (\psi)
	=&-\int_U \partial^\alpha u_k \psi \divv R_\beta^T - \int_U \partial^\alpha u_k \partial^\beta \psi \\
	&+ \int_U ((\partial^\alpha R_\beta - \partial^\beta R_\alpha) \cdot \nabla u_k) \psi \\
	&- \int_U \partial^\beta (u-u_k) \psi \divv R_\alpha^T - \int_U \partial^\beta (u - u_k) \partial^\alpha \psi .
\end{aligned}
\end{equation*}
Now, adding and subtracting $u$ in the second row we get
\begin{equation*}
\begin{aligned}
	T_{\partial^\alpha (\partial^\beta u)} (\psi)
	=&-\int_U \partial^\alpha u_k \psi \divv R_\beta^T - \int_U \partial^\alpha u_k \partial^\beta \psi \\
	& + T_{(\partial^\alpha R_{\beta} - \partial^\beta R_{\alpha}) \nabla u} (\psi) 
	+ \int_U ((\partial^\alpha R_\beta - \partial^\beta R_\alpha) \cdot \nabla (u_k-u)) \psi, \\
	&- \int_U \partial^\beta (u-u_k) \psi \divv R_\alpha^T - \int_U \partial^\beta (u - u_k) \partial^\alpha \psi .
\end{aligned}
\end{equation*}
By symmetry in \rf{distributional alpha-beta derivative} and subtracting we obtain
\begin{equation*}
\begin{aligned}
	T_{\partial^\alpha (\partial^\beta u)} (\psi) - T_{\partial^\beta (\partial^\alpha u)} (\psi) 
	=&\ T_{(\partial^\alpha R_{\beta} - \partial^\beta R_{\alpha}) \nabla u} (\psi) 
	+ \int_U ((\partial^\alpha R_\beta - \partial^\beta R_\alpha) \cdot \nabla (u_k-u)) \psi, \\
	&- \int_U \partial^\beta (u-u_k) \psi \divv R_\alpha^T - \int_U \partial^\beta (u - u_k) \partial^\alpha \psi \\
	&+ \int_U \partial^\alpha (u-u_k) \psi \divv R_\beta^T + \int_U \partial^\alpha (u - u_k) \partial^\beta \psi.
\end{aligned}
\end{equation*}
By the convergence $\lim_{k\to \infty}\|u_k - u\|_{W^{1,1}(\supp \psi)} =0$ we get \rf{directional derivatives does not commute distribution}.

\subsubsection{Properties of directional derivatives and the dual space}

The following claim allows us to move from the initial matrix $A$ to its symmetric part $A_0$, which we will relate to the directional derivatives in \cref{directional_div}.

\begin{claim}\label{div A and div Asim}
	Let $U\subset \widetilde \Omega$ be an open set and $u\in W_{\loc}^{1,2} (U)$. Every $\psi \in W^{1,\infty}_c (U)$ satisfies
	$$
	\left(\divv A\nabla u\right)(\psi) = \left(\divv A_0 \nabla u\right)(\psi) + \left(\sum_{i,j=1}^{n+1} \partial_i \left(\frac{a_{i,j}-a_{j,i}}{2}\right) \partial_j u\right)(\psi).
	$$
	\begin{proof}
		Let $u\in C^\infty (U)$. By definition of the divergence and differentiating we have
		\begin{equation*}
		\divv A\nabla u =  \sum_{i,j=1}^{n+1} \partial_i  a_{i,j} \partial_j u  + \sum_{i,j=1}^{n+1} a_{i,j} \partial_{i,j} u 
		=  \sum_{i,j=1}^{n+1} \partial_i  a_{i,j} \partial_j u  
		+ \sum_{i,j=1}^{n+1} \frac{a_{i,j}+a_{j,i}}{2} \partial_{i,j} u .
		\end{equation*}
		Note that $a^0_{i,j} \coloneqq \frac{a_{i,j}+a_{j,i}}{2}$ are precisely the coefficients of the symmetric matrix $A_0$. The product derivative rule gives
		$$
		a^0_{i,j} \partial_{i,j} u =  \partial_i \left(a^0_{i,j} \partial_j u \right) - \partial_i  a^0_{i,j}  \partial_j u ,
		$$
		and applying this relation to each pair of indexes $i,j \in \{1,\ldots, n+1\}$ we get
		\begin{equation*}
		\sum_{i,j=1}^{n+1} a^0_{i,j} \partial_{i,j} u = \divv A_0 \nabla u - \sum_{i,j=1}^{n+1} \partial_i  a^0_{i,j}  \partial_j u .
		\end{equation*}
		Summing up,
		\begin{align*}
			\divv A\nabla u
			&=  \divv A_0 \nabla u + \sum_{i,j=1}^{n+1} \partial_i  a_{i,j} \partial_j u - \sum_{i,j=1}^{n+1} \partial_i \left( \frac{a_{i,j} +a_{j,i}}{2} \right) \partial_j u \\
			&= \divv A_0 \nabla u + \sum_{i,j=1}^{n+1} \partial_i \left(\frac{a_{i,j}-a_{j,i}}{2}\right) \partial_j u,
		\end{align*}
	as claimed, for functions in $C^\infty(U)$.
	
	Let us check this in the dual sense for a function $u\in W^{1,2}_{\loc} (U)$. Let $\psi \in W^{1,\infty}_c(U)$, and hence $u\in W^{1,2} (\supp \psi)$. Take $\{u_k\}_{k\geq 1} \subset C^\infty_c (U)$ such that $\lim_{k\to \infty} \|u_k - u\|_{W^{1,2}(\supp \psi)} = 0$. As we have the claim for each function $u_k$, in particular
	\begin{align*}
		\int_{U} A \nabla u \nabla \psi 
		=&  \int_{U} A \nabla (u-u_k) \nabla \psi  + \int_{U} A \nabla u_k \nabla \psi \\
		=& \int_{U} A \nabla (u-u_k) \nabla \psi  + \int_{U} A_0 \nabla u_k \nabla \psi - \int_{U} \sum_{i,j=1}^{n+1} \partial_i \left(\frac{a_{i,j}-a_{j,i}}{2}\right) \partial_j u_k \cdot \psi \\
		=& \int_{U} A_0 \nabla u \nabla \psi 
		- \int_{U} \sum_{i,j=1}^{n+1} \partial_i \left(\frac{a_{i,j}-a_{j,i}}{2}\right) \partial_j u \cdot \psi \\
		&+ \int_{U} A \nabla (u-u_k) \nabla \psi  
		+ \int_{U} A_0 \nabla (u_k-u) \nabla \psi 
		- \int_{U} \sum_{i,j=1}^{n+1} \partial_i \left(\frac{a_{i,j}-a_{j,i}}{2}\right) \partial_j (u_k-u) \psi .
	\end{align*}
	Now, since the matrix is Lipschitz and $\lim_{k\to \infty} \|u_k - u\|_{W^{1,2}(\supp\psi)} = 0$, we conclude
	$$
	\int_{U} A \nabla u \nabla \psi =  \int_{U} A_0 \nabla u \nabla \psi 
	- \int_{U} \sum_{i,j=1}^{n+1} \partial_i \left(\frac{a_{i,j}-a_{j,i}}{2}\right) \partial_j u \cdot \psi ,
	$$
	as claimed.
	\end{proof}
\end{claim}

We treat one of the main terms in \rf{key integral log bound absolute} by means of a perturbation argument. To do that, we note that we can write $\divv (A_0 \nabla u)$ as in the diagonal case using the $R$-directional derivatives, plus an error term. More precisely, we have the following claim.

\begin{claim}\label{directional_div}
	Let $A_0=R^T B R$ satisfy the conditions in \cref{key integral log bound absolute lemma}, let $U\subset \widetilde \Omega$ be an open set and $u\in W^{1,2}_{\loc} (U)$. Every $\psi\in W^{1,\infty}_c (U)$ satisfies
	$$
	\left(\divv (A_0 \nabla u)\right) (\psi)
	= \left(\sum_{i=1}^{n+1} \partial^i (b_i \partial^i u)\right) (\psi)
	+ \left(\sum_{i=1}^{n+1} b_i \partial^i u \cdot \divv R^T_i\right) (\psi).
	$$
	\begin{proof}
		Assume first $u\in C^\infty (U)$. Let us write $R\nabla u = ( R_j \cdot \nabla u )_{j=1}^{n+1} = ( \partial^j u )_{j=1}^{n+1}$, and so $BR\nabla u= ( b_j \partial^j u )_{j=1}^{n+1}$. Hence,
		$$
		R^T BR\nabla u = R^T ( b_j \partial^j u )_{j=1}^{n+1} = \left( \sum_{j=1}^{n+1} R^T_{i,j} b_j \partial^j u \right)_{i=1}^{n+1} .
		$$
		Therefore,
		\begin{equation*}
			\divv (A_0 \nabla u) = \sum_{i=1}^{n+1} \partial_i \left( \sum_{j=1}^{n+1} R^T_{i,j} b_j \partial^j u \right)   = \sum_{j=1}^{n+1} b_j \partial^j u \sum_{i=1}^{n+1} \partial_i (R_{j,i}) + \sum_{j=1}^{n+1} \sum_{i=1}^{n+1} R_{j,i} \partial_i (b_j \partial^j u) .
		\end{equation*}
		Note that the sum inside the first term in the right-hand side is precisely $\sum_{i=1}^{n+1} \partial_i (R_{j,i}) = \divv R^T_j$, and the sum inside the second term in the right-hand side is $\sum_{i=1}^{n+1} R_{j,i} \partial_i (b_j \partial^j u) = R_j \cdot \nabla (b_j \partial^j u) = \partial^j (b_j \partial^j u)$. Thus, the claim follows for $u\in C^\infty (U)$.
		
		From the definition of the directional $\partial^i$-derivative in \rf{def:distributional directional derivative} and \rf{def2:distributional directional derivative}, the claim follows by a density argument as in \cref{div A and div Asim}. 
	\end{proof}
\end{claim}

\subsection{Sketch of the proof}\label{sec:sketch of the proof}

Throughout this section the pole $p\in \Omega$ so that $\dist(p,\partial \Omega)>r_0$ is fixed unless it is otherwise stated (see \cref{bound easy terms 1 2} below). In any case, it will be far from $\partial \Omega$ and so from $\partial \widetilde \Omega$. Recall also that $S\coloneqq A_0^{1/2}$, i.e., $S^T S = A_0$.

Fix $R\coloneqq \min \{1, r_0 / 2 \} / 2$, and so we have $0<R < \min \{1, r_0 / 2 \}$. By \cite[Lemma 3.35]{Lewis2008}, see \cref{comparability function gradient} and \cref{comparability function gradient Lipschitz} above, if $\widetilde \Omega$ is Reifenberg flat enough, depending on the ellipticity constant $\lambda$ and the value $\kappa \|A\|_{L^\infty (\R^{n+1})}$, then there exists a constant $c = c (\lambda, \kappa \|A\|_{L^\infty (\R^{n+1})})\geq 1$ such that
\begin{equation*}
c^{-1} |\nabla g (y)| \leq \frac{g (y)}{\dist \left(y, \partial \widetilde \Omega\right)} \leq c |\nabla g (y)|  \text{ for every }y \in B(\xi, R/ c) \text{ with } \xi \in \partial \widetilde \Omega .
\end{equation*}
Here we need to work with $(\delta, r_0)$-Reifenberg flat domains with flatness parameter $\delta$ smaller than some constant depending on the ellipticity constant $\lambda$ and the value $\kappa \|A\|_{L^\infty (\R^2)}$.

Let us now fix the support function $\varphi$ with $\varphi (p)=0$:

\begin{rem}\label{choice of the support function}
	Let $\varphi \in C_c^\infty (\R^2)$ with
	\begin{itemize}
		\item $\varphi = 1$ in $U_{\frac{R}{2000c}} (\partial \Omega)\coloneqq \{x\in \R^2 : \dist(x,\partial \Omega)\leq \frac{R}{2000c}\}$,
		\item $\varphi = 0$ in $U_{\frac{R}{1500c}} (\partial \Omega)^c$, and
		\item $|\nabla \varphi|\lesssim 1$.
	\end{itemize}
	So $\supp \varphi \subseteq U_{\frac{R}{1500c}} (\partial \Omega)$. Note that $\partial \widetilde \Omega \subset \Omega \cap \supp \varphi$ if we assume that $\dist_\HH (\partial \widetilde \Omega, \partial \Omega)$ is small enough. With this choice of $\varphi$ we have the comparability
	\begin{equation}\label{key comparability inside the support}
	|\nabla g(y)| \approx \frac{g(y)}{\dist(y, \partial \widetilde \Omega)} \text{ for }y\in \widetilde\Omega \cap U_{\frac{R}{1500c}} (\partial\Omega) \supset \widetilde\Omega\cap \supp \varphi.
	\end{equation}
\end{rem}

We claim
\begin{equation}\label{log gradient g is sobolev 1-2 loc}
	\log |S\nabla g|^2 \in W^{1,2}_{\loc} (\widetilde\Omega \cap U_{\frac{R}{1500c}} (\partial\Omega)).
\end{equation}
Indeed, this follows as $|\nabla g| \approx g/\dist(\cdot, \partial \widetilde \Omega) >0$ in $\widetilde\Omega \cap U_{\frac{R}{1500c}} (\partial\Omega)$ and $g\in W^{2,2} (\widetilde\Omega \cap U_{\frac{R}{1500c}} (\partial\Omega))$, see \rf{gradient log_a gradient g}.

\begin{notation}
	From now on the variables and the region of integration will not be written unless they are not clear from the context.
\end{notation}

Estimate \rf{key integral log bound absolute symmetric} will follow from the following partial results.

For $a\geq 10$, let $\log_{(a)} x \coloneqq \max \{\log x, -a\}$. Note that in the weak sense $(\log_{(a)} (x))' = \characteristic_{\{x\geq e^{-a}\}}(x) \log x$.

\begin{lemma}[Step 1]\label{integral elliptic measure inicial and truncated}
	For $a\geq 10$ big enough,
	$$
	\int_{\partial\widetilde\Omega} \log |S\nabla g|^2 \, d\widetilde\omega_T^p 
	= \int_{\partial\widetilde\Omega} \log_{(a)} |S\nabla g|^2 \, d\widetilde\omega_T^p 
	+ \sigma(\partial \widetilde\Omega) \OO(ae^{-a/2}).
	$$
	Since $\sigma(\partial\widetilde\Omega)<\infty$, the second term in the right-hand side tends to zero as $a\to \infty$.
\end{lemma}

We prove the preceding lemma in \cref{section:link lemmas}.

\begin{lemma}[Step 2]\label{bound easy terms 1 2}%lemma 1
	For any $a\geq 10$ so that $e^{-a}<\min_{z\in\supp\nabla \varphi} |S\nabla g (z)|^2$, and for a.e.\ $p\in\Omega$ with $\dist(p,\partial\Omega)>r_0$,
	\begin{equation*}
			\int_{\partial \widetilde \Omega} \log_{(a)} |S\nabla g|^2 \, d\widetilde \omega_T^p = - \int_{\widetilde\Omega} \langle  A^T \nabla \log_{(a)} |S\nabla g|^2 , \nabla (\varphi g) \rangle
			+ \OO\left(1+\int_{\supp \varphi} \frac{|\nabla^2 g|}{|\nabla g|} g\right).
	\end{equation*}
\end{lemma}

We prove the preceding lemma in \cref{proof of bound easy terms 1 2}. Note that the left-hand side integral is supported on the boundary, while the one on the right-hand side is a solid integral.

Next, fix $\widetilde\Omega_\varphi \subset \widetilde\Omega$ a smooth domain satisfying
$$
\widetilde\Omega\cap \supp \varphi \subset \widetilde\Omega_\varphi \subset U_{\frac{R}{1000c}} (\partial \Omega) \cap \widetilde \Omega .
$$
With this choice we have the comparability $|\nabla g| \approx g/\dist(\cdot, \partial\widetilde\Omega)$ in $\widetilde\Omega_\varphi$, see \rf{key comparability inside the support}. We prove the following two lemmas in \cref{section:link lemmas}.

\begin{lemma}[Step 3]\label{truncated solid integral to non truncated}
	For $a\geq 10$ big enough,
	\begin{align*}
	- \int_{\widetilde\Omega} \langle A^T \nabla \log_{(a)} |S\nabla g|^2, \nabla (\varphi g)\rangle  
	= &  - \int_{\widetilde\Omega} \langle A^T \nabla \log |S\nabla g|^2, \nabla (\varphi g)\rangle  \\
	& + \OO\left(
	\HH^2 (\Omega) e^{-a/2} + \int_{\widetilde\Omega\cap\supp\varphi\cap\{|S\nabla g|^2 \leq e^{-a}\}} |\nabla^2 g|
	\right).
	\end{align*}
	Note that the last term in the right-hand side also tends to zero as $a \to \infty$ because $g \in W^{2,2} (\widetilde\Omega_\varphi)$ and $\bigcap_{a\geq 10} \{|S\nabla g|^2 \leq e^{-a}\} \cap \widetilde\Omega_\varphi = \emptyset$.
\end{lemma}

For $\varepsilon >0$ as small as desired, we consider a given function $\psi_\varepsilon \in C^\infty (\R^2)$ satisfying
\begin{enumerate}
	\item $0\leq \psi_\varepsilon \leq 1$ everywhere,
	\item $\psi_\varepsilon = 0$ in $U_\varepsilon (\partial \widetilde\Omega)$,
	\item $\psi_\varepsilon = 1$ in $U_{3\varepsilon} (\partial \widetilde\Omega)^c$, and
	\item $|\nabla \psi_\varepsilon| \lesssim \frac{1}{\varepsilon}$.
\end{enumerate}

\begin{lemma}[Step 4]\label{lemma close to the boundary}
	For $\psi_\varepsilon$ as above we have 
	\begin{equation*}
		\int_{\widetilde\Omega} \langle A^T \nabla \log |S\nabla g|^2, \nabla (\varphi g)\rangle
		= \lim_{\varepsilon \to 0} \left\{
		\int_{\widetilde\Omega} \langle A^T \nabla \log |S\nabla g|^2, \nabla (\psi_\varepsilon\varphi g)\rangle 
		+\OO\left(\frac{1}{\varepsilon}\int_{U_{3\varepsilon} (\partial\widetilde\Omega)} \frac{|\nabla^2 g|}{|\nabla g|}g\right)
		\right\}.
	\end{equation*}
\end{lemma}

In the weak sense, we can write the first term in the right-hand side in \cref{lemma close to the boundary} as $\divv (A^T \nabla \log |S\nabla g|^2)$ acting on $\psi_\varepsilon \varphi g \in W_c^{1,\infty}(\widetilde\Omega_\varphi)$, which can be understood as test functions. When studying the harmonic measure in the plane, i.e., $A=Id$, this term is $0$ by the harmonicity of $\log |\nabla g|$ far from the critical points, a property that does not hold in general in higher dimensions. This was a key point to establish \rf{key integral log bound} for the Laplacian in \cite[Lemma 3.1]{Jones1988}.

In the following remark we discuss this argument in the constant coefficient case.

\begin{rem}[Constant matrix and $L_A$-harmonic functions]\label{rem:constant matrix and LA harmonic functions}
	Suppose now that the matrix $A$ is constant. Given any function $u$, by \cref{div A and div Asim},
	$$
	\divv (A^T\nabla u) = \divv (A_0 \nabla u ) ,
	$$
	i.e., we can reduce the study to the symmetric part. Moreover, by means of a linear change of variables, see \cref{linear change of variables}, we have that $\divv (A^T \nabla u)=0$ if and only if $\widetilde u = u \circ D$ satisfies $\divv \left( \widetilde {A_0} \nabla \widetilde u \right)=0$ where $\widetilde {A_0} = D^{-1} A_0 (D^{-1})^T$, for any constant matrix $D$ with $\det D \not = 0$. Since $A$ is constant now, if we choose $D=S^T$ where $S=A_0^{1/2}$, i.e., $S^T S=A_0$, then $\widetilde {A_0} = (S^T)^{-1} A_0 S^{-1} = (S^T)^{-1} S^T S S^{-1} = Id$, and so $\Delta \widetilde u = 0$. It is known that if $\Delta \widetilde u=0$, then $\Delta \log |\nabla \widetilde u|^2 =0$ (only) in the plane, whenever $\nabla \widetilde u\not = 0$. If we undo the previous change of variables then we obtain
	$$
	\divv \left( A_0 \nabla \left( \log \left| S \nabla u \right|^2 \right) \right) = 0 ,
	$$
	off the critical points.
\end{rem}

Back to the general case, when the matrix $A$ is non-constant, instead of moving from our matrix to the Laplacian by means of a change of variables as it is done in the constant case (meaning that we would need a non-constant change of variables), we will work with its symmetric part. Using the particular form $A_0 = R^T BR$, the rotation matrices will play the role of the directional derivatives, so we will be left with the diagonal matrix $B$, and this allows us to apply the previous strategy.

Now we study the first term appearing in the right-hand side in \cref{lemma close to the boundary},
$$
- \int_{\widetilde\Omega} \langle A^T \nabla \log |S\nabla g|^2 , \nabla (\psi_\varepsilon \varphi g) \rangle,
$$
which is the most delicate, and the key point in the different behavior between the planar case and higher dimensions. Recall that, in the dual space, this term is
$$
\left(\divv (A^T\nabla \log |S\nabla g|^2)\right)(\psi_\varepsilon \varphi g) \text{ where } \psi_\varepsilon \varphi g \in W^{1,\infty}_c (\widetilde\Omega_\varphi).
$$

\begin{notation}
	From now on, a boxed term $\boxed{X}$ (representing the functional $T_X$) can be understood as a pointwise function, while a functional written as $T_X$ needs to be understood strictly in the dual sense.
\end{notation}

In the following lemma we decompose this last term.

\begin{lemma}[Step 5]\label{decomposition term 3}%lemma 2
	In the dual space $W_c^{1,\infty}(\widetilde\Omega_\varphi)^\prime$, we have
	$$
	\divv (A^T \nabla \log |S\nabla g|^2) = \boxed{M1} + \boxed{M2} + T_{M3} + T_E,
	$$
	where
	\begin{align*}
		\boxed{M1} & \coloneqq -\sum_{i=1}^2 \frac{4 b_i}{|S\nabla g|^4}  \langle \partial^i \nabla_R g, B\nabla_R g \rangle^2 , \\
		\boxed{M2} & \coloneqq \sum_{i=1}^2 \frac{2 b_i}{|S\nabla g|^2} \langle \partial^i \nabla_R g, B \partial^i \nabla_R g \rangle , \\
		T_{M3} & \coloneqq \frac{2}{|S\nabla g|^2} \left\langle \sum_{i=1}^2 \partial^i  \left(b_i \partial^i \nabla_R g \right), B \nabla_R g \right\rangle ,
	\end{align*}
	and $T_E$ is an error term (involving derivatives of the matrix $A$) of the form
	\begin{equation}\label{T_E simplified}
	T_E = \sum_{i=1}^2 \frac{1}{|S\nabla g|^2} \partial^i \left( \OO(|\nabla g|^2) \right)
	+ \OO\left(1+\frac{|\nabla^2 g|}{|\nabla g|}\right).
	\end{equation}
\end{lemma}

We prove the preceding lemma in \cref{proof of decomposition term 3}.

Note that, using \rf{two directional derivatives}, we have that the two main terms satisfy
\begin{equation}\label{square dependence in two main terms}
	\left|\boxed{M1}\right|+\left|\boxed{M2}\right|\lesssim 1 + \left(\frac{|\nabla^2 g|}{|\nabla g|}\right)^2,
\end{equation}
which is bad for our purposes. A similar behavior occurs with the term $T_{M3}$. Instead, we need to exploit their cancellation, as illustrated in the constant matrix case:

\begin{rem}[Constant matrix and main terms]\label{constant matrix and main terms}
	If the matrix $A$ were constant, and hence $R$ and $B$ were constant as well, then we would have $T_E=0$. That is,
	$$
	\divv \left( A^T \nabla \log |S\nabla g|^2 \right) = \boxed{M1} +  \boxed{M2} +  T_{M3}.
	$$
	This suggests that these 3 terms are the main terms, and the others must be bounded error terms. Moreover, we have the following points:
	\begin{itemize}
		\item The key point in the plane is that $\sum_{i=1}^2 \partial^i (b_i \partial^i g)=0$ would imply $b_1 \partial^{1,1} g=-b_2 \partial^{2,2} g$ (off the critical points). One can show that, as a consequence, the main terms \boxed{M1} and \boxed{M2} would cancel each other in the sense $\boxed{M1} +  \boxed{M2} = 0$ off the critical points of $g$.
		
		\item Since the derivatives would commute, i.e., $\partial^{1,2} = \partial^{2,1}$, we would have
		\begin{equation}\label{extract gradient R outside}
			\sum_{i=1}^2 \partial^i  \left(b_i \partial^i \nabla_R g \right) = \nabla_R \left(\sum_{i=1}^2 \partial^i  \left(b_i \partial^i g \right) \right) = \nabla_R (0) = 0,
		\end{equation}
		where we used that the Green function $g$ satisfies $\sum_{i=1}^2 \partial^i (b_i \partial^i g)=0$. This is why we call $T_{M3}$ the ``zero'' term. In contrast to the previous point, here it is not needed that we are in the plane. However, we use that the Green function is $L_A$-harmonic, that is, this term is zero even in higher dimensions.
	\end{itemize}
	All in all, in the constant matrix case we would have
	$$
	\divv \left( A^T \nabla \log |S\nabla g|^2 \right) =\boxed{M1} +  \boxed{M2} +  T_{M3} = 0
	$$
	off the critical points of $g$.
\end{rem}

The strategy explained in the previous remark is, morally, what we will do to prove the following two lemmas.

\begin{lemma}[Step 6]\label{cancellation main terms M1 M2}%lemma 3
	Pointwise $\left| \boxed{M1} + \boxed{M2} \right| \lesssim 1+|\nabla^2 g| / |\nabla g|$ in $\widetilde\Omega_\varphi$.
\end{lemma}

We prove the preceding lemma in \cref{proof of cancellation main terms M1 M2}.

\begin{lemma}[Step 7]\label{cancellation main term M3}
	In the dual space $W_c^{1,\infty}(\widetilde\Omega_\varphi)^\prime$, the functional $T_{M3}$ is of the form
	$$
	T_{M3} = \sum_{j=1}^2 \frac{2 b_j \partial^j g}{|S\nabla g|^2} \left(
	\OO (|\nabla g|) + \OO (|\nabla^2 g|)
	+\partial^1 \left( \OO (|\nabla g|)\right) 
	+\partial^2 \left( \OO (|\nabla g|) \right)
	\right).
	$$
	Recall that this equality must be understood at a functional level, see \rf{final form after cancellation of M3} below.
\end{lemma}

We prove the preceding lemma in \cref{proof of cancellation main term M3}.

In the following lemma, we bound the error terms that have appeared in the previous computations.

\begin{lemma}[Step 8]\label{bound error and main M3 terms}%lemma 4
	For small enough $\varepsilon >0$,
	$$
	|T_{M3} (\psi_\varepsilon\varphi g)| + |T_E (\psi_\varepsilon\varphi g)|\lesssim 1 + \int_{\supp \varphi} \frac{|\nabla^2 g|}{|\nabla g|} g .
	$$
\end{lemma}

We prove the preceding lemma in \cref{proof of bound error and main M3 terms}.

\begin{lemma}[Step 9]\label{bound fraction derivatives of green function}%lemma 6
	$
	\int_{\supp \varphi} \frac{|\nabla^2 g|}{|\nabla g|} g \leq C < \infty,
	$
	and for small enough $\varepsilon>0$, $\int_{U_\varepsilon (\partial\widetilde\Omega)} \frac{|\nabla^2 g|}{|\nabla g|} g \lesssim \varepsilon$.
\end{lemma}

We prove the preceding lemma at the end of \cref{proof of bound fraction derivatives of green function}.

Using the previous lemmas we obtain \rf{key integral log bound absolute symmetric}.

\begin{proof}[Proof of \rf{key integral log bound absolute symmetric}.]
	Let $p\in \Omega$ with $\dist(p,\partial\Omega)>r_0$ so that \cref{bound easy terms 1 2} holds for all $a>10$ integer. By \cref{integral elliptic measure inicial and truncated,bound easy terms 1 2,truncated solid integral to non truncated} we have
	$$
		\left|
		\int_{\partial\widetilde\Omega} \log |S\nabla g|^2 \, d\widetilde\omega_T^p
		\right|
		\leq \left| - \int \langle A^T \nabla \log |S\nabla g|^2, \nabla (\varphi g)\rangle \right|
		+ \OO\left(1+\int_{\supp \varphi} \frac{|\nabla^2 g|}{|\nabla g|} g \right).
	$$
	By \cref{lemma close to the boundary,decomposition term 3,cancellation main terms M1 M2,bound error and main M3 terms}, we get
	$$
	\left| - \int \langle A^T \nabla \log |S\nabla g|^2, \nabla (\varphi g)\rangle \right| \lesssim 1+\int_{\supp \varphi} \frac{|\nabla^2 g|}{|\nabla g|} g +\OO\left(\lim_{\varepsilon \to 0}\frac{1}{\varepsilon}\int_{U_{3\varepsilon}(\partial\widetilde\Omega)}\frac{|\nabla^2 g|}{|\nabla g|} g\right).
	$$
	In particular,
	$$
	\left|
	\int_{\partial\widetilde\Omega} \log |S\nabla g|^2 \, d\widetilde\omega_T^p
	\right| \lesssim 1 + \int_{\supp \varphi} \frac{|\nabla^2 g|}{|\nabla g|} g +\OO\left(\lim_{\varepsilon \to 0}\frac{1}{\varepsilon}\int_{U_{3\varepsilon}(\partial\widetilde\Omega)}\frac{|\nabla^2 g|}{|\nabla g|} g\right).
	$$
	Finally, \rf{key integral log bound absolute symmetric} follows from this and \cref{bound fraction derivatives of green function}. Note that by Harnack inequality this is true for every $p\in\Omega$ with $\dist(p,\partial\Omega)>r_0$.
\end{proof}

\subsection{Second derivatives of the Green function}

In this subsection, we collect a fundamental property of the Green function in our situation, and some useful equalities and bounds that we will use later on involving the second order derivatives of the Green function.

The following claim points out the main difference between the planar and the higher dimensional case. This is a key point in the proof of \cref{cancellation main terms M1 M2,cancellation main term M3}.

\begin{claim}\label{green function property grad g}
	The Green function $g$ satisfies
	$$
	\sum_{i=1}^2 \partial^i (b_i \partial^i g) 
	= -\sum_{i=1}^{2} b_i \partial^i g  \divv R_i^T 
	- \sum_{i,j=1}^{2} \partial_i \left(\frac{a_{ij}-a_{ji}}{2}\right) \partial_j g
	=\OO (|\nabla g|).
	$$
	a.e.\ in $U_{r_0} (\partial \Omega) \cap \widetilde\Omega$. In particular,
	\begin{equation}
		b_1 \partial^{1,1} g = -b_2 \partial^{2,2} g + \OO(|\nabla g|) \label{property 2 green function}.
	\end{equation}
	\begin{proof}
		By \cref{directional_div,div A and div Asim},
		\begin{align*}
			\divv A\nabla g &=\sum_{i=1}^{2} \partial^i (b_i \partial^i g) + \sum_{i=1}^{2} b_i \partial^i g  \divv R_i^T + \sum_{i,j=1}^{2} \partial_i \left(\frac{a_{ij}-a_{ji}}{2}\right) \partial_j g \\
			&= \sum_{i=1}^{2} \partial^i (b_i \partial^i g) + \OO (|\nabla g|),
		\end{align*}
		where we used that the matrices $A$, $B$ and $R$ are Lipschitz. The claim follows from the fact that the Green function satisfies $\divv A\nabla g = 0$ a.e.\ in $U_{r_0} (\partial\Omega)\cap\widetilde\Omega$, see \cref{pointwise A-solution}.
		
		As a consequence of $\partial^i (b_i \partial^i g) = \partial^i b_i \partial^i g + b_i \partial^{i,i} g$ and $\partial^i b_i \partial^i g \in \OO(|\nabla g|)$, we have \rf{property 2 green function}.
	\end{proof}
\end{claim}

A close look to the proof reveals that if the matrix $A$ were constant, then we would have $b_1 \partial^{1,1} g = -b_2 \partial^{2,2} g$ in the planar case, as it was explained in \cref{constant matrix and main terms}. Instead, we get the analog with some error terms in \rf{property 2 green function}.

In the following claim we collect some basic equalities and estimates. These equalities will be useful in the decomposition procedure in \cref{decomposition term 3}. 

\begin{claim}
	Using the notation $\nabla_R f = R\cdot \nabla f$, see \cref{sec:set up and notation}, we have
	\begin{equation}\label{derivative product of dot product}
		\partial^i |S\nabla g|^2 = 2 \langle \partial^i \nabla_R g, B\nabla_R g \rangle  +  \langle \nabla_R g, \partial^i B \nabla_R g \rangle .
	\end{equation}
	Consequently,
	\begin{equation}\label{derivative quotient}
		\partial^i \left( \frac{1}{|S\nabla g|^2} \right) 
		= \frac{- \partial^i |S\nabla g|^2}{|S\nabla g|^4}
		\overset{\text{\rf{derivative product of dot product}}}{=} -\frac{2 \langle \partial^i \nabla_R g, B\nabla_R g \rangle  +  \langle \nabla_R g, \partial^i B \nabla_R g \rangle}{|S\nabla g|^4},
	\end{equation}
	and
	\begin{equation}\label{gradient log_a gradient g}
		|\nabla \log |S\nabla g|^2|
		= \frac{|\nabla |S\nabla g|^2|}{|S\nabla g|^2} 
		\lesssim 1 + \frac{|\nabla^2 g|}{|\nabla g|} \text{ in }\widetilde\Omega_\varphi .
	\end{equation}
	\begin{proof}
		Since $S^T S = A_0$ and $A_0 = R^T B R$, we can write
		$$
		|S \nabla g|^2 = \langle S\nabla g, S\nabla g\rangle = \langle \nabla g, A_0 \nabla g \rangle = \langle R\nabla g, BR\nabla g \rangle = \langle \nabla_R g, B\nabla_R g \rangle .
		$$
		Recall that $R$ and $B$ are Lipschitz regular, thus we can apply the usual derivative rules a.e.\ to obtain
		$$
		\partial^i |S\nabla g|^2 = \langle \partial^i \nabla_R g, B\nabla_R g \rangle + \langle \nabla_R g, \partial^i B \nabla_R g \rangle + \langle \nabla_R g, B \partial^i \nabla_R g \rangle,
		$$
		and \rf{derivative product of dot product} follows from the symmetry of $B$.
		
		The inequality in \rf{gradient log_a gradient g} follows from
		\begin{align*}
			\left| \nabla |S\nabla g|^2 \right| 
			\overset{\text{\rf{derivative product of dot product}}}&{\approx}
			\sum_{i=1}^2 \left| 2\langle \partial^i \nabla_R g,  B \nabla_R g \rangle + \langle \nabla_R g, \partial^i B \nabla_R g \rangle \right| \\
			& \leq \sum_{i=1}^2 2|\langle \partial^i \nabla_R g,  B \nabla_R g \rangle|
			+\sum_{i=1}^2 |\langle \nabla_R g, \partial^i B \nabla_R g \rangle | \\
			& \lesssim |\nabla g|^2 + \sum_{i=1}^2 | \partial^i \nabla_R g| | \nabla g|
			\overset{\text{\rf{two directional derivatives}}}{\lesssim} |\nabla g|^2 + |\nabla g| | \nabla^2 g| ,
		\end{align*}
	and this concludes the proof of the claim.
	\end{proof}
\end{claim}

\begin{rem}\label{quotient gradient green is in sobolev 12}
	$1/|S\nabla g|^2 \in W^{1,2}_{\loc} (\widetilde\Omega_\varphi)$ because of \rf{derivative quotient}, $|\nabla g|\approx g/\dist(\cdot,\partial\widetilde\Omega) >0$ in $\widetilde\Omega_\varphi$ and $g\in W^{2,2}(\widetilde\Omega_\varphi)$. By the same reason, $\log |S\nabla g|^2\in W^{1,2}_{\loc} (\widetilde \Omega_\varphi)$.
\end{rem}

\subsection{Whitney cubes and proof of Step 9}\label{proof of bound fraction derivatives of green function}
A classical way to compute an integral over a given set is to discretize it. We will do so using Whitney cubes, that is, we divide the domain into regions (cubes) which have diameter comparable to their distance to the boundary, so that Harnack's and Caccioppoli's inequalities can be used locally.

Let us start by defining the Whitney covering of $\widetilde \Omega$, and then we will move to study the properties of those cubes touching $\supp \varphi$.

\begin{definition}
	There exists a collection $W_{\widetilde \Omega}$ of dyadic cubes satisfying the following properties:
	\begin{enumerate}[label=({W\arabic*}),ref={W\arabic*}]
		\item\label{whitney p1} $\sqrt{2} \ell(Q) \leq \dist (Q,\partial \widetilde \Omega) \leq 4\sqrt{2} \ell (Q)$, equivalently, $\diam Q \leq \dist (Q,\partial \widetilde \Omega) \leq 4 \diam Q$.
		\item\label{whitney p2} $1.5Q \cap \partial \widetilde \Omega = \emptyset$, since $\dist (1.5Q, \partial \widetilde \Omega) \geq \dist (Q, \partial \widetilde \Omega)-  \sqrt{2} \cdot 0.25 \ell (Q) \geq  0.75 \sqrt{2} \ell (Q) > 0$,
	\end{enumerate}
	which we call \emph{Whitney cubes}, see \cite[Appendix J]{Grafakos2009_classical}. We define $W_\varphi \coloneqq \{Q\in W_{\widetilde \Omega} : Q\cap \supp \varphi \not = \emptyset \}$.
\end{definition}

\begin{lemma}\label{whitney cubes properties}
	Every $Q\in W_\varphi$ satisfies also:
	\begin{enumerate}[label=({W\arabic*}),ref={W\arabic*}]
		\setcounter{enumi}{2}
		\item\label{whitney cubes properties points are near} $\dist (x,\partial \widetilde \Omega) \leq 7 \frac{R}{2000c} \leq R/ c$ for $x\in 1.5Q$.
		\item\label{whitney cubes properties lenght bounded above} $\ell (Q) \lesssim 1$.
		\item If $Q \cap  \supp \nabla \varphi \not = \emptyset$, then $\ell (Q) \gtrsim 1$.
		\item\label{green bounded by measure and 1 in the supp of varphi} For $x\in 1.5Q$ we have 
		$$
		g(x) \approx g(A(\xi, 10\ell (Q))) \lesssim \widetilde\omega_T (B(\xi, 10\ell (Q))) \leq \widetilde\omega_T (50Q) \leq 1 ,
		$$
		where $\xi\in \partial \widetilde \Omega$ is such that $\dist(\xi, Q) = \dist(\partial \widetilde \Omega, Q)$, and $A(\xi, 10\ell (Q))$ is the Corkscrew point at $\xi$ with radius $10\ell (Q)$.
	\end{enumerate}
	
	\begin{proof}
		Note that for $x\in 1.5Q$, Harnack's inequality gives $g(x) \approx g(A(\xi, 10\ell (Q)))$. Using the relation between the elliptic measure and Green's function in NTA domains (see \cite[Lemma 1.3.3]{Kenig1994}), we get
		$$
		g(A(\xi, 10\ell (Q))) \lesssim \widetilde\omega_T (B(\xi, 10\ell (Q))).
		$$
		Here we need $p\not\in B(\xi,3\cdot 10\ell (Q))$, which is granted as long as $\partial \widetilde \Omega$ and $\partial \Omega$ are close enough.
		
		The remaining estimates in the lemma follow from the definition by standard arguments. The details are left to the reader.
	\end{proof}
\end{lemma}

In order to bound the error terms arising in the proofs of \cref{sec:sketch of the proof}, we will use (without mention) the following remark.

\begin{rem}
	If $f \in \OO (|\nabla g|^k)$ then $\left| \int \frac{f}{|\nabla g|^k} \varphi g \right| \leq C < \infty$.
	\begin{proof}
		Since $g\lesssim 1$ in $\supp \varphi$ (see \cref{whitney cubes properties}\rf{green bounded by measure and 1 in the supp of varphi}), we obtain
		\begin{equation*}
			\left| \int \frac{f}{|\nabla g|^k} \varphi g \right| \leq \int \frac{|f|}{|\nabla g|^k} |\varphi g|  \lesssim \int \left| \varphi g \right| \lesssim \HH^2(\Omega) \leq (\diam \partial\Omega)^2 <\infty,
		\end{equation*}
		where we used that $\Omega$ is bounded.
	\end{proof}
\end{rem}

Next we continue by controlling the Green function over Whitney cubes:
\begin{lemma}\label{sum over whitney region}\label{integrals involving green function near boundary}
	We have:
	\begin{enumerate}
		\item\label{sum over whitney region p1} $\sum_{Q \in W_\varphi} \ell(Q) \cdot \widetilde \omega_T (50Q) \lesssim 1$.
		\item\label{sum over whitney region p1 epsilon} For $\varepsilon >0$, 
		$
		\sum_{\left\{Q\in W_{\widetilde\Omega} : Q\cap U_\varepsilon (\partial\widetilde\Omega)\not=\emptyset\right\}} \ell (Q) \cdot \widetilde\omega_T (50Q) \lesssim \varepsilon
		$.
		\item\label{integral g near boundary}For $\varepsilon>0$, $\int_{U_\varepsilon (\partial\widetilde\Omega)} g \lesssim \varepsilon^2$.
		\item\label{sum over whitney region p2} For each $Q\in W_\varphi$, $\left(\int_{1.5Q} g^2\right)^{1/2} \lesssim \ell(Q) \cdot \widetilde \omega_T (50Q)$.
		\item\label{sum over whitney region p3} For each $Q\in W_\varphi$, $\int_Q |\nabla g| \lesssim \ell(Q) \cdot \widetilde \omega_T (50Q)$.
		\item\label{sum over whitney region p4} For each $Q\in W_\varphi$, $\int_Q |\nabla^2 g| \lesssim \widetilde \omega_T (50Q)$.
	\end{enumerate}
\begin{proof}
	\Cref{sum over whitney region p1}: As $\ell (Q) \lesssim 1$ for every $Q\in W_\varphi$ (see \cref{whitney cubes properties}\rf{whitney cubes properties lenght bounded above}), there is $k_0$ such that every $Q\in W_\varphi$ satisfies $\ell(Q) = 2^{-k}$ with $k\geq k_0$. Note also that, for each scale $k\geq k_0$, the family $\{50Q\}_{\{Q\in W_\varphi \, : \, \ell(Q) = 2^{-k}\}}$ has finite overlapping. Therefore,
\begin{align*}
	\sum_{Q \in W_\varphi} \ell(Q) \cdot \widetilde \omega_T (50Q) & =  \sum_{k\geq k_0} \sum_{\substack{Q \in W_\varphi \\ \ell(Q)=2^{-k}}} \ell(Q) \cdot \widetilde\omega_T (50{Q})
	= \sum_{k\geq k_0} 2^{-k} \sum_{\substack{Q \in W_\varphi \\ \ell(Q)=2^{-k}}}  \widetilde\omega_T (50{Q}) \\
	& \lesssim \sum_{k\geq k_0} 2^{-k} \widetilde\omega_T 
	\Biggl( \bigcup_{\substack{Q \in W_\varphi \\ \ell(Q)=2^{-k}}}  50{Q} \Biggr) 
	\leq  \sum_{k\geq k_0} 2^{-k} < \infty.
\end{align*}

	\Cref{sum over whitney region p1 epsilon}: Recall that the cubes $Q\in W_{\widetilde\Omega}$ with $Q\cap U_{\varepsilon} (\partial\widetilde\Omega)\not=\emptyset$ satisfy $\ell(Q) \lesssim \varepsilon$, and $\ell (Q) \approx \varepsilon$ if $Q\cap (\partial U_{\varepsilon} (\partial\widetilde\Omega) \setminus \partial\widetilde\Omega)\not=\emptyset$, see \rf{whitney p1}. Hence, as in the proof of \cref{sum over whitney region p1}, taking $k_0 (\varepsilon) \approx -\log_2 \varepsilon$ we have
	$$
	\sum_{\substack{Q\in W_{\widetilde\Omega} \\ Q\cap U_{\varepsilon} (\partial\widetilde\Omega)\not=\emptyset}} \ell (Q) \widetilde\omega_T (50Q)\\
	\leq \sum_{k\geq k_0 (\varepsilon)} 2^{-k} \sum_{\substack{Q \in W_{\widetilde\Omega} \\ \ell(Q)=2^{-k}}} \widetilde\omega_T(50Q)
	\lesssim \sum_{k\geq k_0 (\varepsilon)} 2^{-k}
	\approx 2^{\log_2 \varepsilon} = \varepsilon.
	$$
	
	\Cref{integral g near boundary}: Using that $g\lesssim \widetilde\omega_T (50Q)$ in $Q$, see \rf{green bounded by measure and 1 in the supp of varphi}, that the cubes $Q\in W_{\widetilde\Omega}$ with $Q\cap U_\varepsilon (\partial\widetilde\Omega)\not=\emptyset$ satisfy $\ell (Q) \lesssim \varepsilon$, see \rf{whitney p1}, and \cref{sum over whitney region p1 epsilon}, we have
	\begin{equation*}
		\int_{U_{\varepsilon} (\partial\widetilde\Omega)} g
		\leq \sum_{\substack{Q\in W_{\widetilde\Omega} \\ Q\cap U_{\varepsilon} (\partial\widetilde\Omega)\not=\emptyset}} \int_Q g
		\overset{\text{\rf{green bounded by measure and 1 in the supp of varphi},\rf{whitney p1}}}{\lesssim} \varepsilon \sum_{\substack{Q\in W_{\widetilde\Omega} \\ Q\cap U_{\varepsilon} (\partial\widetilde\Omega)\not=\emptyset}} \ell (Q) \widetilde\omega_T (50Q) \overset{\text{\rf{sum over whitney region p1 epsilon}}}{\lesssim} \varepsilon^2 .
	\end{equation*}

	\Cref{sum over whitney region p2} follows since $g(x) \lesssim \widetilde \omega_T (50Q)$ for $x\in 1.5Q$ (see \cref{whitney cubes properties}\rf{green bounded by measure and 1 in the supp of varphi}). \Cref{sum over whitney region p3} follows from Cauchy-Schwarz and Caccioppoli's inequalities and \cref{sum over whitney region p2}.
	
	\Cref{sum over whitney region p4}: Since the Green function $g$ of $\widetilde \Omega$ is solution of $\divv (A \nabla \cdot) = 0$ in $1.1Q$, the function $\widehat{g} (\cdot) = g(\ell(Q) \cdot)$ is solution of $\divv (\widehat{A} \nabla \cdot)$ in $1.1\widehat Q$, where $\widehat{Q}\coloneqq \{x/\ell(Q) : x\in Q\}$ has side-length one, and $\widehat{A} (\cdot) = A(\ell (Q) \cdot)$. Moreover, since $\partial_i \partial_j \widehat g (\cdot) = \ell (Q)^2 \partial_i \partial_j g (\ell (Q) \cdot)$, we get
	$$
	\int_Q |\nabla^2 g (x)|^2 \, dx = \ell(Q)^2 \int_{\widehat Q} |\nabla^2 g (\ell (Q) x)|^2 \, dx = \frac{1}{\ell(Q)^2 }\int_{\widehat Q} |\nabla^2 \widehat g (x)|^2 \, dx .
	$$

	The matrix $\widehat A$ has the same ellipticity constant as $A$, but as $\ell (Q) \lesssim 1$ for every $Q\in W_\varphi$, see \cref{whitney cubes properties}\rf{whitney cubes properties lenght bounded above}, the Lipschitz norm cannot grow too much:
	$$
	\| \widehat A \|_{C^{0,1} \left(1.1\widehat Q\right)} \coloneqq \|\widehat A\|_{L^\infty \left(1.1\widehat Q\right)} + [\widehat A]_{C^{0,1} \left(1.1\widehat Q\right)} 
	=  \| A\|_{L^\infty ( 1.1 Q )} + \ell(Q) [A]_{C^{0,1} ( 1.1 Q)} \leq C \| A \|_{C^{0,1} \left(1.1 Q\right)} .
	$$
	Moreover, $\dist\left(\widehat{Q}, (1.1\widehat{Q})^c \right) \approx 1$. Therefore we can apply \rf{cacciopoli type ineq for weak twice diff solutions} in \cref{twice weak dif solutions} to obtain
	\begin{equation*}
		\int_Q |\nabla^2 g|^2 \, dx = \frac{1}{\ell (Q)^2} \int_{\widehat Q} |\nabla^2 \widehat g|^2 \, dx \overset{\text{\rf{cacciopoli type ineq for weak twice diff solutions}}}{\lesssim} \frac{1}{\ell (Q)^2} \left[ \left( \int_{1.1 \widehat Q} |\nabla \widehat g|^2 \right)^{1/2} + \left( \int_{1.1 \widehat Q} \widehat g ^2 \right)^{1/2} \right]^2 .
	\end{equation*}
	Applying Caccioppolli's inequality and \cref{sum over whitney region p2}, this implies
	\begin{equation*}
		\int_Q |\nabla^2 g|^2 \, dx \lesssim 
		\frac{1}{\ell (Q)^2}  \int_{1.5 \widehat Q} \widehat g ^2 
		= \frac{1}{\ell (Q)^4}  \int_{1.5  Q}  g^2 \overset{\text{\rf{sum over whitney region p2}}}{\lesssim} \frac{\widetilde \omega_T (50Q)^2}{\ell(Q)^2}.
	\end{equation*}
	Using the Cauchy-Schwarz inequality and this, we get
	\begin{equation*}
		\int_Q |\nabla^2 g| \leq \ell(Q) \left( \int_Q |\nabla^2 g|^2 \right)^{1/2} \lesssim \widetilde\omega_T (50{Q}) .
	\end{equation*}
\end{proof}
\end{lemma}

Using this lemma we obtain \cref{bound fraction derivatives of green function}.

\begin{proof}[Proof of \cref{bound fraction derivatives of green function}]%lemma 6
Recall that, by \rf{key comparability inside the support} we have
$$
|\nabla g(x)| \gtrsim \frac{g(x)}{\dist \left(x, \partial \widetilde \Omega \right)} \text{ for }x\in \supp \varphi,
$$
provided the domain $\widetilde \Omega$ is $(\delta, r_0 /2)$-Reifenberg flat with $\delta$ small enough depending only on $\lambda$ and $\kappa \|A\|_{L^\infty (\R^2)}$. With this it suffices to show
$$
\int_{\supp \varphi} |\nabla^2 g (x)| \dist\left( x, \partial \widetilde \Omega \right) \, dx \leq C .
$$

Summing in $W_\varphi$, by \rf{sum over whitney region p4} and \rf{sum over whitney region p1} in \cref{sum over whitney region} we obtain
\begin{equation*}
	\int_{\supp \varphi} |\nabla^2 g (x)| \dist(x, \partial \Omega) \, dx  \lesssim \sum_{Q \in W_\varphi} \ell (Q) \int_Q |\nabla^2 g|
	\overset{\text{\rf{sum over whitney region p4}}}{\lesssim} \sum_{Q \in W_\varphi} \ell(Q) \cdot \widetilde\omega_T (50{Q}) \overset{\text{\rf{sum over whitney region p1}}}{\lesssim} 1,
\end{equation*}
as claimed.

Let us see the localized version. As above, using \rf{key comparability inside the support}, and \rf{sum over whitney region p4} and \rf{sum over whitney region p1 epsilon} in \cref{sum over whitney region},
\begin{align*}
	\int_{U_{\varepsilon}(\partial\widetilde\Omega)} \frac{|\nabla^2 g|}{|\nabla g|} g
	\overset{\text{\rf{key comparability inside the support}}}&{\lesssim} \int_{U_{\varepsilon} (\partial\widetilde\Omega)} |\nabla^2 g| \dist (x,\partial\widetilde\Omega)
	\leq \sum_{\substack{Q\in W_{\widetilde\Omega} \\ Q\cap U_{\varepsilon} (\partial\widetilde\Omega) \not=\emptyset}} \ell (Q) \int_Q |\nabla^2 g| \\
	\overset{\text{\rf{sum over whitney region p4}}}&{\lesssim} \sum_{\substack{Q\in W_{\widetilde\Omega} \\ Q\cap U_{\varepsilon} (\partial\widetilde\Omega) \not=\emptyset}} \ell (Q) \widetilde\omega_T (50Q) \overset{\text{\rf{sum over whitney region p1 epsilon}}}{\lesssim} \varepsilon .
\end{align*}
\end{proof}

\subsection{Proof of Steps 1 to 8}

In this subsection we prove remaining Steps 1 to 8 in \cref{sec:sketch of the proof}.

\subsubsection{Proof of \cref{integral elliptic measure inicial and truncated,truncated solid integral to non truncated,lemma close to the boundary}}\label{section:link lemmas}

\begin{proof}[Proof of \cref{integral elliptic measure inicial and truncated}]
	We have
	\begin{multline*}
		\left|\int_{\partial\widetilde\Omega} \log |S(\xi)\nabla g(\xi)|^2 \, d\widetilde\omega_T^{p} (\xi)
		-\int_{\partial\widetilde\Omega} \log_{(a)} |S(\xi)\nabla g(\xi)|^2 \, d\widetilde\omega_T^{p} (\xi)\right|\\
		= \left|\int_{\{\xi\in\partial\widetilde\Omega : |S(\xi)\nabla g(\xi)|^2\leq e^{-a} \}} (\log |S(\xi)\nabla g(\xi)|^2 + a) \, d\widetilde\omega_T^{p} (\xi) \right|
	\end{multline*}
	Since the domain is smooth, we can use \rf{elliptic measure density} to write and bound the right-hand side term as
        \begin{align*}
		\Biggl| \int_{\{\xi\in\partial\widetilde\Omega : |S(\xi)\nabla g(\xi)|^2\leq e^{-a} \}} &-\langle A(\xi) \nabla g(\xi), \nu\rangle(\log |S(\xi)\nabla g(\xi)|^2 + a) \, d\sigma (\xi) \Biggr| \\
		&\lesssim \int_{\{\xi\in\partial\widetilde\Omega : |S(\xi)\nabla g(\xi)|^2\leq e^{-a} \}}  |S(\xi)\nabla g(\xi)|(\left|\log |S(\xi)\nabla g(\xi)|\right|) \, d\sigma (\xi) \\
		&\leq \sigma(\partial\widetilde\Omega) \left(e^{-a/2}|{\log e^{-a/2}}|\right) .
	\end{align*}
\end{proof}

\begin{proof}[Proof of \cref{truncated solid integral to non truncated}]
	We have
	\begin{multline*}
		- \int \langle A^T \nabla \log_{(a)} |S\nabla g|^2, \nabla (\varphi g)\rangle
		+ \int \langle A^T \nabla \log |S\nabla g|^2, \nabla (\varphi g)\rangle\\
		=  \int \langle A^T \nabla \min\{0, a+\log |S\nabla g|^2\}, \nabla (\varphi g)\rangle .
	\end{multline*}
	On the other hand, by \rf{key comparability inside the support} we have that if $a\geq 10$ is big enough then $\varphi = 1$ in $\widetilde\Omega_\varphi \cap \{|S\nabla g|^2 \leq e^{-a}\}$, and hence
	$$
	\int \langle A^T \nabla \min\{0, a+\log |S\nabla g|^2\}, \nabla (\varphi g)\rangle = \int \langle A^T \nabla \log |S\nabla g|^2, \nabla g \rangle \characteristic_{\{|S\nabla g|^2 \leq e^{-a}\}} .
	$$
	By \rf{gradient log_a gradient g}, this is controlled in absolute value by
	$$
	\int_{\{|S\nabla g|^2 \leq e^{-a}\}} \left(1+\frac{|\nabla^2 g|}{|\nabla g|}\right) |\nabla g| 
	\leq \int_{\{|S\nabla g|^2 \leq e^{-a}\} } |\nabla g| 
	+ \int_{\{|S\nabla g|^2 \leq e^{-a}\}} |\nabla^2 g| .
	$$
	The lemma follows because the first term in the right-hand side is bounded by a constant times $\HH^2 (\Omega) e^{-a/2}$.
\end{proof}

\begin{proof}[Proof of \cref{lemma close to the boundary}]
	For $\varepsilon >0$ small we have $(1-\psi_\varepsilon) \varphi = (1-\psi_\varepsilon)$. Hence
	\begin{multline*}
		- \int \langle A^T \nabla \log |S\nabla g|^2, \nabla ((1-\psi_\varepsilon) \varphi g)\rangle
		= - \int \langle A^T \nabla \log |S\nabla g|^2, \nabla ((1-\psi_\varepsilon) g)\rangle \\
		= \int \langle A^T \nabla \log |S\nabla g|^2, \nabla \psi_\varepsilon \rangle g
		- \int \langle A^T \nabla \log |S\nabla g|^2, \nabla g \rangle (1-\psi_\varepsilon).
	\end{multline*}
	Using \rf{gradient log_a gradient g}, the second term in the right-hand side is controlled in absolute value by
	$$
	\int_{U_{3\varepsilon} (\partial\widetilde\Omega)} \left(1+\frac{|\nabla^2 g|}{|\nabla g|}\right) |\nabla g|
	= \int_{U_{3\varepsilon} (\partial\widetilde\Omega)} |\nabla g|
	+\int_{U_{3\varepsilon}(\partial\widetilde\Omega)} |\nabla^2 g| \overset{\varepsilon \to 0}{\longrightarrow} 0,
	$$
	and the first term is controlled in absolute value by
	$$
	\int_{U_{3\varepsilon} (\partial\widetilde\Omega)} \left(1+\frac{|\nabla^2 g|}{|\nabla g|}\right) g |\nabla \psi_\varepsilon|
	= \frac{1}{\varepsilon}\int_{U_{3\varepsilon} (\partial\widetilde\Omega)} g
	+\frac{1}{\varepsilon} \int_{U_{3\varepsilon}(\partial\widetilde\Omega)} \frac{|\nabla^2 g|}{|\nabla g|} g.
	$$
	The lemma follows by applying \rf{integral g near boundary} in \cref{integrals involving green function near boundary} in the last line.
\end{proof}

\subsubsection{Proof of \cref{bound easy terms 1 2}}\label{proof of bound easy terms 1 2}%lemma 1

Recall that the Green function $g$ in $\widetilde \Omega$ for the operator $L_A u = \divv A \nabla u$ satisfies, for every $\phi \in C^0 (\overline{\widetilde\Omega}) \cap W^{1,2} (\widetilde\Omega)$,
\begin{equation}\label{recall boundary integral to interior integral}
\int_{\partial \widetilde \Omega} \phi (\xi) \, d \widetilde \omega_T^{x} (\xi)-\phi(x)=-\int_{\widetilde \Omega} \langle A^T(y) \nabla \phi(y),\nabla g_x (y) \rangle \, d y, \text { for a.e.\ } x \in \widetilde \Omega ,
\end{equation}
see \rf{boundary integral to interior integral}.

We claim that
$$
\varphi \log_{(a)} |S\nabla g|^2 \in C^0 (\overline{\widetilde\Omega}) \cap W^{1,2} (\widetilde\Omega).
$$
Indeed, $\varphi \log_{(a)} |S\nabla g|^2 \in C^0 (\overline{\widetilde\Omega})$ because the function $\varphi$ avoids the pole and $\nabla g$ is continuous up to the boundary, see \cref{C^a to the boundary}. Let us see that $\varphi \log_{(a)} |S\nabla g|^2 \in W^{1,2} (\widetilde\Omega)$. Using Jensen's inequality as the function $x\mapsto \left|\log x\right|^2$ is concave for $x>e$, the $L^2$-norm is controlled (depending on $a$) by
\begin{equation}\label{L^2 norm of log_a}
	\begin{aligned}
		\int_{\widetilde\Omega\cap\supp \varphi} \left|\log_{(a)} |S\nabla g|^2\right|^2
		&\leq \int_{\widetilde\Omega\cap \{|S\nabla g|^2\leq e\}} \left|\log_{(a)} |S\nabla g|^2\right|^2
		+ \int_{\widetilde\Omega\cap\supp\varphi\cap\{|S\nabla g|^2>e\}} \left|\log |S\nabla g|^2\right|^2 \\
		& \lesssim \HH^2 (\widetilde\Omega) \left( a^2 + \left|\log\left( \avint_{\widetilde\Omega\cap\supp\varphi\cap\{|S\nabla g|^2>e\}} |S\nabla g|^2\right) \right|^2\right) <\infty,
	\end{aligned}
\end{equation}
which is bounded since $g\in W^{1,2}(\widetilde\Omega\cap\supp\varphi)$. The $L^2$ norm of $\nabla(\varphi \log_{(a)} |S\nabla g|^2)$ is
\begin{equation*}
	\int_{\widetilde\Omega} |\nabla (\varphi \log_{(a)} |S\nabla g|^2)|^2
	\leq \int_{\widetilde\Omega} |\nabla \varphi|^2 \left|\log_{(a)} |S\nabla g|^2\right|^2
	+\int_{\widetilde\Omega} \varphi^2 |\nabla \log_{(a)} |S\nabla g|^2|^2.
\end{equation*}
The first term is finite by \rf{L^2 norm of log_a} because $\varphi \in C^\infty_c (\R^2)$. From the definition of $\log_{(a)} x \coloneqq \max \{\log x, -a\}$ and by \rf{gradient log_a gradient g}, we get
\begin{equation*}
	|\nabla \log_{(a)} |S\nabla g|^2|
	\lesssim \characteristic_{\{x\geq e^{-a}\}}(|S\nabla g|^2)\left(1 + \frac{|\nabla^2 g|}{|\nabla g|}\right),
\end{equation*}
and hence the second term is controlled by
\begin{equation*}
	\int_{\widetilde\Omega\cap\supp\varphi} |\nabla \log_{(a)} |S\nabla g|^2|^2
	\lesssim
	\int_{\widetilde\Omega\cap\supp\varphi\cap\{|S\nabla g|^2 \geq e^{-a}\}} 1 + \frac{|\nabla^2 g|^2}{|\nabla g|^2}
	\lesssim \HH^2 (\widetilde\Omega) 
	+ e^a \int_{\widetilde\Omega\cap\supp\varphi} |\nabla^2 g|^2 <\infty,
\end{equation*}
which is bounded as $g\in W^{2,2} (\widetilde\Omega\cap\supp\varphi)$.

By the choice of $\varphi$ with $\varphi (p)=0$ and $\varphi |_{\partial\widetilde\Omega}=1$, and plugging $\phi=\varphi \log_{(a)} |S\nabla g|^2 \in C^0 (\overline{\widetilde\Omega}) \cap W^{1,2} (\widetilde\Omega)$ in \rf{recall boundary integral to interior integral}, for a.e.\ $p\in \Omega$ with $\dist(p,\partial\Omega)>r_0$ we have
\begin{equation}\label{integral with 1 2 3}
	\begin{aligned}
		\int_{\partial \widetilde \Omega}  \log_{(a)} |S\nabla g |^2 \, d\widetilde \omega_T^p &= - \int_{\widetilde \Omega \cap \supp \varphi}  \langle A^T \nabla \left(\varphi \log_{(a)} |S\nabla g|^2 \right), \nabla g \rangle \\
		&= - \int  \langle A^T \nabla \varphi, \nabla g \rangle \cdot \log_{(a)} |S\nabla g|^2 - \int  \langle  A^T \nabla \log_{(a)} |S\nabla g|^2 , \varphi \nabla g \rangle \\
		&= \boxed{1} + \boxed{2}
		+ \boxed{3} ,
	\end{aligned}
\end{equation}
with
\begin{align*}
	\boxed{1} &\coloneqq - \int  \langle A^T \nabla \varphi, \nabla g \rangle \cdot \log_{(a)} |S\nabla g|^2 ,\\
	\boxed{2} &\coloneqq \int \langle A^T \nabla \log_{(a)} |S\nabla g|^2,\nabla \varphi \rangle \cdot g, \text{ and}\\
	\boxed{3} &\coloneqq - \int \langle  A^T \nabla \log_{(a)} |S\nabla g|^2 , \nabla (\varphi g) \rangle .
\end{align*}

We claim that $\left|\boxed{1}\right| \lesssim 1$ and $\left|\boxed{2}\right| \lesssim 1+\int_{\supp \varphi} \frac{|\nabla^2 g|}{|\nabla g|} g$, and hence \cref{bound easy terms 1 2} follows.

We start by controlling the term $\boxed{1}$. Recall that we have fixed $\varphi$ so that $|\nabla \varphi| \lesssim 1$, and assumed that $a$ is big enough so that $e^{-a}<\min_{z\in\supp\nabla\varphi} |S\nabla g(z)|^2$, that is, $\log_{(a)} |S\nabla g(z)|^2=\log |S\nabla g(z)|^2$ for $z\in \supp\nabla\varphi$. Thus,
\begin{equation*}\label{first term step 1}
	\left|\boxed{1}\right| \lesssim \int_{\supp \nabla \varphi} |\nabla g| \left|\log |S\nabla g|\right|
	\approx  \int_{\supp \nabla \varphi} |S\nabla g| \left|\log |S\nabla g|\right| .
\end{equation*}
Since $\supp \nabla \varphi$ is far from $\partial \widetilde \Omega$, in particular $\dist(x, \partial \widetilde \Omega) \approx 1$ for $x\in \supp \nabla \varphi$. Consequently, from this and \rf{key comparability inside the support} we have $|\nabla g| \approx g/\dist(\cdot, \partial \widetilde \Omega) \approx g$ in the support of $\nabla \varphi$, and hence
$$
|S\nabla g(x)| \approx |\nabla g(x)| \approx g(x) \overset{\text{\rf{green bounded by measure and 1 in the supp of varphi}}}{\lesssim} 1 \text{ for } x \in \supp \nabla \varphi .
$$
Note that for $0<t<t_0$ we have $t\left|\log t\right| \leq \max\{1, t_0 \left|\log t_0\right|\}$. Therefore,
\begin{equation*}\label{first term step 3.2}
	\boxed{1} \lesssim \int_{\supp \nabla \varphi} |S\nabla g| \left|\log |S \nabla g|\right| \lesssim \HH^2 (\supp \nabla \varphi) \leq C <\infty .
\end{equation*}

Now we study the term $\boxed{2}$. Since $\varphi$ is fixed, we have $|\nabla \varphi| \lesssim 1$ and so 
\begin{equation*}
	\left| \boxed{ 2 } \right| \lesssim \int_{\supp \nabla \varphi}  g \cdot | \nabla \log_{(a)} |S\nabla g|^2 | 
	\leq \int_{\supp \varphi}  g \cdot | \nabla \log_{(a)} |S\nabla g|^2 |.
\end{equation*}
Plugging \rf{gradient log_a gradient g} in here we get
$$
\left| \boxed{2 } \right| \lesssim \int_{\supp \varphi} g +  \int_{\supp \varphi} \frac{|\nabla^2 g|}{|\nabla g|} g \overset{\text{\rf{green bounded by measure and 1 in the supp of varphi}}}{\lesssim} 1 +  \int_{\supp  \varphi} \frac{|\nabla^2 g|}{|\nabla g|} g .
$$
\qed

\subsubsection{Proof of \cref{decomposition term 3}}\label{proof of decomposition term 3}%lemma 2

In this subsection we study, via a perturbation argument, the functional term
$$
\divv \left(A^T\nabla \log |S\nabla g|^2\right) \in W^{1,\infty}_c (\widetilde\Omega_\varphi)^\prime.
$$
Note that the action of this functional on $\varphi g$ gives rise to a modified version of the term $\boxed{3}$ in \rf{integral with 1 2 3}.

First, we move from the matrix $A$ to its symmetric part $A_0$ by \cref{div A and div Asim}, and we write its divergence in terms of the directional derivatives using \cref{directional_div} (see \cref{quotient gradient green is in sobolev 12}):
\begin{align*}
	\divv  \left(  A^T  \nabla \log |S\nabla g|^2 \right)  = &
	\divv \left( A_0 \nabla \log |S\nabla g|^2\right) + \sum_{i,j=1}^2 \partial_i \left(\frac{a_{ij}^T-a_{ji}^T}{2}\right) \partial_j \log |S\nabla g|^2 \\
	=& \sum_{i=1}^2 \partial^i (b_i \partial^i (\log |S\nabla g|^2)) 
	+ \sum_{i=1}^2 b_i \partial^i (\log |S\nabla g|^2) \divv R^T_i \\
	&+ \sum_{i,j=1}^{2} \partial_i \left(\frac{a_{ji}-a_{ij}}{2}\right) \partial_j (\log |S\nabla g|^2) 
	\eqqcolon  T_{3.A} + \boxed{3.B} + \boxed{3.C}.
\end{align*}

Note that the terms \boxed{3.B} and \boxed{3.C} contain derivatives of the matrix $A$. In particular, if the matrix were constant then these two terms would be zero, which suggests that the terms \boxed{3.B} and \boxed{3.C} must be bounded error terms.

Next we deal with $T_{3.A}$. As $1/|S\nabla g|^2 \in W^{1,2}_{\loc} (\widetilde\Omega_\varphi)$ (see \cref{quotient gradient green is in sobolev 12}), by \rf{product rule distributional directional derivative} we have
\begin{equation}\label{term 3A}
	\begin{aligned}
		T_{3.A} =& \sum_{i=1}^2 \partial^i \left(b_i \frac{\partial^i |S\nabla  g|^2}{|S\nabla g|^2} \right) \\
		=& \sum_{i=1}^2 \partial^i \left( \frac{1}{|S\nabla g|^2} \right) b_i \partial^i |S\nabla g|^2
		+ \sum_{i=1}^2 \frac{1}{|S\nabla g|^2} \partial^i (b_i \partial^i |S\nabla g|^2).
	\end{aligned}
\end{equation}

We infer that each element in the sum in the first term in the right-hand side in \rf{term 3A} is
\begin{equation}\label{using derivative product in the first term}
	\partial^i \left( \frac{1}{|S\nabla g|^2} \right) b_i \partial^i |S\nabla g|^2 \overset{\text{\rf{derivative quotient}}}{=} \frac{- b_i \left( \partial^i |S\nabla g|^2 \right)^2}{|S\nabla g|^4} .
\end{equation}
Therefore, using identities \rf{derivative product of dot product} and \rf{using derivative product in the first term} and expanding, the first term in the right-hand side of \rf{term 3A} can be written as
$$
\sum_{i=1}^2 \partial^i \left(\frac{1}{|S\nabla g|^2} \right) b_i \partial^i |S\nabla g|^2 = \boxed{M1} + \boxed{E1} + \boxed{E2},
$$
where
\begin{align*}
	\boxed{M1} & \coloneqq -\sum_{i=1}^2 \frac{4 b_i}{|S\nabla g|^4}  \langle \partial^i \nabla_R g, B\nabla_R g \rangle^2 , \\
	\boxed{E1} & \coloneqq -\sum_{i=1}^2 \frac{4 b_i}{|S\nabla g|^4}  \langle \partial^i \nabla_R g, B\nabla_R g \rangle  \langle \nabla_R g, \partial^i B \nabla_R g \rangle , \text{ and}\\
	\boxed{E2} & \coloneqq -\sum_{i=1}^2 \frac{b_i}{|S\nabla g|^4}   \langle \nabla_R g, \partial^i B \nabla_R g \rangle ^2 .
\end{align*}

Note that \boxed{E1} and \boxed{E2} can be considered ``error terms'' from a perturbation point of view, since in case $B$ were constant we would get $\boxed{E1} = \boxed{E2} = 0$, while \boxed{M1} can be considered a ``main term''.

Next we perform a similar decomposition for $\partial^i (b_i \partial^i |S\nabla g|^2)$, in the second term in the right-hand side in \rf{term 3A}. In this case, we need to pay attention to the first term in the right-hand side of \rf{derivative product of dot product}:
\begin{equation}\label{using 2 derivative product of dot product}
		\partial^i \left( b_i \langle \partial^i \nabla_R g, B \nabla_R g \rangle \right)
		=  \left\langle \partial^i  \left(b_i \partial^i \nabla_R g \right),B \nabla_R g \right\rangle 
		+ b_i \langle \partial^i \nabla_R g, \partial^i B \nabla_R g \rangle
		+ b_i \langle \partial^i \nabla_R g, B \partial^i \nabla_R g \rangle .
\end{equation}

Then using \rf{derivative product of dot product} and \rf{using 2 derivative product of dot product} we get that the second term in the right-hand side of \rf{term 3A} can be written as
$$
\sum_{i=1}^2 \frac{1}{|S\nabla g|^2} \partial^i (b_i \partial^i |S\nabla g|^2 ) 
= \boxed{M2} + T_{M3} + \boxed{E3} + T_{E4},
$$
where
\begin{align*}
	\boxed{M2} & \coloneqq \sum_{i=1}^2 \frac{2 b_i}{|S\nabla g|^2} \langle \partial^i \nabla_R g, B \partial^i \nabla_R g \rangle , \\
	T_{M3} & \coloneqq \frac{2}{|S\nabla g|^2} \left\langle \sum_{i=1}^2 \partial^i  \left(b_i \partial^i \nabla_R g \right), B \nabla_R g \right\rangle , \\
	\boxed{E3} & \coloneqq \sum_{i=1}^2 \frac{2 b_i}{|S\nabla g|^2}  \langle \partial^i \nabla_R g, \partial^i B \nabla_R g \rangle , \text{ and}\\
	T_{E4} & \coloneqq  \sum_{i=1}^2 \frac{1}{|S\nabla g|^2} \partial^i \left( b_i \langle \nabla_R g, \partial^i B \nabla_R g \rangle \right) .
\end{align*}

All together we have
\begin{align*}
	\divv \left( A^T \nabla \log |S\nabla g|^2 \right) =&\ T_{3.A} + \boxed{3.B} + \boxed{3.C} & \\
	=&\ \boxed{M1} + \boxed{M2} + T_{M3} & \text{(Main terms)}\\
	&+ \boxed{3.B} + \boxed{3.C} + \boxed{E1} +\boxed{E2} + \boxed{E3} + T_{E4}, & \text{(Error terms)}
\end{align*}
and the first part of the lemma follows taking
$$
T_E \coloneqq \boxed{3.B} + \boxed{3.C} + \boxed{E1} +\boxed{E2} + \boxed{E3} + T_{E4}.
$$

Let us see now the second part of the lemma. Indeed, $\boxed{E2}$ is $\OO(1)$ because $|\langle \nabla_R g, \partial^i B \nabla_R g \rangle| \lesssim |\nabla g|^2 $. By \rf{gradient log_a gradient g} we get that $\boxed{3.B}$ and $\boxed{3.C}$ are in $\OO\left(1 + \frac{|\nabla^2 g|}{|\nabla g|}\right)$, and \rf{two directional derivatives} implies that $\boxed{E1}$ and $\boxed{E3}$ are of the form $\OO\left( 1 + \frac{|\nabla^2 g|}{|\nabla g|}\right)$. All in all,
$$
T_E = T_{E4} + \OO\left( 1+\frac{|\nabla^2 g|}{|\nabla g|} \right) ,
$$
as claimed.\qed

\subsubsection{Proof of \cref{cancellation main terms M1 M2}}\label{proof of cancellation main terms M1 M2}%lemma 3

In this subsection we prove that the two main terms cancel out in the sense 
$$
\left| \boxed{M1} + \boxed{M2} \right| \lesssim 1 + \frac{|\nabla^2 g|}{|\nabla g|}.
$$
This cancellation is related to what happens in the constant case, see \cref{constant matrix and main terms}.

We will see that they both have a common term in opposite sign, and exploiting this cancellation we will obtain a sum of error terms. The key identity is \rf{property 2 green function}, which only applies in the plane.

The key idea to prove the lemma is that using \rf{derivative_commutator}, which relates $\partial^{1,2} g$ and $\partial^{2,1} g$, and \rf{property 2 green function}, which relates $b_1 \partial^{1,1} g$ and $b_2 \partial^{2,2} g$, we will be able to write $\partial^{1,1}$, $\partial^{2,2}$, $\partial^{1,2}$ and $\partial^{2,1}$ in terms of $\partial^{1,1}$ and $\partial^{2,1}$ only.

We start by studying the term $\boxed{M1} \coloneqq -\sum_i \frac{4 b_i}{|S\nabla g|^4}  \langle \partial^i \nabla_R g, B\nabla_R g \rangle^2$ in $\widetilde\Omega_\varphi$. Expanding the numerator,
\begin{align*}
	\sum_{i=1}^2  b_i \langle \partial^i \nabla_R g, B \nabla_R g \rangle^2 &= b_1 \langle \partial^1 \nabla_R g, B \nabla_R g \rangle^2 + b_2 \langle \partial^2 \nabla_R g, B \nabla_R g \rangle^2 \\
	&=b_1 (b_1 \partial^1 g \partial^{1,1} g + b_2 \partial^2 g \partial^{1,2} g)^2 + b_2 (b_1 \partial^1 g \partial^{2,1} g + b_2 \partial^2 g \partial^{2,2} g)^2.
\end{align*}
By \rf{derivative_commutator} and \rf{property 2 green function}, which read as $\partial^{1,2} g = \partial^{2,1} g + \OO(|\nabla g|)$ and $b_2 \partial^{2,2} g = - b_1 \partial^{1,1} g + \OO(|\nabla g|)$ respectively, we have 
\begin{align*}
	\sum_{i=1}^2  b_i \langle \partial^i \nabla_R g, B \nabla_R g \rangle^2 
	=&\ b_1 (b_1 \partial^1 g \partial^{1,1} g + b_2 \partial^2 g\partial^{2,1} g+\OO(|\nabla g|^2))^2 \\
	&+ b_2 (b_1 \partial^1 g \partial^{2,1} g - b_1 \partial^2 g\partial^{1,1} g+\OO(|\nabla g|^2))^2.
\end{align*}
Expanding, and using $b_1 (\partial^1 g)^2 + b_2 (\partial^2 g)^2 =\langle \nabla_R g, B \nabla_R g \rangle = |S\nabla g|^2$ and the cancellation of cross terms, we get
\begin{equation*}
	\sum_{i=1}^2  b_i \langle \partial^i \nabla_R g, B \nabla_R g \rangle^2 
	=|S\nabla g|^2 \left( b_1^2 (\partial^{1,1} g)^2 +  b_1 b_2 (\partial^{2,1} g)^2 \right)
	+ (\partial^{1,1} g+\partial^{2,1} g) \OO(|\nabla g|^3)
	+\OO (|\nabla g|^4).
\end{equation*}
In conclusion, \boxed{M1} is of the form
\begin{equation}\label{M1 simplified}
		\boxed{M1} 
		= -4\frac{ b_1^2 (\partial^{1,1} g)^2 +  b_1 b_2 (\partial^{2,1} g)^2  }{|S\nabla g|^2} + \frac{ \OO (\partial^{1,1} g + \partial^{2,1} g )  }{|S\nabla g|} + \OO(1).
\end{equation}

With the same strategy we study the term $\boxed{M2} \coloneqq \sum_i \frac{2 b_i}{|S\nabla g|^2} \langle \partial^i \nabla_R g, B \partial^i \nabla_R g \rangle$ in $\widetilde\Omega_\varphi$. Using \rf{derivative_commutator} and \rf{property 2 green function} as before and expanding, the numerator can be written as
\begin{align*}
	\sum_{i=1}^2 & b_i \langle  \partial^i \nabla_R g , B \partial^i \nabla_R g \rangle = b_1^2 (\partial^{1,1} g)^2 + b_1 b_2 (\partial^{1,2} g)^2 + b_1 b_2 (\partial^{2,1} g)^2 + b_2^2 (\partial^{2,2} g)^2 \\
	&=  b_1^2 (\partial^{1,1} g)^2 + b_1 b_2 \left( \partial^{2,1} g + \OO(|\nabla g|)\right)^2 + b_1 b_2 (\partial^{2,1} g)^2 + \left(-b_1 \partial^{1,1} g + \OO(|\nabla g|)\right)^2 \\
	&= 2\left( b_1^2 (\partial^{1,1} g)^2 + b_1 b_2 (\partial^{2,1}g)^2 \right)
	+ \left(\partial^{1,1} g + \partial^{2,1} g\right) \OO (|\nabla g|)
	+\OO (|\nabla g|^2).
\end{align*}
Hence, \boxed{M2} is of the form
\begin{equation}\label{M2 simplified}
		\boxed{M2}=4\frac{b_1^2 (\partial^{1,1} g)^2 + b_1 b_2 (\partial^{2,1}g)^2}{|S\nabla g|^2}+ \frac{\OO (\partial^{1,1} g + \partial^{2,1} g)}{|S\nabla g|} + \OO(1).
\end{equation}

Note that both \boxed{M1} in \rf{M1 simplified} and \boxed{M2} in \rf{M2 simplified} have the term
$$
4\frac{b_1^2 (\partial^{1,1} g)^2 + b_1 b_2 (\partial^{2,1}g)^2}{|S\nabla g|^2} 
$$
in common with opposite sign, which allows us to remove the square in \rf{square dependence in two main terms}. That is, adding up both terms,
$$
\boxed{M1} + \boxed{M2}  =\frac{\OO(\partial^{1,1} g + \partial^{2,1} g)}{|S\nabla g|} + \OO(1),
$$
whence we obtain
$$
\left| \boxed{M1} + \boxed{M2} \right| \overset{\text{\rf{two directional derivatives}}}{\lesssim} 1 + \frac{|\nabla^2 g|}{|\nabla g|},
$$
as claimed.
\qed

\subsubsection{Proof of \cref{cancellation main term M3}}\label{proof of cancellation main term M3}

We study, for $\psi\in W^{1,\infty}_c (\widetilde\Omega_\varphi)$, the term
$$
T_{M3} (\psi) \coloneqq \left(\frac{2}{|S\nabla g|^2} \left\langle \sum_{i=1}^2 \partial^i  \left(b_i \partial^i \nabla_R g \right), B \nabla_R g \right\rangle \right)(\psi) .
$$

In this subsection we don't use the fact that we are in the plane. Instead, the key ingredient is the $L_A$-harmonicity of the Green function far from the pole, and so the computations in this subsection could be done also in higher dimensions.

The term $T_{M3}$ must be read as
\begin{equation}\label{decomposition term M3 distribution}
	\begin{aligned}
	T_{M3} &= \frac{2}{|S\nabla g|^2} \left\langle
	\sum_{i=1}^2
	\begin{pmatrix}
		\partial^i (b_i\partial^{i,1} g ) \\
		\partial^i (b_i\partial^{i,2} g )
	\end{pmatrix} ,
	\begin{pmatrix}
		b_1 \partial^1 g \\
		b_2 \partial^2 g
	\end{pmatrix}
	\right\rangle \\
	&= \frac{2}{|S\nabla g|^2} \sum_{i=1}^2  \left\{
		\partial^i (b_i\partial^{i,1} g ) b_1 \partial^1 g
		+ \partial^i (b_i\partial^{i,2} g ) b_2 \partial^2 g
	\right\} \\
	&= \frac{2}{|S\nabla g|^2} \sum_{i=1}^2 \sum_{j=1}^2
	\partial^i (b_i\partial^{i,j} g ) b_j \partial^j g 
	\eqqcolon \sum_{i=1}^2 \sum_{j=1}^2 T_{M3}^{(i,j)} .
	\end{aligned}
\end{equation}
Here, for each $i,j \in \{1,2\}$,
\begin{equation}\label{distribution M3-ij}
T_{M3}^{(i,j)} (\psi) = T_{\partial^i (b_i \partial^{i,j} g )} \left( \frac{2}{|S\nabla g|^2} b_j \partial^j g \psi \right) .
\end{equation}

We want to study the functionals $T_{M3}^{(i,j)}$ for $i,j\in \{1,2\}$. First note that
$$
\frac{2}{|S\nabla g|^2} b_j \partial^j g \psi \in W_c^{1,\infty}(\widetilde\Omega_\varphi),
$$
while $b_i \partial^{i,j} g\in L^{2} (\widetilde\Omega_\varphi)$. Hence, \rf{distribution M3-ij} makes sense and it suffices to study the functionals $T_{\partial^i (b_i \partial^{i,j} g )}\in W_c^{1,\infty}(\widetilde\Omega)^\prime$. In fact, fixed $j\in \{1,2\}$, we will exploit the cancellation of the functional $\sum_{i=1}^2 T_{\partial^i (b_i \partial^{i,j} g )}$. For simplicity we will write $\partial^i (b_i \partial^{i,j} g )$ instead of $T_{\partial^i (b_i \partial^{i,j} g)}$. 

Compared to the strategy seen in \cref{constant matrix and main terms}, now the matrix is not constant and hence the directional derivatives do not commute. For this reason, some error terms will appear in the procedure of extracting the gradient $\nabla_R$ outside from
\begin{equation}\label{gradient R is inside}
	\sum_{i=1}^2 \partial^i  \left(b_i \partial^i \nabla_R g \right) = 
	\begin{pmatrix}
		\partial^1 (b_1 \partial^{1, 1} g ) + \partial^2 (b_2 \partial^{2, 1} g ) \\
		\partial^1 (b_1 \partial^{1, 2} g ) 
		+ \partial^2 (b_2 \partial^{2, 2} g )
	\end{pmatrix} ,
\end{equation}
as it is done in the constant matrix case in \rf{extract gradient R outside}. The idea is the same as in the proof of \cref{cancellation main terms M1 M2}, see \cref{proof of cancellation main terms M1 M2}. That is, using the `almost'-commutative property \rf{derivative_commutator} (relating $\partial^{1,2} g$ and $\partial^{2,1} g$) of the directional derivatives, and its functional version in \rf{directional derivatives does not commute distribution}, we manage to extract the gradient $\nabla_R$ outside.

Let us fix $j\in \{1,2\}$. If $i=j$ (clearly $\partial^{i,j} = \partial^{j,i}$ in this case), then
\begin{equation}\label{gradient R outside step 1}
	\partial^i (b_i \partial^{i,j} g) 
	= \partial^j \left(\partial^i (b_i \partial^i g) \right) - \partial^i ( \partial^j b_i \partial^i g ) ,
\end{equation}
and when $i\not= j$, using \rf{derivative_commutator} and \rf{directional derivatives does not commute distribution} we get
\begin{equation}\label{gradient R outside step 2}
	\begin{aligned}
		\partial^i (b_i \partial^{i, j} g) \overset{\text{\rf{derivative_commutator}}}&{=} \partial^i \left(b_i \partial^{j,i} g  \right) 
		+ \partial^i \left(b_i (\partial^i R_j - \partial^j R_i)\nabla g\right) \\
		&=\partial^i \left(\partial^j (b_i  \partial^i g )\right) 
		- \partial^i \left( \partial^j b_i \partial^i g \right)
		+ \partial^i \left(b_i (\partial^i R_j - \partial^j R_i)\nabla g\right) \\
		\overset{\text{\rf{directional derivatives does not commute distribution}}}&{=} \partial^j \left(\partial^i  (b_i  \partial^i g ) \right)
		- (\partial^j R_i - \partial^i R_j ) \nabla (b_i \partial^i g)
		- \partial^i \left( \partial^j b_i \partial^i g \right)
		+ \partial^i \left(b_i (\partial^i R_j - \partial^j R_i)\nabla g\right).
	\end{aligned}
\end{equation}

Plugging \rf{gradient R outside step 1} and \rf{gradient R outside step 2} in \rf{gradient R is inside}, we get, for $j=1$,
\begin{equation}\label{gradient inside to outside j=1}
	\begin{aligned}
	\sum_{i=1}^2 \partial^i (b_i \partial^{i,1} g )
	=&\ \partial^1 \left( \sum_{i=1}^2 \partial^i  \left(b_i \partial^i g \right) \right)
	- \sum_{i=1}^2\partial^i ( \partial^1 b_i \partial^i g ) 
	- (\partial^1 R_2 - \partial^2 R_1 ) \cdot \nabla (b_2 \partial^2 g)\\
	& + \partial^2 \left(b_2 (\partial^2 R_1 - \partial^1 R_2)\cdot\nabla g\right),
	\end{aligned}
\end{equation}
and for $j=2$,
\begin{equation}\label{gradient inside to outside j=2}
	\begin{aligned}
	\sum_{i=1}^2 \partial^i (b_i \partial^{i,2} g )
	=&\ \partial^2 \left( \sum_{i=1}^2 \partial^i  \left(b_i \partial^i g \right) \right)
	- \sum_{i=1}^2\partial^i ( \partial^2 b_i \partial^i g ) 
	- (\partial^2 R_1 - \partial^1 R_2 ) \cdot \nabla (b_1 \partial^1 g)\\
	& + \partial^1 \left(b_1 (\partial^1 R_2 - \partial^2 R_1)\cdot\nabla g\right).
	\end{aligned}
\end{equation}

Note that in both cases there is the term $\sum_{i=1}^2 \partial^i  \left(b_i \partial^i g \right)$, which is studied in \cref{green function property grad g}. By \rf{decomposition term M3 distribution} and using \cref{green function property grad g} in both \rf{gradient inside to outside j=1} and \rf{gradient inside to outside j=2}, in the dual space we get
\begin{equation}\label{final form after cancellation of M3}
	\begin{aligned}
	T_{M3} =\ & \frac{2}{|S\nabla g|^2} \sum_{j=1}^2 b_j \partial^j g \left(\sum_{i=1}^2 \partial^i (b_i \partial^{i,j} g)\right) \\
	=\ & \sum_{j=1}^2 \frac{2b_j \partial^j g}{|S\nabla g|^2} \left\{
	\partial^j \left( 
	-\sum_{i=1}^{2} b_i \partial^i g  \divv R_i^T 
	- \sum_{i,k=1}^{2} \partial_i \left(\frac{a_{ik}-a_{ki}}{2}\right) \partial_k g
	\right)
	- \sum_{i=1}^2\partial^i ( \partial^j b_i \partial^i g )\right\}\\
	&+ \frac{2b_1 \partial^1 g}{|S\nabla g|^2} \left\{ - (\partial^1 R_2 - \partial^2 R_1 ) \cdot \nabla (b_2 \partial^2 g)
	+ \partial^2 \left(b_2 (\partial^2 R_1 - \partial^1 R_2)\cdot\nabla g\right) \right\} \\
	&+ \frac{2b_2 \partial^2 g}{|S\nabla g|^2} 
	\left\{ - (\partial^2 R_1 - \partial^1 R_2 ) \cdot \nabla (b_1 \partial^1 g)
	+ \partial^1 \left(b_1 (\partial^1 R_2 - \partial^2 R_1)\cdot\nabla g\right) \right\}.
	\end{aligned}
\end{equation}

As a consequence of this, the right-hand equality in \cref{green function property grad g} and $|\nabla \partial^i g| = \OO (|\nabla g|) + \OO (|\nabla^2 g|)$, see \rf{two directional derivatives}, we conclude that the functional $T_{M3}$ is of the form
\begin{equation*}
	T_{M3} = \sum_{j=1}^2 \frac{2 b_j \partial^j g}{|S\nabla g|^2} \left(
	\OO (|\nabla g|) + \OO (|\nabla^2 g|)
	+\partial^1 \left( \OO (|\nabla g|)\right) 
	+\partial^2 \left( \OO (|\nabla g|) \right)
	\right),
\end{equation*}
as claimed.
\qed

\subsubsection{Proof of \cref{bound error and main M3 terms}}\label{proof of bound error and main M3 terms}%lemma 4

First, we need the following claim, which we will prove later.

\begin{claim}\label{equalities and bounds to apply density argument}
	Let $\psi \in W^{1,\infty}_c (\widetilde\Omega_\varphi)$. Then, for $i,j\in \{1,2\}$, both
	$$
	T_{\frac{1}{|S\nabla g|^2} \partial^i \left( \OO (|\nabla g|^2) \right)} (\psi) 
	\text{ and }
	T_{\frac{b_j \partial^j g}{|S\nabla g|^2} \partial^i \left( \OO (|\nabla g|) \right)} (\psi)
	$$
	are of the form
	\begin{equation*}
		\int \frac{\OO (|\nabla^2 g|)}{|\nabla g|} \psi
		+ \int \OO (1) \psi
		+ \int \OO (1) \partial^i \psi .
	\end{equation*}
\end{claim}

Granted this, by \cref{cancellation main term M3} and \rf{T_E simplified} we have
\begin{equation*}
	|T_{M3} (\psi_\varepsilon\varphi g)| + |T_E (\psi_\varepsilon\varphi g)|
	\lesssim \int \frac{|\nabla^2 g|}{|\nabla g|} \psi_\varepsilon\varphi g
	+ \int \psi_\varepsilon\varphi g
	+ \int |\partial^i (\psi_\varepsilon\varphi g)| .
\end{equation*}
The second term in the right-hand side is bounded by a constant times $\HH^2(\Omega)<\infty$ as $\Omega$ is bounded and $g\lesssim 1$ in $\supp\varphi$, see \rf{green bounded by measure and 1 in the supp of varphi}. Hence, it suffices to prove
$$
\int \left|\partial^i (\psi_\varepsilon\varphi g) \right| \lesssim 1, \text{ for }i\in \{1,2\}.
$$

By the product derivative rule and the notation $\partial^i f = R_i \cdot \nabla f $, we can use
\begin{equation*}
	\int \left| \partial^i \left( \psi_\varepsilon \varphi g\right) \right|
	\lesssim
	\int g|\nabla\varphi|
	+\int_{\supp\varphi} |\nabla g|
	+ \frac{1}{\varepsilon} \int_{\supp \nabla \psi_\varepsilon} g.
\end{equation*}
Since $|\nabla \varphi| \lesssim 1$ and $g(x)\lesssim 1$ in $\supp \varphi$ (see \rf{green bounded by measure and 1 in the supp of varphi}), we get $\int g |\nabla \varphi| \lesssim 1$. For the second integral on the right-hand side we sum over Whitney cubes in $W_\varphi$ and apply \cref{sum over whitney region p3,sum over whitney region p1} in \cref{sum over whitney region}
\begin{equation*}
	\int_{\supp \varphi} |\nabla g| 
	\leq \sum_{Q \in W_\varphi} \int_Q |\nabla g| 
	\overset{\text{\rf{sum over whitney region p3}}}{\lesssim} \sum_{Q \in W_\varphi} \ell(Q) \cdot \widetilde \omega_T (50Q) 
	\overset{\text{\rf{sum over whitney region p1}}}{\lesssim} 1.
\end{equation*}
On the other hand, using $\supp \nabla \psi_\varepsilon \subset U_{3\varepsilon} (\partial \widetilde\Omega)$ and \cref{integrals involving green function near boundary}\rf{integral g near boundary} respectively, we have
$$
\frac{1}{\varepsilon} \int_{\supp \nabla \psi_\varepsilon} g
\leq \frac{1}{\varepsilon} \int_{U_{3\varepsilon} (\partial \widetilde\Omega)} g \lesssim \varepsilon \leq 1,
$$
and \cref{bound error and main M3 terms} follows.\qed

We now turn to the proof of \cref{equalities and bounds to apply density argument}.
\begin{proof}[Proof of \cref{equalities and bounds to apply density argument}]
        First we study $T_{\frac{1}{|S\nabla g|^2} \partial^i \left( \OO (|\nabla g|^2) \right)} (\psi)$. We have
	\begin{align*}
		T_{\frac{1}{|S\nabla g|^2} \partial^i (\OO (|\nabla g|^2))} (\psi) &=
		T_{\partial^i (\OO (|\nabla g|^2))} \left( \frac{1}{|S\nabla g|^2} \psi \right) \\
		\overset{\text{\rf{def2:distributional directional derivative}}}&{=}
		-\int \OO (|\nabla g|^2) \frac{1}{|S\nabla g|^2} \psi \divv R_i^T 
		- \int \OO (|\nabla g|^2) \partial^i \left( \frac{1}{|S\nabla g|^2} \psi \right).
	\end{align*}
	The first term in the right-hand side is of the form $\int \OO (1) \psi$. Let us study the second term in the right-hand side. By \rf{derivative quotient},
	\begin{equation*}
		\partial^i \left( \frac{1}{|S\nabla g|^2} \psi \right)
		= - \frac{2\langle \partial^i \nabla_R g, B \nabla_R g \rangle}{|S\nabla g|^4} \psi
		- \frac{\langle \nabla_R g, \partial^i B \nabla_R g \rangle}{|S\nabla g|^4} \psi
		+ \frac{\partial^i \psi}{|S\nabla g|^2} ,
	\end{equation*}
	and hence we get
	\begin{align*}
		- \int \OO (|\nabla g|^2) \partial^i \left( \frac{1}{|S\nabla g|^2} \psi \right) 
		&= \int \frac{\OO (\partial^i \nabla_R g)}{|S\nabla g|} \psi
		+ \int \OO (1) \psi
		+ \int \OO (1) \partial^i \psi \\
		\overset{\text{\rf{two directional derivatives}}}&{=}
		\int \frac{\OO (|\nabla^2 g|)}{|\nabla g|} \psi
		+ \int \OO (1) \psi
		+ \int \OO (1) \partial^i \psi .
	\end{align*}
	
	Next we study $T_{\frac{b_j \partial^j g}{|S\nabla g|^2} \partial^i \left( \OO (|\nabla g|) \right)} (\psi)$, similarly as we did with the previous term. Now we have
	\begin{align*}
		T_{\frac{b_j\partial^j g}{|S\nabla g|^2} \partial^i (\OO(|\nabla g|))} (\psi)
		&=T_{\partial^i (\OO(|\nabla g|))} \left(\frac{b_j \partial^j g}{|S\nabla g|^2} \psi\right) \\
		\overset{\text{\rf{def2:distributional directional derivative}}}&{=} 
		-\int \frac{\OO (|\nabla g|) b_j \partial^j g}{|S\nabla g|^2} \psi \divv R_i 
		-\int \OO(|\nabla g|) \partial^i \left( \frac{b_j \partial^j g}{|S\nabla g|^2} \psi \right) .
	\end{align*}
	As before, the first term in the right-hand side is of the form $\int \OO(1) \psi$ and the second term is
	\begin{equation*}
		-\int \OO(|\nabla g|) \partial^i \left( \frac{b_j \partial^j g}{|S\nabla g|^2} \psi \right) 
		=
		\int \frac{\OO(1) \psi}{|S\nabla g|} \partial^i \left( b_j \partial^j g \right) 
		+\int \OO(|\nabla g|^2) \partial^i \left( \frac{1}{|S\nabla g|^2} \right) \psi
		+ \int \OO(1) \partial^i \psi.
	\end{equation*}
	The claim follows by using \rf{two directional derivatives} and \rf{derivative quotient} in the previous line.
\end{proof}

\vv

\renewcommand{\abstractname}{Acknowledgements}
\begin{abstract}
Part of this work was carried out while I.G.M. and X.T. were visiting the Hausdorff Institute for Mathematics in Bonn (Germany) during the trimester program \emph{Interactions between Geometric measure theory, Singular integrals, and PDE}.    

I.G.M. was supported by the PIF UAB 2020/21 from the Universitat Autònoma de Barcelona, and by Generalitat de Catalunya’s Agency for Management of University and Research Grants (AGAUR) (2021 FI\_B 00637 and 2023 FI-3 00151). I.G.M. and M.P. were partially supported by the Spanish State Research Agency (AEI) project PID2021-125021NAI00 (Spain). M.P. was partially supported by the projects PID2021-123405NB-I00 (AEI), and 2021-SGR-00087 (AGAUR). M.P. and X.T. were supported by the AEI, through the Severo Ochoa and Mar\'{i}a de Maeztu Program for Centers and Units of Excellence in R\&D (CEX2020-001084-M). X.T. was supported by the European Research Council (ERC) under the European Union’s Horizon 2020 research and innovation programme (grant agreement 101018680), and project 2021-SGR-00071 (AGAUR). All the authors were also supported by PID2020-114167GB-I00. We thank CERCA Programme/Generalitat de Catalunya for institutional support.
\end{abstract}

% ********************************************************************************
% ********************************************************************************

\vv

\vv

%\printbibliography
\bibliographystyle{alpha}
%\bibliography{./../../References-bibtex-JabRef/references-phd.bib} %texstudio
\bibliography{references} %overleaf

\end{document}